\documentclass[11pt, pdftex]{amsart}

\usepackage{amsmath}
\usepackage{amscd}
\usepackage{amssymb}
\usepackage{amsbsy}
\usepackage{amsfonts}
\usepackage{color}
\usepackage{latexsym}
\usepackage{graphics}
\usepackage{graphicx}
\usepackage{amsmath,amscd,latexsym}
\usepackage{multirow}
\usepackage{array}
\usepackage{paralist}
\usepackage{titletoc}
\usepackage{url}
\usepackage{epsfig, eucal,float}
\usepackage{CJKutf8}

\graphicspath{{Pic_Thurston/}}

\theoremstyle{plain}
\newtheorem{theorem}{Theorem}[section]
\newtheorem{proposition}[theorem]{Proposition}
\newtheorem{lemma}[theorem]{Lemma}
\newtheorem{corollary}[theorem]{Corollary}
\newtheorem{remark}[theorem]{Remark}
\newtheorem{definition}[theorem]{Definition}

\newtheorem{example}[theorem]{Example}
\newtheorem{main theorem}[theorem]{Main Theorem}
\newtheorem{problem}[theorem]{Problem}
\newtheorem{question}[theorem]{Question}

\newtheorem{conjecture}[theorem]{Conjecture}

\newlength\savewidth

\newcommand{\interior}{\operatorname{int}}

\newcommand{\vol}{\operatorname{vol}}
\newcommand{\Comm}{\operatorname{Comm}}

\newcommand{\lk}{\operatorname{lk}}

\newcommand{\NN}{\mathbb{N}}
\newcommand{\ZZ}{\mathbb{Z}}
\newcommand{\QQ}{\mathbb{Q}}
\newcommand{\RR}{\mathbb{R}}
\newcommand{\CC}{\mathbb{C}}
\newcommand{\HH}{\mathbb{H}}
\newcommand{\EE}{\mathbb{E}}
\newcommand{\Sphere}{\mathbb{S}}

\newcommand{\zb}{\mbox{\boldmath$z$}}
\newcommand{\ub}{\mbox{\boldmath$u$}}
\newcommand{\vb}{\mbox{\boldmath$v$}}
\newcommand{\taub}{\mbox{\boldmath$\tau$}}
\newcommand{\nub}{\mbox{\boldmath$\nu$}}
\newcommand{\inftyb}{\mbox{\boldmath$\infty$}}

\newcommand{\PSL}{\mbox{$\mathrm{PSL}$}}
\newcommand{\SL}{\mbox{$\mathrm{SL}$}}
\newcommand{\GL}{\mbox{$\mathrm{GL}$}}
\newcommand{\PGL}{\mbox{$\mathrm{PGL}$}}
\newcommand{\DD}{\mathcal{D}}
\newcommand{\foliation}{\mathcal{F}}
\newcommand{\Subgp}{\mathcal{C}}

\newcommand{\OO}{{\mathcal{O}}}

\newcommand{\Ford}{{\mathcal{F}}}

\newcommand{\HypTorP}{{\mathcal{J}}}
\newcommand{\QF}{{\mathcal{QF}}}
\newcommand{\LL}{{\mathcal{L}}}
\newcommand{\XX}{{\mathcal{X}}}
\newcommand{\Rep}{{\mathcal{R}}}
\newcommand{\UniversalOrb}{{\mathcal{U}}}
\newcommand{\UniversalGp}{{\mathrm{U}}}

\newcommand{\fix}{\operatorname{Fix}}

\newcommand{\Orth}{\operatorname{O}}
\newcommand{\Isom}{\operatorname{Isom}}
\newcommand{\Aut}{\operatorname{Aut}}
\newcommand{\Ker}{\operatorname{Ker}}

\newcommand{\Sym}{\operatorname{Sym}}
\newcommand{\Out}{\operatorname{Out}}
\newcommand{\Diff}{\operatorname{Diff}}
\newcommand{\Hom}{\operatorname{Hom}}
\newcommand{\tr}{\operatorname{tr}}
\newcommand{\inj}{\operatorname{inj}}
\newcommand{\Teich}{\operatorname{Teich}}

\newcommand{\MCG}{\operatorname{MCG}}
\newcommand{\PML}{\operatorname{PML}}
\newcommand{\PMF}{\operatorname{PMF}}
\newcommand{\Tor}{\operatorname{Tor}}
\newcommand{\axis}{\operatorname{axis}}

\newcommand{\CT}{\mbox{$\kappa$}}

\makeatletter
\renewcommand\subsection{\@startsection{subsection}{2}{0mm}
    {-10.5dd plus-8pt minus-4pt}{10.5dd}
     {\normalsize\upshape}}
\makeatother

\begin{document}

\title{A survey of the impact of Thurston's work on Knot Theory}

\author{Makoto Sakuma}
\address{Osaka City University Advanced Mathematical Institute\\
3-3-138 Sugimoto, Sumiyoshi-ku, Osaka, 558-8585, Japan}
\address{Department of Mathematics\\
Hiroshima University\\
Higashi-Hiroshima, 739-8526, Japan}
\email{sakuma@hiroshima-u.ac.jp}

\date{\today}


\begin{abstract}
This is a survey of the impact of 
Thurston's work on knot theory,
laying emphasis on the two characteristic features,
rigidity and flexibility, of $3$-dimensional hyperbolic structures.
We also lay emphasis on the role of the classical invariants,
the Alexander polynomial and the homology
of finite branched/unbranched coverings.
\end{abstract}

\maketitle

\setcounter{tocdepth}{3}
\tableofcontents

\section{Introduction}

Knot theory is the analysis of pairs $(S^3,K)$,
where $K$ is a knot (i.e., an embedded circle) in the $3$-sphere $S^3$,
and classification of knots
has been one of the main problems in knot theory.
The Alexander polynomial is an excellent invariant of knots,
and it had been a dominating tool and theme in knot theory, 
until knot theory was influenced by Thurston's work and 
the Jones polynomial was discovered.
In fact, the classical textbook by Crowell and Fox \cite{Crowell-Fox} is devoted to 
the calculation of the knot group and
the definition of the Alexander polynomial by using the free differential calculus.
The influential textbook by Rolfsen \cite{Rolfsen} lies emphasis on 
geometric understanding of the Alexander polynomial
through surgery description of the infinite cyclic cover
(cf. \cite{Hironaka2016}).
However, the Alexander polynomial is far from being complete:
there are infinitely many nontrivial knots with trivial Alexander polynomial.
The famous Kinoshita--Terasaka knot and the Conway knot 
are related by mutation, and therefore 
no skein polynomial, including the Alexander polynomial, 
can distinguish between them.
The first proof of their inequivalence was given by 
Riley \cite{Riley1971}
by studying parabolic representations of the knot groups into
the finite simple group $\PSL(2,\ZZ/7\ZZ)$
and the homology of corresponding finite branched/unbranched coverings.
(This work led him to the discovery of the hyperbolic structure 
of the figure-eight knot complement,
which inspired Thurston.)
Riley called this a {\it universal method} for obtaining algebraic invariants of knots.
The method turned out to be, 
at least experimentally, a very powerful tool in knot theory,
due to the development of computer technology. 
However
theoretical background of the universal method has not been given yet.

In 1976, around the time Rolfsen's book was published,
William Thurston started a series of lectures
on \lq\lq The geometry and toplogy of $3$-manifolds''.
His lecture notes \cite{Thurston0} begin with the following words.

{\it
The theme I intend to develop is that topology and geometry,
in dimensions up through $3$,
are intricately related.
Because of this relation, many questions which seem utterly hopeless 
from a purely topological point of view can be fruitfully studied.
It is not totally unreasonable to hope that
eventually all $3$-manifolds will be understood in a systematic way.
}

This prophecy turned out to be true.
Thurston's work has revolutionized $3$-dimensional topology,
and it has had tremendous impact on knot theory.
The first major impact was 
the proof of the Smith conjecture \cite{Morgan-Bass},
a result of the efforts by
Thurston, Meeks and Yau, Bass, Shalen, Gordon and Litherland, and Morgan.
As Morgan predicted in \cite[p.6]{Morgan-Bass},
this was just the beginning of the saga. 

\medskip

In this chapter, we give a survey of the impact of Thurston's work
on knot theory.
However, the impact is huge, whereas 
both my ability and knowledge are poor.
Moreover, there already exist excellent surveys,
including 
Callahan-Reid \cite{CR},
Adams \cite{Adams1} and
Futer-Kalfagianni-Purcell \cite{FKP_survey}.
So, I decided to lay emphasis on the two characteristic features,
{\it rigidity} and {\it flexibility}, of hyperbolic $3$-manifolds.

As the title of Section 5 of
Thurston's lecture notes \cite{Thurston0} represents,
hyperbolic structures on $3$-manifolds
have two different features, rigidity and flexibility.

The Mostow-Prasad rigidity theorem says that a complete hyperbolic structure of finite volume
on an $n$-manifold with $n\ge 3$ is rigid: 
it does not admit local deformation, and moreover,
such a structure is unique.
Thus any geometric invariant determined by the complete hyperbolic structure
of an $n$-manifold $M$ with $n\ge 3$
is automatically a topological invariant of $M$.
Thurston's uniformization theorem for Haken manifolds
implies that almost every knot $K$ in $S^3$ is hyperbolic,
namely the complement $S^3-K$ admits a complete hyperbolic structure
of finite volume.
Thus we obtain plenty of topological invariants of hyperbolic knots,
including the volume, the maximal cusp volume, 
the Euclidean modulus of the cusp torus,
the length spectrum, 
the lengths of geodesic paths joining the cusp to itself,
the invariant trace field, the invariant quaternion algebra, etc.
In particular, 
the {\it canonical decomposition} (see Subsection \ref{subsec:canonical})
gives 
a complete combinatorial invariant
for hyperbolic knots, by virtue of
the Gordon-Luecke knot complement theorem.
The computer program, SnapPea, developed by J. Weeks
enables us to calculate the canonical decomposition of hyperbolic knot complements.
For example, we can easily detect the inequivalence of 
the Kinoshita--Terasaka knot
and the Conway knot,
by checking with SnapPea that the number of $3$-cells 
in the canonical decompositions of the knot complements
are $12$ and $14$, respectively.
The rigidity theorem provides us a number of 
powerful invariants, and
it has enriched knot theory by opening new directions of research,
namely the study of the behavior of the geometric invariants.
An enormous amount of deep research have been made in these new directions.
(See Sections \ref{sec:computation-canonical}, \ref{sec:volume} and
\ref{sec:arithmetic-invariant}.)

There are two kinds of flexibility of hyperbolic structures on $3$-manifolds.
One of them is that of cusped hyperbolic manifolds:
the complete hyperbolic structure
admits nontrivial continuous deformations into incomplete hyperbolic structures.
By considering the metric completions of incomplete hyperbolic structures,
Thurston established the hyperbolic Dehn filling (surgery) theorem,
which says that \lq\lq almost all'' Dehn fillings of an
orientable cusped hyperbolic $3$-manifold
produce complete hyperbolic manifolds.
Since every closed orientable $3$-manifold is obtained by 
Dehn surgery of a hyperbolic link,
the theorem implies that
\lq\lq almost all closed orientable $3$-manifolds'' are hyperbolic.
This gave strong evidence for Thurston's geometrization conjecture,
which was eventually proved by Perelman.
The natural and important problem of the study of the exceptional surgeries
of hyperbolic knot complements attracted the attention
of many mathematicians and numerous research was made on this problem.
Due to the development of Heegaard-Floer homology,
this problem now attracts renewed interest.

The other flexibility of the $3$-dimensional hyperbolic structure
is that of complete hyperbolic structures of infinite volume,
in other words, the flexibility of complete hyperbolic structures
on the interior of a compact orientable $3$-manifold 
whose boundary contains a component with negative Euler characteristic.
The deformation theory of such structures is the heart of Kleinian group theory,
and it is this flexibility that enabled Thurston to prove
the hyperbolization theorem of atoroidal Haken $3$-manifolds.
In particular, the complete hyperbolic structure of a 
surface bundle over $S^1$ (with pseudo-Anosov mondromy)
was constructed by developing the deformation theory
of the complete hyperbolic structures on $\Sigma\times\RR$,
where $\Sigma$ is the fiber surface.
The idea of the Cannon--Thurston map,
a $\pi_1(\Sigma)$-equivariant sphere filling curve,
naturally arose from this construction.
Thurston produced various astonishing pictures 
of (approximations of) Cannon--Thurston maps.
(See \cite[Figures 8 and 10]{Thurston2}, \cite[Figure 1]{Thurston5}
and the beautiful book \cite{MSW} by Mumford-Series-Wright.)
It was indeed a shocking event for the author of this survey
(who was ignorant of deformation theory
and had no idea that it has something to do with knot theory)
to learn that a simple topological object, such as the figure-eight knot,
carries such mysterious mathematics under cover.

In conclusion, the flexibility of $3$-dimensional hyperbolic structure
has enriched knot theory
by bringing 
the concept of deformation
into knot theory.
(See Sections \ref{sec:Hyperbolic-Dehn-filling} and \ref{sec:flexibility}.)

\medskip

In this review, we also consider the role of
the classical knot invariants,
the Alexander polynomial and the homology 
of finite branched/unbranched coverings.
After the appearance of Thurston's work and 
the Jones polynomials,
the role of these invariants in knot theory might have decreased.
However, they continue to be important themes in knot theory.
For the Alexander polynomial,
its twisted version was
defined by Lin \cite{Lin} for classical knots
and by Wada \cite{Wada} in a general setting.
For a hyperbolic knot, we can consider 
the {\it hyperbolic torsion polynomial} 
(see \cite{Dunfield-Friedl-Jackson})
as the most natural twisted Alexander polynomial,
and a
beautiful {\it Thurstonian connection}
(cf. \cite[Section 1.2]{Agol-Dunfield})
between the topology and geometry of knots is found
(see Subsection \ref{subsec:twisted-A-plynomial}).
For the homology 
of finite branched/unbranched coverings of a knot,
Thang Le \cite{Le2018} proved a mysterious relation 
between the asymptotic growth of the order of the torsion part
and the Gromov norm of the knot (see Subsection \ref{subsec:homology-growth}).
This result is particularly surprising
to the author of this survey,
for whom homology of finite coverings is a favorite invariant, 
but who had never imagined that the whole family of the familiar invariant
could contain such deep geometric information.

\medskip

{\bf Acknowledgements.}
The author would like to thank Ken'ichi Ohshika
for giving him the challenging opportunity
to survey the tremendous impact of Thurston's work on knot theory.
He is also grateful to
Fran\c{c}ois Gu\'eritaud, Luisa Paoluzzi, 
Kenneth A. Perko Jr., 
and Han Yoshida
for correcting errors and providing valuable comments
for
Sections \ref{sec:flexibility}, 
\ref{sec:Bonahon-Siebenmann_orbifold-thm},
\ref{sec:Before-Thurston} and \ref{sec:arithmetic-invariant},
respectively.
The author would also like to thank
Yuya Koda, Gaven Martin, and Hitoshi Murakami 
for reading through an early version and 
for sending him a large number of valuable suggestions and corrections.
The author's thanks also go to 
Hirotaka Akiyoshi,
Warren Dicks,
Hiroshi Goda,
Kazuhiro Ichihara, 
Yuichi Kabaya, 
Takuya Katayama,
Akio Kawauchi, 
Eiko Kin,
Thang Le,
Hidetoshi Masai, 
Jos\'e Mar\'ia Montesinos,
Kimihiko Motegi,
Kunio Murasugi,
Shunsuke Sakai,
Masakazu Teragaito,
Ken'ichi Yoshida,
and
Bruno Zimmermann
for their valuable information and suggestions
on early versions of this survey.
The author is grateful to the referee
for his/her very careful reading and valuable suggestions,
including Remark \ref{rem.parameter-space}(1). 
Finally, the author would like to thank Athanase Papadopoulos 
for his extremely careful check of the final draft.

The author was supported by JSPS Grants-in-Aid 15H03620, 20K03614, and by Osaka City University Advanced Mathematical Institute (MEXT Joint Usage/Research Center on Mathematics and Theoretical Physics JPMXP0619217849).

\section{Knot theory before Thurston}
\label{sec:Before-Thurston}

In this section, 
we recall basic definitions 
and the classical results in knot theory,
mostly obtained before 
knot theory was influenced by Thurston's work:
(i) genera of knots,
(ii) Schubert's unique prime decomposition theorem,
(iii) knot groups, consequences of Waldhausen's work on Haken manifolds,
and the Gordon-Lueke knot complement theorem,
(iv) fibered knots and open book decompositions,
(v) the definition of the Alexander polynomial
and its effectiveness and weakness, and 
(vi) representations of knot groups in finite groups.

The book of Adams \cite{Adams2004} is a wonderful introduction to knot theory.
For classical results in knot theory, see the textbooks
Crowell-Fox \cite{Crowell-Fox},
Rolfsen \cite{Rolfsen},
Kauffman \cite{Kauffman1987}, 
Burde-Zieschang \cite{Burde-Zieschang}, 
Kawauchi \cite{Kawauchi1996},
Murasugi \cite{Murasugi4}, 
Lickorish \cite{Lickorish}, Livingston \cite{Livingston}, 
Prasolov-Sossinsky \cite{Prasolov-Sossinsky},
Cromwell \cite{Cromwell2004} and
Burde-Zieschang-Heusener \cite{Burde-Zieschang-Heusener}.
See also
the special issue edited by Adams \cite{Adams1998}
and the handbook Menasco-Thistlethwaite \cite{Menasco-Thistlethwaite}

\subsection{The fundamental problem in knot theory}
A {\it knot}\index{knot} $K$ is a smoothly (or piecewise-linearly) 
embedded circle in the $3$-sphere 
$S^3=\{(z_1,z_2)\in\CC^2 \ | \ |z_1|^2+|z_2|^2=1\}$.
Two knots $K$ and $K'$ are said to be {\it equivalent}, denoted by $K\cong K'$,
if there is a self-homeomorphism $f$ of $S^3$
such that $f(K)=K'$, 
i.e., the pair $(S^3,K)$ is homeomorphic to the pair $(S^3,K')$.
If the homeomorphism $f$ preserves the orientation of $S^3$
and hence is isotopic to the identity homeomorphsim $1_{S^3}$,
then $K$ and $K'$ are said to be {\it isotopic}.
We do not distinguish between a knot $K$
and the equivalence/isotopy class represented by $K$.
A knot is {\it trivial}\index{knot!trivial} if it is isotopic
to a standard circle $O:=\{(z_1,0)\in S^3 \ | \ |z_1|=1\}$.

Every knot is represented by a {\it knot diagram}\index{diagram}\index{knot!diagram},
a $4$-valent plannar graph whose vertices are endowed with over/under information.
A vertex of a knot diagram with over/under information is called a {\it crossing}. 

For a knot $K$, we denote by $K^*$
the image of $K$ by an orientation-reversing homeomorphism of $S^3$,
and call it the {\it mirror image} of $K$.
$K^*$ is represented by the knot diagram
which is obtained from that of $K$ by reversing the over/under information
at every crossing.
A knot $K$ is {\it achiral}\index{achiral}\index{amphicheiral}\index{knot!achiral} 
(or {\it amphicheiral})\footnote
{
This follows  
\cite{Crowell-Fox},
though the terminology \lq\lq amphichiral'' seems to be more popular.
}
if $K^*$ is isotopic to $K$,
otherwise it is {\it chiral}\index{chiral}\index{knot!chiral}.

\begin{figure}[ht]
\begin{center}
 {
\includegraphics[height=8cm]{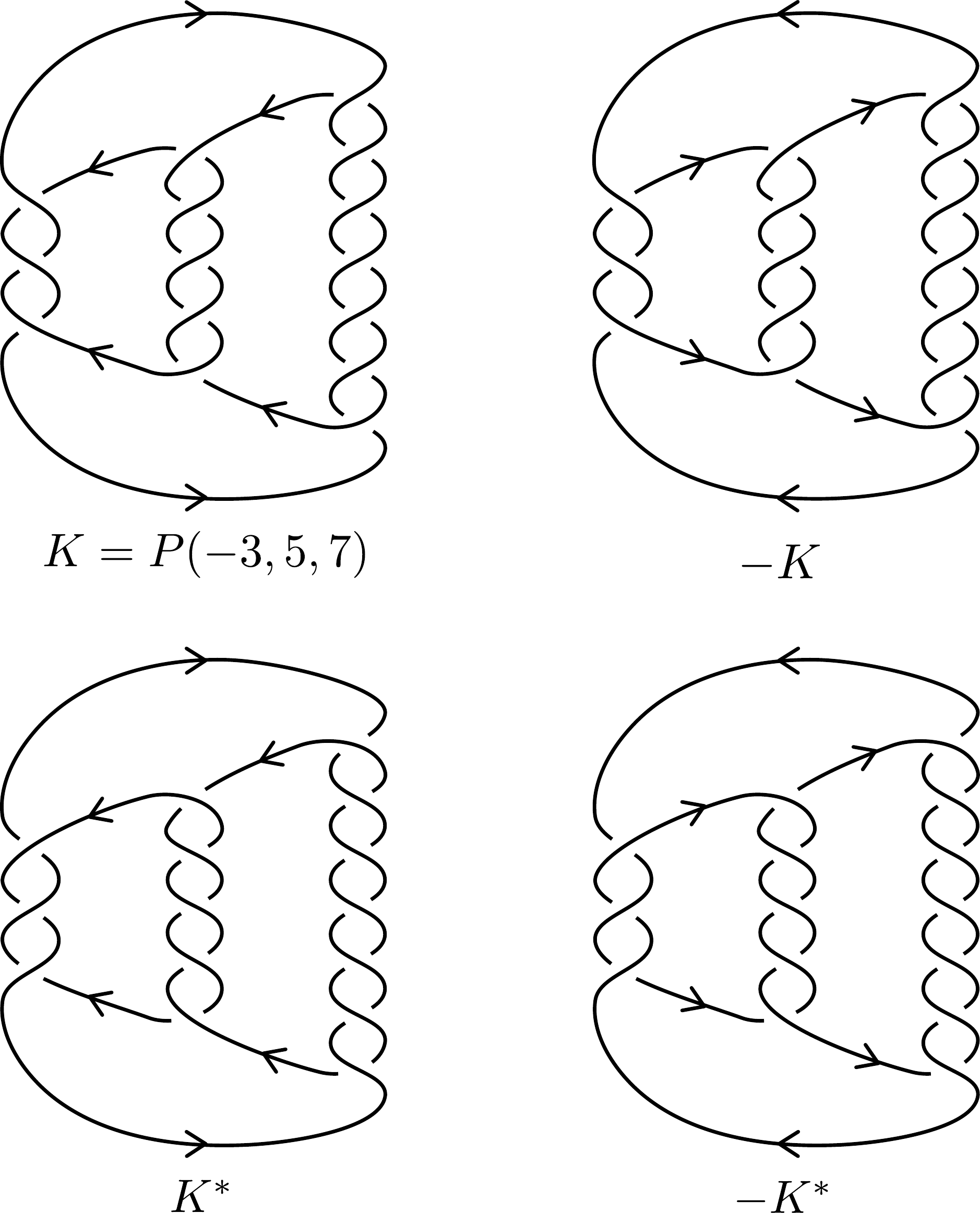}
 }
\end{center}
\caption{The pretzel knot $P(-3,5,7)$ with four different orientations.
These oriented knots are not non-isotopic to each other.}
\label{fig:oriented-knots}
\end{figure}

An {\it oriented knot}\index{knot!oriented} is a knot $K$ where the circle $K$ is 
also endowed with an orientation.
(We assume that $S^3$ is endowed with the standard orientation.)
Two oriented knots $K$ and $K'$ are said to be {\it isotopic}
if there is an orientation-preserving
self-homeomorphism $f$ of $S^3$ with $f(K)=K'$
such that $f|_{K}:K\to K'$ is also orientation-preserving.
This is equivalent to the condition that 
there is an isotopy of $S^3$ which carries the oriented circle $K$
to the oriented circle $K'$.
For a given oriented knot $K$,
we obtain the following three (possibly isotopic) oriented knots,
by reversing one or both of the orientations of $S^3$ and the circle $K$
(see Figure \ref{fig:oriented-knots}):
\[
-K:=(S^3,-K), \quad K^*:=(-S^3,K)\cong (S^3,K^*), \quad -K^*:=(-S^3,-K)\cong (S^3,-K^*)
\]

\begin{figure}[ht]
\includegraphics[height=4cm]{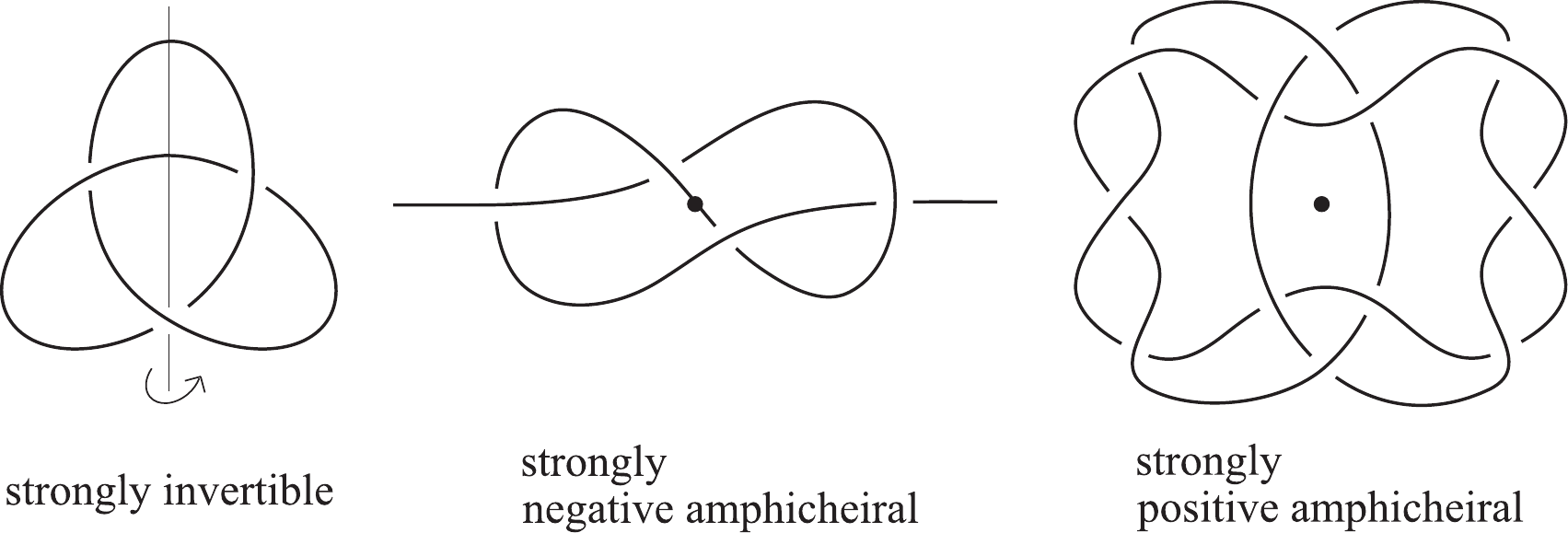}
\caption{Symmetries of knots realized by involutions}\label{Figure 1}
\end{figure}

A knot $K$ is {\it invertible}\index{invertible}\index{knot!invertible}, 
{\it positive-amphicheiral}, or 
{\it negative-amphicheiral}, respectively,
if $K$ is isotopic to $-K$, $K^*$, or $-K^*$.
If the symmetry can realize by an involution,
then we say that $K$ is
{\it strongly invertible}, {\it  strongly positive-amphicheiral}, or 
{\it strongly negative-amphicheiral}, respectively (see Figure \ref{Figure 1})..

It is one of the most fundamental problems in knot theory
to detect whether two given knots $K$ and $K'$ are equivalent or not,
in particular if a given knot $K$ is trivial or not.
The problem of detecting whether a given knot is chiral (or invertible)
is a special case of a refinement for oriented knots
of this fundamental problem.

To end this section, we note that the first proof of the existence 
of non-invertible knots due to Trotter \cite{Trotter} 
essentially uses $2$-dimensional hyperbolic geometry
(see the paragraph after Theorem \ref{thm:pi-orbifold}).

\subsection{Seifert surface}
\begin{figure}[ht]
\begin{center}
 {
\includegraphics[height=4.5cm]{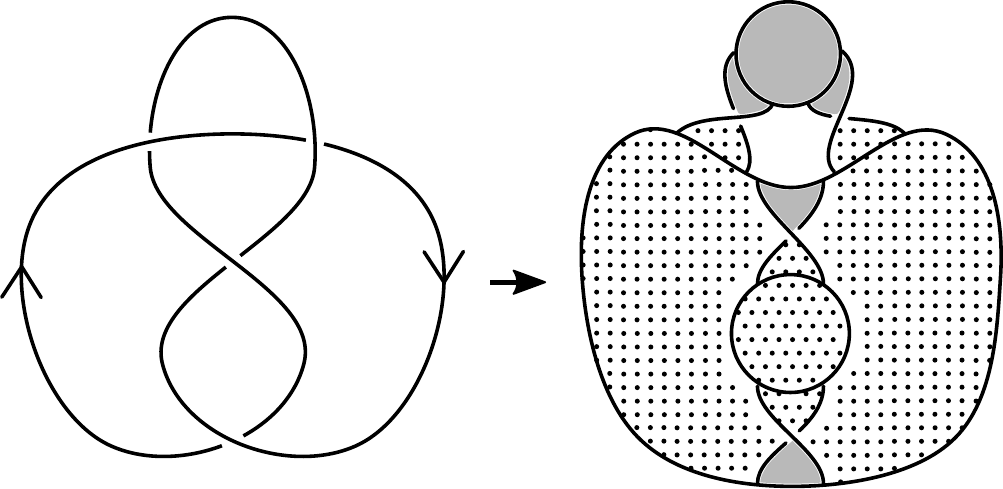}
 }
\end{center}
\caption
{
Seifert algorithm:
By smoothing all crossings of a knot diagram,
we obtain {\it Seifert circles} (mutually disjoint circles in the plane).
Construct mutually disjoint disks in $\RR^3$ bounded by the Seifert circles,
and join them by bands. The resulting surface is a Seifert surface.
}
\label{fig:Seifert-surface}
\end{figure}

A {\it Seifert surface}\index{Seifert surface} of a knot $K$ in $S^3$ is a 
connected compact orientable surface $\Sigma$ in $S^3$ with $\partial \Sigma=K$.
The existence of a Seifert surface was first proved by 
Frankel and Pontryagin \cite{Frankel-Pontryagin},
through a smooth map $f:S^3-K \to S^1$ which represents a generator 
of $H^1(S^3-K;\ZZ)\cong \ZZ$ as the closure
of the inverse image $f^{-1}(b)$ of a regular point $b\in S^1$.
Later, Seifert gave a simple effective method, 
called the {\it Seifert algorithm}\index{Seifert algorithm},
for constructing a Seifert surface from an oriented knot diagram
(see Figure \ref{fig:Seifert-surface}).
The {\it genus}\index{genus}\index{knot!genus} $g(K)$ of a knot $K$ is the minimum of the genera of
Seifert surfaces for $K$.
This is one of the most fundamental invariants of a knot,
generalized by Thurston to the concept of Thurston norm. 
The trivial knot $O$ is characterized by the property $g(O)=0$.

\subsection{The unique prime decomposition of a knot}
\begin{figure}[ht]
\begin{center}
 {
\includegraphics[height=2.5cm]{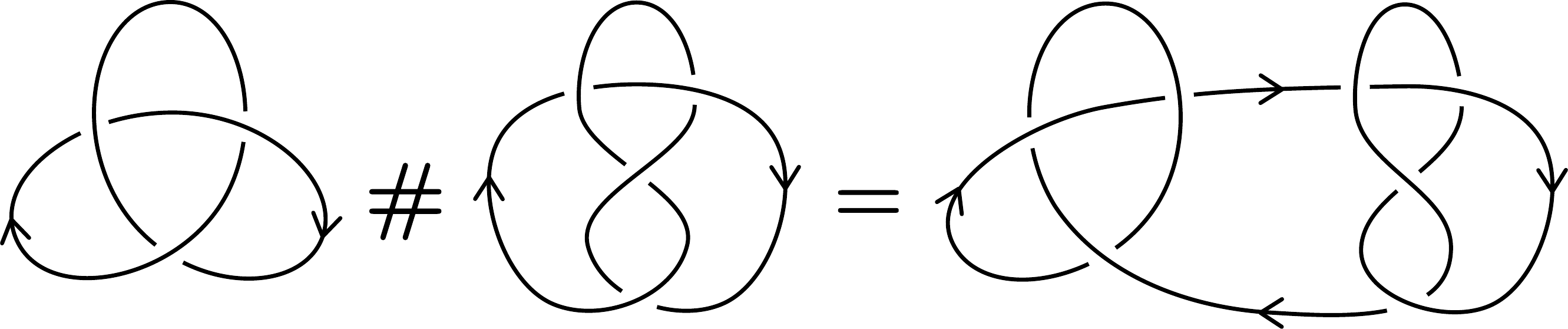}
 }
\end{center}
\caption
{
Connected sum of knots
}
\label{fig:connected-sum}
\end{figure}

We recall Shubert's unique prime decomposition theorem,
which reduces the classification problem of knots to that of prime knots.
Given two oriented knots $K_1$ and $K_2$,
we can define the {\it composition}
$K_1\# K_2$ as the pairwise connected sum
$(S^3,K_1)\# (S^3,K_2)$ of oriented manifold pairs,
as in Figure \ref{fig:connected-sum}.
With respect to the connected sum,
the set of all oriented knots up to isotopy
becomes a {\it commutative} semi-group 
having the trivial knot $O$ as the unit.

A knot $K$ is {\it prime}\index{prime}\index{knot!prime}
if $K\cong K_1\# K_2$ implies $K_1\cong O$ or $K_2\cong O$.
It is a classical theorem due to Schubert \cite{Schubert}
that every oriented knot has a unique prime decomposition.

\begin{theorem}
[The unique prime decomposition of knots]\index{unique prime decomposition}
\label{thm:upd-knot}
Every nontrivial oriented knot $(S^3,K)$ can be decomposed as the sum
of finitely many nontrivial prime oriented knots.
Moreover if $K\cong K_1\# K_2\# \cdots \# K_n$ and 
$K\cong J_1\# J_2\# \cdots \# J_m$
with each $K_i$ and $J_i$ nontrivial prime knots,
then $m=n$, and after reordering, $K_i\cong J_i$ as oriented knots.
\end{theorem}

The existence of a prime decomposition is guaranteed by 
the additivity of genus with respect to connected sum, i.e., 
\[
g(K_1\#K_2)= g(K_1)+g(K_2) \quad\mbox{for any oriented knots $K_1$ and $K_2$.}
\]
The uniqueness of the prime decomposition is proved by a simple cut and paste argument.

\subsection{Knot complements and knot groups}
The {\it exterior}\index{exterior}\index{knot!exterior} of a knot $K$ is defined by
$E(K):=S^3 - \interior N(K)$,
where $N(K)$ is a regular neighborhood of $K$.
The {\it knot complement}\index{knot complement} $S^3-K$ is homeomorphic to the interior of $E(K)$, and
the fundamental group $\pi_1(S^3-K)\cong \pi_1(E(K))$ 
is called the {\it knot group}\index{knot group}, and denoted by $G(K)$.
By using the sphere theorem
\cite{Papakyriakopoulos} (cf. \cite[Chapter 4]{Hempel1}),
we can see that $E(K)$ is aspherical,
and hence the homotopy type of $E(K)$ is completely determined by
the knot group $G(K)$. 
A group presentation, called the
{\it Wirtinger presentation}, of $G(K)$ 
can be obtained from a knot diagram of $K$
(see \cite{Crowell-Fox, Rolfsen}).

The {\it peripheral subgroup}\index{peripheral subgroup}\index{knot group!peripheral subgroup} 
$P(K)$ of the knot group $G(K)$
is defined as (the conjugacy class of) the image of the homomorphism
$j_*:\pi_1(\partial E(K))\to \pi_1(E(K))$
induced by the inclusion map $j:\partial E(K) \to E(K)$.
For the trivial knot $O$, $E(O)$ is homeomorphic to the solid torus,
and so $G(O)=P(O)\cong \ZZ$. 
Dehn's lemma, established by Papakyriakopoulos \cite{Papakyriakopoulos},
gives the following characterization 
of the trivial knot.

\begin{theorem}
A knot $K$ is trivial if and only if 
the following mutually equivalent
conditions hold.
\begin{enumerate}
\item
$\Ker[j_*:\pi_1(\partial E(K))\to \pi_1(E(K))]$ is nontrivial.
\item
The peripheral subgroup $P(K)$ is isomorphic to $\ZZ$.
\item
The knot group $G(K)$ is isomorphic to $\ZZ$.
\end{enumerate}
\end{theorem} 
 
If $K$ is a nontrivial knot, then 
the peripheral subgroup $P(K)\cong \pi_1(\partial E(K))\cong \ZZ\oplus\ZZ$ is generated by two special elements,
a {\em meridian}\index{meridian}\index{knot group!meridian} $\mu$ and 
a {\em longitude}\index{knot}\index{knot group!longitude} $\lambda$,
represented by the simple loops $\mu:=\partial D^2 \times\{*\}$ 
and $\lambda:=\{*\}\times S^1$ respectively, 
in $\partial E(K)=\partial N(K)=\partial (D^2\times S^1)$.
Here the framing $N(K)\cong D^2\times S^1$ is chosen so that
the linking number $\lk(K,\lambda)=0$.
When $K$ is oriented, the orientations of $\mu$ and $\lambda$ are chosen
so that $\lk(K, \mu)=+1$ and that $K$ and $\lambda$ are homologous in $N(K)$.  

The classical work of Waldhausen \cite{Waldhausen} on Haken manifolds
implies the following theorem which reduces the equivalence problem for knots
to a problem of knot groups.

\begin{theorem}
For two knots $K$ and $K'$, the following hold.

{\rm (1)} $E(K)$ and $E(K')$ are homeomorphic
if and only if $(G(K), P(K))$ and $(G(K'),P(K'))$ are isomorphic,
i.e., there is an isomorphism $\varphi:G(K)\to G(K')$ 
such that $\varphi(P(K))=P(K')$ up to conjugacy.

{\rm (2)} $K$ and $K'$ are equivalent
if and only if $(G(K), P(K), \mu)$ and $(G(K'),P(K'),\mu'^{\pm 1})$ are isomorphic,
i.e., there is an isomorphism $\varphi:(G(K), P(K))\to (G(K'),P(K'))$
such that $\varphi(\mu)=\mu'^{\pm 1}$ up to conjugacy.
\end{theorem}

For nontrivial oriented knots $K_1$ and $K_2$,
the knot groups $G(K_1\#K_2)$ and $G(K_1\#(-K_2^*))$ are isomorphic.
In fact, both $E(K_1\#K_2)$ and $E(K_1\#(-K_2^*))$ are obtained from
$E(K_1)$ and $E(K_2)$ by gluing annuli in their boundaries,
and so homotopy equivalent to the space obtained from 
$E(K_1)$ and $E(K_2)$ by identifying the meridians $\mu_1$ and $\mu_2$.
On the other hand, by the unique prime decomposition Theorem \ref{thm:upd-knot},
the oriented knots $K_1\#K_2$ and $K_1\#(-K_2^*)$ are isotopic
if and only if $K_2$ is negative amphicheiral
(i.e., $-K_2^*$ is isotopic to $K_2$). 
Thus, in general, the knot group alone is not a complete invariant for knots.

Building on the cyclic surgery theorem (Theorem \ref{thm:cyclic-surgery})
by Culler, Gordon, Luecke and Shalen \cite{CGLS},
Whitten \cite{Whitten}
proved that prime knots with isomorphic knot groups have homeomorphic
exteriors. 
On the other hand, we have
the following celebrated theorem of Gordon and Luecke \cite{Gordon-Luecke}.

\begin{theorem}[Knot Complement Theorem]\index{knot complement theorem}
\label{thm:knot-complement-thm}
Two knots are equivalent if and only if they have homeomorphic complements.
\end{theorem}

Thus we have the following theorem.

\begin{theorem}
Two prime knots are equivalent 
if and only if they have isomorphic knot groups.
\end{theorem}

\subsection{Fibered knots}
A knot $K$ is {\it fibered}\index{fibered}\index{knot!fibered} if $E(K)$ has the structure of a bundle over the circle,
namely, there is a connected compact orientable surface $\Sigma$
and an orientation-preserving homeomorphism $\varphi:\Sigma\to \Sigma$, such that
\[
E(K)\cong \Sigma \times [0,1]/(x,0)\sim (\varphi(x),1). 
\quad
\]
The homeomorphism $\varphi$ is called the {\it monodromy}\index{monodromy} 
of the fiber structure.
Each fiber $\Sigma$ of the bundle structure is a compact orientable surface
in $E(K)$ such that $\Sigma\cap \partial E(K)= \partial \Sigma$
is a longitude of $K$.
The union of $\Sigma$ and an
annulus in $N(K)$ cobounded by $\partial \Sigma$ and $K$ is a minimal genus Seifert surface for $K$.
This is the unique minimal genus Seifert surface for $K$ up to isotopy fixing $K$.

We may choose $\varphi$ so that its restriction to $\partial \Sigma$ is the identity map
and thus the image of $y\times [0,1]$ in $E(K)$ is a meridian of $K$ for every $y\in \partial \Sigma$.
Then 
\[
(S^3,K) \cong (\Sigma,\partial \Sigma) \times [0,1]/
[(x,0)\sim (\varphi(x),1);\ y\times [0,1]\sim y \ \mbox{(for $y\in \partial \Sigma$)}].
\]
This structure is called an {\it open book decomposition} with {\it binding} $K$,
and the homeomorphism $\varphi$ is called the {\it monodromy} of the fibered knot $K$. 
It was proved by Alexander \cite{Alexander} that every connected closed orientable $3$-manifold
admits an open book decomposition.
Later, Giroux \cite{Giroux} found a very important correspondence 
between the open book decompositions (up to positive stabilization)
of a given closed oriented $3$-manifold $M$
and oriented contact structures on $M$ up to isotopy
(see \cite{Etnyre} for details).
The following characterization of fibered knots in terms of knot groups
was proved by Stallings \cite{Stallings1},
and attracted the attention of researchers at the time.

\begin{theorem}
A knot $K$ in $S^3$ is a fibered knot
if and only if the commutator subgroup $G(K)'=[G(K),G(K)]$
is finitely generated.
\end{theorem}

The only if part follows from the fact that
the infinite cyclic covering $E_{\infty}(K)$ of $E(K)$,
introduced in the subsection below,
is identified with $\Sigma\times \RR$
and so $G(K)'\cong \pi_1(E_{\infty}(K))\cong \pi_1(\Sigma)$ is 
a free group of rank $2g(K)$.
The heart of the theorem is that the converse also holds.
  
\subsection{Alexander invariants}
\label{subsec:Alexander-polynomial}
Though the knot group is a complete invariant for prime knots,  
it is, in general, not easy to distinguish two given knot groups.
The Alexander polynomial serves as a convenient and tractable tool  
for this problem, even though it is not almighty.

Let $K$ be an oriented knot.
Then the first integral homology group $H_1(E(K);\ZZ)$ is the infinite cyclic group generated by the image, $t$,
of the meridian $\mu$. 
Thus there is a unique infinite cyclic covering $p_{\infty}:E_{\infty}(K) \to E(K)$,
and the covering transformation group is identified
with the infinite cyclic group $\langle t\rangle$ generated by $t$.
$H_1(E_{\infty}(K);\ZZ)$ has the structure of
a module over the integral group ring $\ZZ\langle t\rangle$.
This module is called the {\it knot module}\index{knot module}.
As an abelian group, $H_1(E_{\infty}(K);\ZZ)$ is identified with $G(K)'/G(K)''$
where $G(K)'$ and $G(K)''$, respectively, 
are the first and second commutator subgroups of $G(K)$.
Moreover the action of the generator $t$ is given by
$t[\alpha]=[\mu \alpha \mu^{-1}]$ for $\alpha\in G(K)'$,
where $\mu$ is a meridian.
Thus the knot module is determined by $G(K)$.
In fact, a presentation matrix is obtained from a presentation of 
the knot group, via Fox's free differential calculus (see \cite{Crowell-Fox},
\cite[Chapter 7]{Kawauchi1996}).
The {\it Alexander polynomial}\index{Alexander polynomial}
 $\Delta_K(t)$ of $K$ is defined as the generator 
of the first elementary ideal of the knot module.

A more conceptual definition can be given
by using the $\QQ\langle t\rangle$-module
$H_1(E_{\infty}(K);\QQ)$ as follows.
Since the rational group ring $\QQ\langle t\rangle$ is a principal ideal domain
and since $H_1(E_{\infty}(K);\QQ)$ is a finitely generated torsion module over 
$\QQ\langle t\rangle$, we have
\[
H_1(E_{\infty}(K);\QQ)\cong 
\frac{\QQ\langle t\rangle}{(f_1(t))}\oplus\cdots\oplus 
\frac{\QQ\langle t\rangle}{(f_r(t))},
\]
where $f_i(t)$ are elements of $\ZZ\langle t\rangle$
whose coefficients are relatively prime. 
Then $\Delta_K(t)\doteq f_1(t)\cdots f_r(t)$,
where $\doteq$ means equality up to multiplication by a unit
$\pm t^i$ of the integral Laurent polynomial ring $\ZZ\langle t\rangle$.
The Alexander polynomial $\Delta_K(t)$ is an integral Laurent polynomial in the variable $t$,
defined up to multiplication by a unit.
For the trivial knot $O$, we have $\Delta_O(t) \doteq 1$.
We summarize basic properties of the Alexander polynomial.

\begin{theorem}
\label{thm:Alexander-polynomial}
{\rm (1)}
For any knot $K$, its Alexander polynomial $\Delta_K(t)$ satisfies the following condition.
\[
\Delta_K(1)=\pm 1, \quad \Delta_K(t^{-1})\doteq \Delta_K(t)
\]
Conversely, for any Laurent polynomial $\Delta(t)$ satisfying the above condition,
there is a knot $K$ whose Alexander polynomial is equal to $\Delta(t)$.

{\rm (2)}
For every knot $K$ in $S^3$, we have the following estimate of the genus:
\[
g(K)\ge \deg \Delta_K(t).
\]

{\rm (3)}
For any fibered knot $K$, the Alexander polynomial $\Delta_K(t)$ is monic,
and the equality hold in the estimate (2).
\end{theorem}

Proof of Theorem \ref{thm:Alexander-polynomial}.
The proof 
relies on an analysis of the manifold $M:=E(K)\backslash \Sigma$,
the manifold obtained from $E(K)$ by cutting along a Seifert surface $\Sigma$;
in other words, $M$ is the complement of an open regular neighborhood of $\Sigma$ in $E(K)$.
Let $\Sigma_+$ and $\Sigma_-$ be copies of $\Sigma$ on $\partial M_{\Sigma}$, and 
consider the annulus $A:=M\cap \partial E(K)$.
Then $(M,\Sigma_+,\Sigma_-,A)$ is 
a {\it sutured manifold}\index{sutured manifold} 
(see \cite{Gabai1983a, Gabai1984}, \cite[Chapter 5]{Kawauchi1996}),
and this together with the natural homeomorphism $\Sigma_+\to \Sigma_-$
recovers $E(K)$.
The infinite cyclic covering $E_{\infty}(K)$ is obtained 
from the set of copies $\{M_n\}$ of $M$ indexed with $n\in\ZZ$,
by gluing the copy of $\Sigma_-$ in $M_n$ with the copy of $F_+$ in $M_{n+1}$.
The homological glueing information is given 
by the {\it Seifert matrix}\index{Seifert matrix} 
$V=(\mathrm{lk}(\alpha_i,\alpha_j^+))_{1\le i,j \le 2g}$,
where $\{\alpha_i\}_{1\le i,j \le 2g}$ with $g=2g(\Sigma)$ 
is a set of oriented simple loops on $\Sigma$
which forms a basis of $H_1(\Sigma)$,
$\alpha_j^+$ is a copy of $\alpha_j$ on the $+$-side of $\Sigma$,
and $\mathrm{lk}(\cdot,\cdot)$ denotes the linking number.
The matrix $tV-V^{T}$ gives a presentation matrix of $H_1(E_{\infty}(K))$
as a $\ZZ\langle t\rangle$-module,
and hence $\Delta_K(t)=\mathrm{det}(tV-V^{T})$.
Using this formula we can prove Theorem \ref{thm:Alexander-polynomial}.

For knots with small crossing numbers,
the Alexander polynomial is quite efficient.
For any prime knot $K$ up to 10 crossings, 
equality holds in the estimate of the genus
in Theorem \ref{thm:Alexander-polynomial}(2).
Moreover, such a knot $K$ is fibered if and only if
$\Delta_K(t)$ is monic (see Kanenobu \cite{Kanenobu}).

The Alexander polynomial is also very efficient for alternating knots.
A knot $K$ is said to be {\it alternating}\index{alternating}\index{knot!alternating} 
if 
it is represented by an {\it alternating diagram}\index{diagram!alternating},
namely a diagram
in which the crossings alternate 
under and over
as one travels along 
the diagram.
A knot diagram is said to be {\it reduced}\index{reduced}\index{diagram!reduced}
if there is no circle in the plane which intersects the diagram 
only at a single crossing.
Any alternating digram can be deformed into a (possible trivial)
reduced alternating diagram.
When $K$ is an alternating knot and $\Sigma$ is a Seifert surface
obtained by the Seifert algorithm from a reduced alternating diagram of $K$,
the complementary sutured manifold $(M,\Sigma_+,\Sigma_-,A)$ has a nice structure,
which in particular implies $\mathrm{det}(V)\ne 0$.
This shows that the estimate Theorem \ref{thm:Alexander-polynomial}(2) 
is sharp for alternating knots
(see Crowell \cite{Crowell} and Murasugi \cite{Murasugi1}).
Moreover, Murasugi \cite{Murasugi2} proved
that the converse to Theorem \ref{thm:Alexander-polynomial}(3) also holds
for alternating knots.

\begin{theorem}
\label{thm:alternating-Alexander}
For any alternating knot $K$, the following hold.
\begin{enumerate}
\item
$g(K)=\deg \Delta_K(t)$.
\item
$K$ is fibered if and only if $\Delta_K(t)$ is monic.
\end{enumerate}
\end{theorem}

In order to prove the above results,
Murasugi introduced the concept of 
a {\it Murasugi sum}\index{Murasugi sum}\index{Seifert surface!Murasugi sum} 
of two Seifert surfaces.
The simplest case corresponds to the connected sum of knots
and the second simplest case corresponds to 
{\it plumbing}\index{plumbing}\index{Seifert surface!plumbing}
introduced by Stallings \cite{Stallings2}. 
It was later shown by Gabai \cite{Gabai1983b} that
the Murasugi sum is a natural geometric operation in the following sense:
If $\Sigma$ is a Murasugi sum of $\Sigma_1$ and $\Sigma_2$,
then the following hold.

\begin{enumerate}
\item
$\Sigma$ is of minimal genus if and only if $\Sigma_1$ and $\Sigma_2$ are of minimal genus.
\item
$\Sigma$ is a fiber surface if and only if $\Sigma_1$ and $\Sigma_2$ are fiber surfaces.
\end{enumerate}

\medskip

In addition to Theorems \ref{thm:Alexander-polynomial} and \ref{thm:alternating-Alexander},
various applications of the Alexander polynomials were found.
Among them, we explain a theorem by Kinoshita \cite{Kinoshita},
which gives a condition on the Alexander polynomial
that a counter-example to the Smith Conjecture must satisfy.
As described in Subsection \ref{subsec:orbifold-theorem},
the Smith conjecture (Theorem \ref{thm:SmithConj}) was later proved using Thurston's geometrization theorem for Haken manifolds.

\begin{theorem}
If $K$ is a fixed point of an orientation-preserving periodic diffeomorphism of period $n$,
then there is an integral Laurent polynomial $f(t)$ such that
\begin{enumerate}
\item
$\Delta_K(t^n)=\Pi_{i=0}^{n-1}f(\xi^i t)$ where $\xi$ is a primitive $n$-th root of unity,
and 
\item
$f(1)=\pm 1, \quad f(t^{-1})\doteq f(t)$
\end{enumerate}
\end{theorem}

See \cite[Chapter 10]{Kawauchi1996} for other applications of the Alexander polynomial
to the study of symmetry of knots,
including the first proof of the non-invertibility of the knot $8_{17}$
by Kawauchi \cite{Kawauchi1979},
answering to a question of Fox (cf. \cite[Problem 10]{Fox_problem}).
(Another proof of the non-invertibility of $8_{17}$
was announced almost at the same time by Bonahon and Siebenmann,
based on their characteristic splitting theory
(see Section \ref{sec:Bonahon-Siebenmann_orbifold-thm}).)
It should be noted that 
though the definition of the Alexander polynomial depends on the
orientation of $K$, by Theorem \ref{thm:Alexander-polynomial}(1)
the resulting $\Delta_K(t)$ does not depend on the orientation.
It is interesting that, despite this fact,
the Alexander polynomial can be used for the study
of invertibility and chirality of knots. 
Finally, we point out that the Alexander module $H_1(E_{\infty}(K))$
does depend on the orientation of $K$,
though it is not easy to detect the dependence
(see \cite{Fox_ideal-classs}, \cite{Hillman}).

\medskip

\begin{figure}[ht]
\begin{center}
 {
\includegraphics[height=5cm]{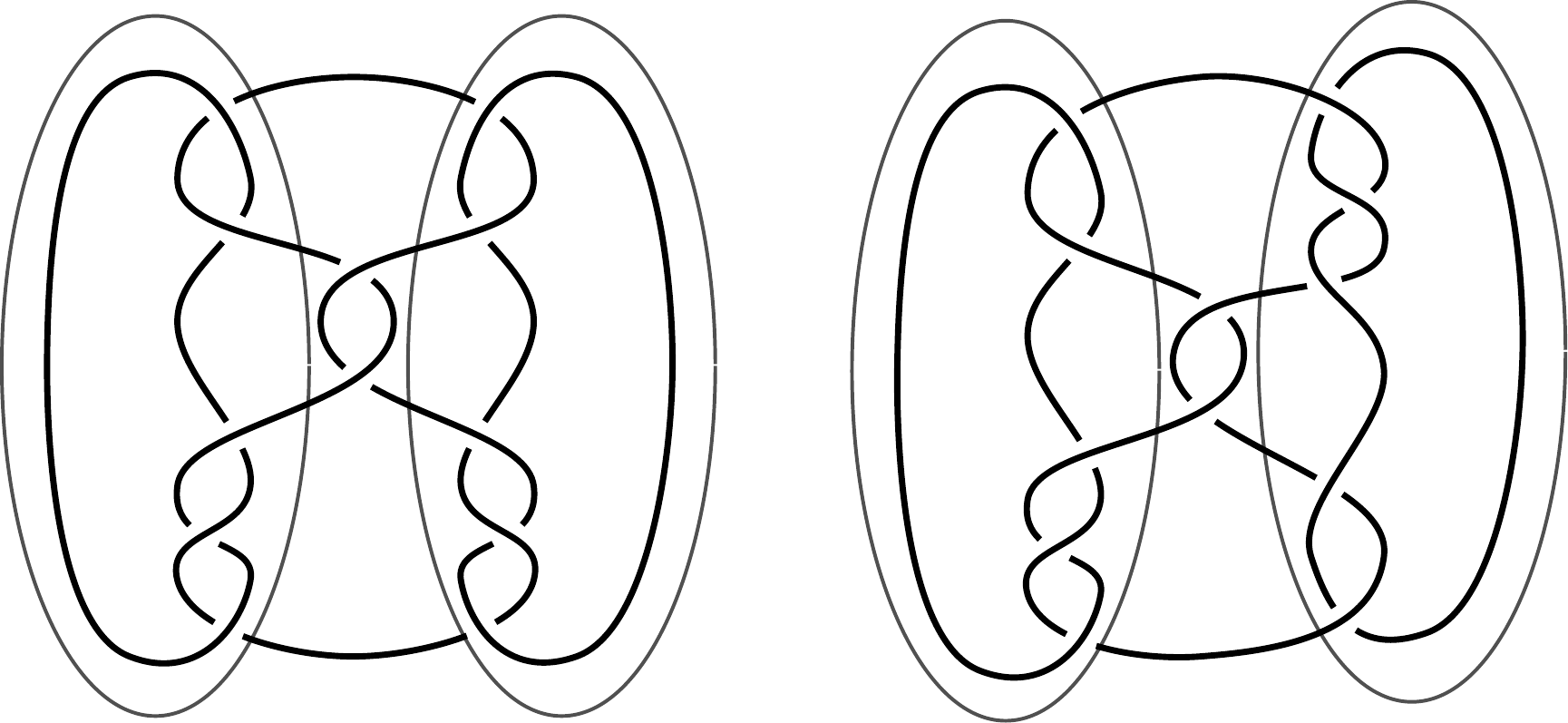}
 }
\end{center}
\caption
{
The Kinoshita--Terasaka knot and the Conway knot:
The circles in the figure represent the Conway spheres
which determine the Bonahon--Siebenmann decompostions,
described in Subsection \ref{subsec:Bonahon-Siebenmann}.
From this picture, we can see that the Conway knot is a mutant of the Kinoshita Terasaka knot.
}
\label{fig:KT-Conway}
\end{figure}

Though we have observed the effectiveness of the Alexander polynomial,
there are a lot of knots 
for which the Alexander polynomial is useless.
In fact, H. Seifert \cite{Seifert}, J. H.C. Whitehead \cite{Whitehead}, 
and Kinoshita-Terasaka \cite{Kinoshita-Terasaka}
gave systematic construction of nontrivial knots
with trivial Alexander polynomial.
For example, the pretzel knot $K(-3,5,7)$ in Figure \ref{fig:oriented-knots},
the Whitehead double of any nontrivial knot (cf. Figure \ref{fig:satellite-knot}), 
the 
{\it Kinoshita--Terasaka knot}\index{Kinoshita--Terasaka knot} and the 
{\it Conway knot}\index{Conway knot}
 in Figure \ref{fig:KT-Conway}
have trivial Alexander polynomial.

The Conway knot is a {\it mutant}\index{mutant}\index{knot!mutant} of the Kinoshita--Terasaka knot
(see Construction 4 in Subsection \ref{subsec:branched-covering}
for the precise definition).
It is known that various invariants coincide for a knot and its mutant,
including the Alexander polynomial, the Jones polynomial,
the Homflypt polynomial, the double branched covering, and Gromov norm.
So it is not easy to distinguish a knot from its mutant.

\subsection{Representations of knot groups onto finite groups}
\label{subsec:universal-method}
The definition of the Alexander polynomial is based on the fact that
the knot group $G(K)$ of an oriented knot $K$ 
admits a unique preferred epimorphism onto the infinite cyclic group $\langle t\rangle$.
By replacing $\ZZ$ with an arbitrary group $\Gamma$,
we obtain the following family of invariants of knots.
Let $R(G(K),\Gamma)$ be the set of homomorphisms from $G(K)$ to $\Gamma$, 
up to conjugacy (i.e., modulo post composition of inner-automorphisms of $\Gamma$),
is an invariant of $G(K)$.
Then its cardinality $|R(G(K),\Gamma)|$ is an invariant of 
the knot group $G(K)$.
We may also consider the quotient of $R(G(K),\Gamma)$
by the action of the automorphism group of $\Gamma$.

Fix a conjugacy class $\gamma$ of an element of $\Gamma$,
and let $R(G(K),\Gamma,\gamma)$ be
the subset of $R(G(K),\Gamma)$ consisting of the homomorphisms
which map the meridian to an element in the conjugacy class $\gamma$.
Then the cardinality $|R(G(K),\Gamma,\gamma)|$ is again an invariant of the 
oriented knot $K$.
If $\Gamma$ is the dihedral group $D_{2p}=\langle a, t \ | \ a^p, t^2, tat^{-1}=a^{-1} \rangle$
of order $2p$ and if $\gamma$ is the conjugacy class of the element $t$,
the {\it Fox $p$-coloring number}\index{Fox coloring} \cite{Fox_quick-trip}
is essentially equal to $|R(G(K),D_{2p},t)|$.

Fix a transitive representation of $\Gamma$ into the symmetric group $S_n$ of degree $n$,
where $n$ is possibly infinite.
Then for each $\phi\in R(G(K),\Gamma)$,
we have a (possibly disconnected) $n$-fold covering $E_{\phi}(K)$ of $E(K)$.
Then the family of homology groups $\{H_1(E_{\phi}(K);\ZZ)\}_{\phi\in R(G(K),\Gamma)}$ 
forms an invariant of the knot $K$.
Furthermore, if the image of the element $\gamma$ in $S_n$ is of finite order,
then we obtain a branched covering $M_{\phi}(K)$ of $S^3$ branched over $K$.
The family of homology groups $\{H_1(M_{\phi}(K);\ZZ)\}_{\phi\in R(G(K),\Gamma,\gamma)}$ 
forms another invariant of $K$.
We can also consider the torsion linking numbers 
among the components of the inverse image of $K$.

Riley \cite{Riley1971} applied this method by choosing $\Gamma$ to be 
the simplest finite simple group $\PSL(2,p)$, with $p$ a prime $\ge 5$,
and setting $\gamma$ to be the parabolic transformation
$\begin{pmatrix}
1&1\\
0&1
\end{pmatrix}$.
This enabled him to prove that the Kinoshita--Terasaka knot and the Conway knot
are different.
His proof also showed that none of them is amphicheiral.
He then considered parabolic representations of knot groups to $\PSL(2,\CC)$,
and this led him to the discovery of the complete hyperbolic structure
of the Figure-eight knot complement in \cite{Riley1975}.

Hartley \cite{Hartley2} realized that one can apply 
this method to the problem
of identifying noninvertible knots, as follows.
Suppose no automorphism of $\Gamma$ maps $\gamma$ to $\gamma^{-1}$.
Then the set $R(G(K),\Gamma,\gamma)$ is possibly different from
the set $R(G(K),\Gamma,\gamma^{-1})$,
and there is a chance to show noninvertibility of $K$ 
by comparing the homology invariants
associated with 
$\phi\in R(G(K),\Gamma,\gamma)$ 
with those associated with 
$\phi'\in R(G(K),\Gamma,\gamma^{-1})$.
Hartley showed that this method is quite effective:
he completely determined the 36 non-invertible knots
up to 10 crossings
claimed by Conway to be noninvertible. 

A variation of the method is to consider the subset
$R_{t}(G(K),S_n,\gamma)$ of $R(G(K),S_n,\gamma)$
consisting of the transitive representations, where $n$ is a finite positive integer.
By virtue of the development of computers,
this turns out to be an extremely efficient method
for distinguishing knots.
In fact, Thistlethwaite \cite{Thistlethwaite}
succeeded in distinguishing prime knots up to 13 crossings,
and later the same method was applied successfully to prime knots up to 16 crossings
in \cite{HTM}.
The tabulation was recently extended by Burton \cite{Burton} 
to prime knots up to 19 crossings, where 
he counts 352,152,252 distinct non-trivial prime knots through 19 crossings.

To end this subsection, 
we note that Perko \cite{Perko1980, Perko1982},
prior to the work of Thistlethwaite \cite{Thistlethwaite},
succeeded in classifying the prime knots up to 11 crossings
through hand calculation of homological invariants
(in particular linking invariants) 
of finite branched coverings
for those knots 
that are not covered by Bonahon and Siebenmann's result described in 
Subsection \ref{subsec:Bonahon-Siebenmann}.
See \cite{Perko2014} for an interesting historical note. 

\section{The geometric decomposition of knot exteriors}
\label{sec:geometric-decomposition}

The purpose of this section is to explain 
the geometric decompositions of knot exteriors
into Seifert pieces and hyperbolic pieces,
obtained as a special case of 
the Thurston's geometrization theorem
of Haken manifolds.

We recall (i) the prime decomposition theorem of 
general compact orientable $3$-manifolds,
(ii) the torus decomposition theorem
of compact irreducible orientable $3$-manifolds,
(iii) the eight homogeneous $3$-dimensional geometries,
and (iv) the geometrization conjecture, which was finally 
established by Perelman.
In the final subsection,
we give detailed exposition of the geometric decompositions of knot exteriors.

For standard facts in $3$-manifold theory,
see the short note Hatcher \cite{Hatcher} and 
the textbooks 
Hempel \cite{Hempel1},
Jaco \cite{Jaco},
Jaco-Shalen \cite{Jaco-Shalen}, 
Johannson \cite{Johannson} and
Schultens \cite{Schultens}.
For an introduction to geometric structures, 
see the surveys
Scott \cite{Scott1983} and Bonahon \cite{Bonahon}
and the textbook Martelli \cite{Martelli*}.

\subsection{Prime decomposition of $3$-manifolds}
In this subsection, we recall the canonical decomposition of
compact orientable $3$-manifolds by $2$-spheres.
Let $M$ be a compact connected orientable $3$-manifold.
A $2$-sphere in $M$ is {\it essential}\index{essential} if it does not bound a $3$-ball in $M$.
$M$ is {\it irreducible}\index{irreducible} if it contains no essential $2$-sphere.
Suppose $M$ is not irreducible,
and let $S$ be an essential $2$-sphere in $M$.
If $S$ is {\it separating}\index{separating}
(i.e., $M-S$ consists of two components),
then $M$ is the {\it connected sum}\index{connected sum} $M_1\#M_2$
of the two compact orientable $3$-manifolds $M_1$ and $M_2$,
which are obtained from the closures of the components
by capping off the resulting sphere boundaries
by adding $3$-balls.
If $S$ is {\it non-separating}\index{non-separating}
(i.e., $M-S$ is connected),
then $M$ is expressed as the connected sum $(S^2\times S^1)\# M'$
of $S^2\times S^1$ with some compact orientable $3$-manifold $M'$
(possibly $S^3$).

$M$ is {\it prime}\index{prime} 
if whenever $M= M_1\# M_2$ we have $M_i\cong S^3$ for $i=1$ or $2$.
Then we have the following Kneser--Milnor unique prime decomposition theorem
\cite{Kneser1929, Milnor1962} (cf. \cite{Hempel1}).

\begin{theorem}
[Unique prime decomposition of compact orientable $3$-manifolds]\index{unique prime decomposition}
\label{thm:upd-mfd}
Any compact orientable $3$-manifold admits a decomposition
$M=P_1\#\cdots \# P_n$ into prime manifolds $\{P_i\}$.
Moreover, the prime factors $\{P_i\}$ are uniquely determined by $M$,
up to change of the indices.
\end{theorem}

\subsection{Torus decomposition of irreducible $3$-manifolds}
In this subsection,
we explain the torus decomposition theorem for 
compact orientable irreducible $3$-manifolds.
Torus decomposition is a simple version of more intricate
JSJ (Jaco--Shalen--Johannson) decomposition,
in which decompositions along annuli are   also involved.
The JSJ decomposition theory grew out of the study to understand
homotopy equivalences among $3$-manifolds,
and its simplified version, the torus decomposition, 
turned out to be a complete obstruction 
for the hyperbolization of a $3$-manifold.

Throughout this subsection,
$\Sigma$ denotes a compact orientable surface in $M$
which is {\it properly embedded}\index{properly embedded} in $M$,
i.e., $\Sigma\cap M=\partial\Sigma$.
We also assume that $\Sigma\not\cong S^2$.
Then $\Sigma$ is {\it incompressible}\index{incompressible} in $M$ 
if for any disk $D$ in $M$ such that $D\cap \Sigma=\partial \Sigma$,
the simple loop $\partial D$ bounds a disk in $\Sigma$.
By the loop theorem (see \cite{Hempel1}), this is equivalent to the algebraic condition that
the homomorphism $j_*:\pi_1(\Sigma)\to \pi_1(M)$ 
induced by the inclusion is injective. 

$M$ is {\it Haken}\index{Haken} if
it is irreducible and contains a properly embedded compact orientable surface 
which is incompressible.

A surface $\Sigma$ in $M$ is {\it essential}\index{essential}
if it is incompressible and is not {\it $\partial$-parallel}\index{$\partial$-parallel},
i.e., $\Sigma$ does not cut off a $3$-manifold in $M$
homeomorphic to $\Sigma\times I$.
$M$ is {\it atoroidal}\index{atoroidal}
if it does not contain an essential torus.

$M$ is a {\it Seifert fibered space}\index{Seifert fibered space}
if it is expressed as a union of disjoint circles, in a particular way.
The quotient of $M$ obtained by collapsing each fiber into a point
has the structure of a $2$-dimensional orbifold, 
and is called the {\it base orbifold}\index{base orbifold}.  
If $M$ admits a smooth $S^1$ action without a fixed point
(i.e., the stabilzer of any point is not the whole group $S^1$),
then $M$ is a Seifert fibered space whose base orbifold is 
the orbit space $M/S^1$.
Seifert fibered spaces are regarded as $S^1$-bundles
over $2$-dimensional orbifolds,
and are completely described by 
the {\it Seifert invariants}.
See \cite{Hempel1, Scott1983} for details.

Now we state the torus decomposition theorem,
which is a simplified version of the JSJ decomposition theorem\index{JSJ decomposition} 
due to Jaco and Shalen \cite{Jaco-Shalen} and Johannson \cite{Johannson}.
(See \cite{Neumann-Swarup} and \cite{Costantino2002} for an alternative proof,
and see \cite{Hatcher} for a simple proof of the torus decomposition theorem.)

\begin{theorem}
[Torus decomposition theorem]\index{torus decomposition theorem}
\label{thm:JSJ-mfd}
For a compact orientable irreducible $3$-manifold $M$,
there is a unique (up to isotopy) family $\mathcal{T}$
of disjoint essential tori,
satisfying the following properties.
\begin{enumerate}[(a)]
\item
Each closed up component of $M-\mathcal{T}$
is either a Seifert fibered space
or atoroidal.
\item
If any component of $\mathcal{T}$ is deleted, 
Property (a) fails.
\end{enumerate}
\end{theorem}

In the above theorem,
by a {\it closed up component} $M-\mathcal{T}$,
we mean the closure of a component of 
$M-N(\mathcal{T})$,
where $N(\mathcal{T})$ is a regular neigborhood of $\mathcal{T}$.
The subsurface $\mathcal{T}$ is called 
the {\it characteristic toric family}\index{characteristic toric family} of $M$,
and each closed up component of $M-\mathcal{T}$
is called a {\it JSJ piece}\index{JSJ piece} of $M$.

It should be noted that the family $\mathcal{T}$ is not only unique up to homeomorphism
but also unique up to {\it isotopy}.
This forms a sharp contrast to the fact that in the prime decomposition theorem,
the family of the splitting spheres is not unique even up to homeomorphisms.
(It only says that the resulting prime manifolds are unique up to homeomorphisms.)

\subsection{The Geometrization Conjecture of Thurston}
\label{subsection:geometrization-theorem}

Thurston's geometrization conjecture says that
any compact orientable irreducible $3$-manifold has a canonical splitting,
by tori, into pieces
which admit one of the following eight homogeneous geometries.
\begin{itemize}
\item[$\circ$]
The spaces of constant curvature, $\Sphere^3$, $\EE^3$ and $\HH^3$;
\item[$\circ$]
The product spaces $\Sphere^2\times \EE^1$ and $\HH^2\times \EE^1$; and
\item[$\circ$]
The $3$-dimensional Lie groups $\mathrm{Nil}$, $\mathrm{Sol}$, and $\tilde{\mathrm{SL}}_2(\RR)$.
\end{itemize}
Here a compact connected orientable $3$-manifold $M$ is {\em geometric}\index{geometric}
if either it is a $3$-ball or its interior can be presented as the quotient 
$\interior M=X/\Gamma$
of one of the above homogeneous spaces, $X$,
by a discrete group $\Gamma$ of isometries
acting freely and discontinuously on $X$.
If $X=\mathrm{Sol}$ then $M$ is a bundle over $S^1$ 
or the $1$-dimensional orbifold $S^1/(z\sim \bar z)$ with torus fiber;
if $X$ is neither $\mathrm{Sol}$ nor $\HH^3$,
then $M$ is a Seifert fibered space
and it is completely described by Seifert invariants.
Conversely any Seifert fibered space admits one of the 6 remaining geometris.
See the nice expositions \cite{Bonahon, Scott1983} for details.

Thus 
we have a complete topological classification
of the geometric manifolds with $X$ geometry for $X\ne \HH^3$,
and the study of hyperbolic manifolds forms the crucial part 
in $3$-manifold theory and knot theory.

Thurston proposed the following geometrization theorem as a conjecture,
which says that the torus decomposition gives a complete obstruction 
for a compact orientable $3$-manifold to be hyperbolic.
Here $M$ is said to be {\it hyperbolic}\index{hyperbolic} if its interior 
can be presented as the quotient 
$\interior M=\HH^3/\Gamma$
by a discrete torsion-free group $\Gamma$ of isometries of $\HH^3$
(cf. Section \ref{sec:hyperbolic-mfd}).
Thurston proved the conjecture for various cases,
including the case when $M$ is Haken, and 
the whole conjecture was finally proved by Perelman 
\cite{Perelman2002, Perelman2003a, Perelman2003b}.
(See \cite{Morgan2004, McMullen2011} for a survey,
and see \cite{CZ2006, Kleiner-Lott_2008, Morgan-Tian2007, Morgan-Tian2014, BBMBP} 
for detailed expositions.)

\begin{theorem}
[Geometrization theorem]\index{geometrization theorem}
\label{thm:geometrization-mfd}
Let $M$ be a connected irreducible atoroidal compact orientable $3$-manifold.
Then $M$ is either a Seifert fibered space or a hyperbolic manifold.
\end{theorem}

By combining Theorems \ref{thm:JSJ-mfd} and \ref{thm:geometrization-mfd},
we obtain the following geometric decomposition theorem.

\begin{theorem}
[Geometric decomposition theorem]\index{geometric decomposition}
\label{thm:geometric-decomposition-mfd}
For a compact orientable irreducible $3$-manifold $M$,
there is a unique (up to isotopy) family $\mathcal{T}$
of disjoint essential tori
satisfying the following properties.
\begin{enumerate}[(a)]
\item
Each closed up component of $M-\mathcal{T}$
is either a Seifert fibered space
or a hyperbolic manifold.
\item
If any component of $\mathcal{T}$ is deleted, 
Property (a) fails.
\end{enumerate}
\end{theorem}

In the above theorem,
a closed up component of $M-\mathcal{T}$
is called a {\it Seifert piece}\index{Seifert piece} or 
a {\it hyperbolic piece}\index{hyperbolic piece}
according to whether it is a Seifert fibered space
or a hyperbolic manifold.

\subsection{Geometric decompositions of knot exteriors}
We describe a consequence for knot exteriors
of the torus decomposition theorem and the geometrization theorem 
described in the previous subsection.
Let $K$ be a knot and consider its exterior $E(K)$.
Then $E(K)$ is irreducible by 
the Sch\"onflies theorem.
Moreover, $E(K)$ is Haken,
because a minimal genus Seifert surface is an incompressible surface in $E(K)$.

\begin{theorem}
[The geometric decomposition of knot exteriors]\index{geometric decomposition}
\label{thm:JSJ-Knot1}
Given a knot $K$ in $S^3$,
there is a unique (up to isotopy)
compact subsurface $\mathcal{T}$ in the interior $E(K)$
satisfying the following properties.
\begin{enumerate}[(a)]
\item
Each component of $\mathcal{T}$ is an essential torus.
\item
Each closed up component of $E(K)-\mathcal{T}$
is either a Seifert fibered space
or a hyperbolic manifold of finite volume.
\item
If any component of $\mathcal{T}$ is deleted, 
Property (b) fails.
\end{enumerate}
\end{theorem}

\begin{figure}[ht]
\begin{center}
 {
\includegraphics[height=5cm]{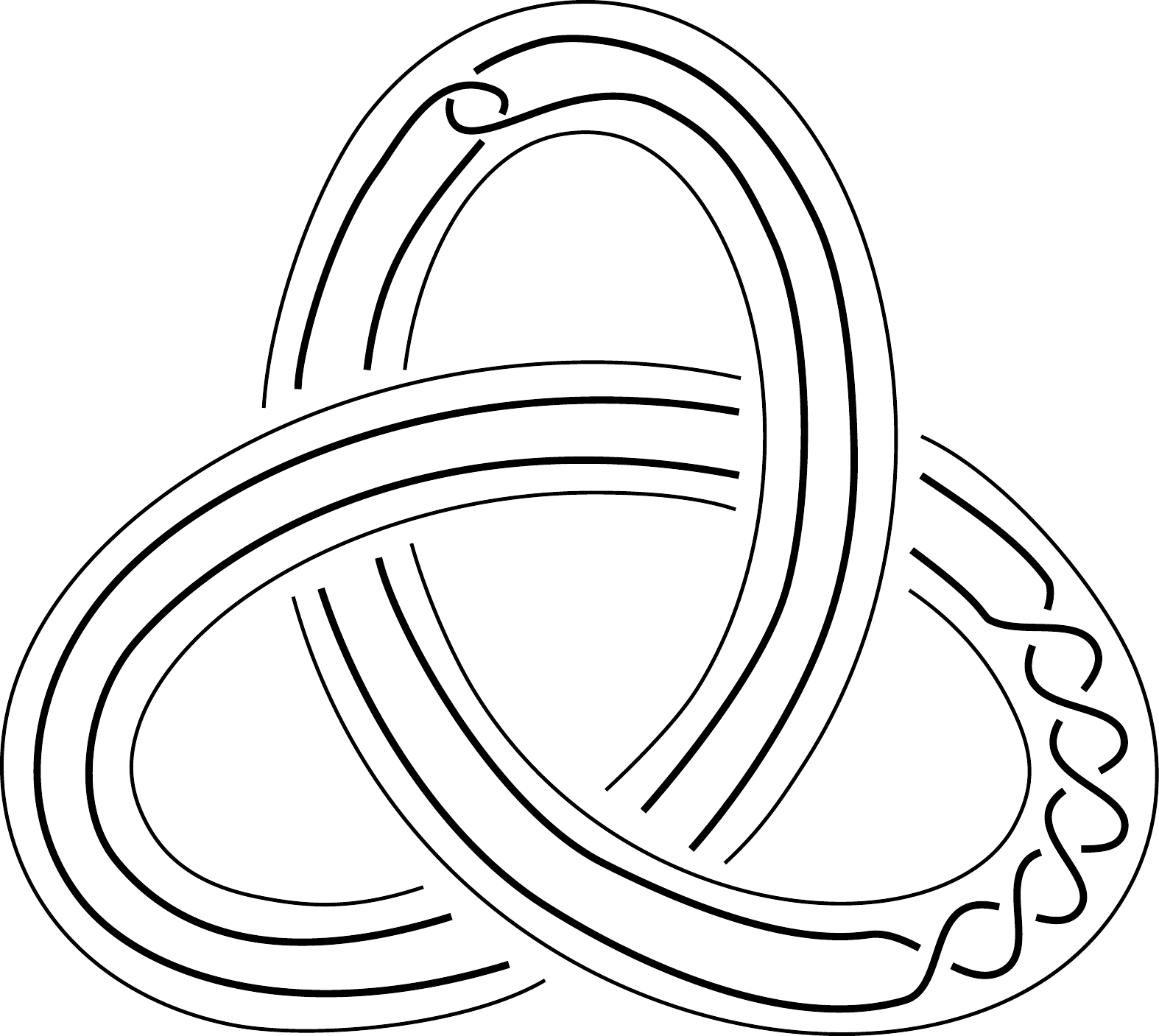}
 }
\end{center}
\caption
{
A {\it Whitehead double} of a trefoil knot is a satellite knot whose companion knot is 
a trefoil knot
and whose pattern knot is represented by the Whitehead link.
}
\label{fig:satellite-knot}
\end{figure}

We call the JSJ piece of $E(K)$ containing $\partial E(K)$ the {\it root JSJ piece}
\index{root JSJ piece}.
The JSJ-decomposition is intimately related 
to Schubert's satellite operation \cite{Schubert_satellite}.
To see this, 
assume that $\mathcal{T}\ne\emptyset$
and pick a component $T$ of $\mathcal{T}$.
Then $T$ bounds a solid torus $V=S^1\times D^2$ in $S^3$, 
which satisfies the following conditions.
\begin{enumerate}
\item
The core $k:=S^1\times 0$ of $V$ forms a nontrivial knot in $S^3$.
\item
$K$ is contained in $V$ geometrically essentially, 
i.e., there is no $3$-ball $B$ such that $K\subset B\subset V$.
Moreover, $K$ is not isotopic in $V$ to the core $k$ of $V$.
\end{enumerate}
Thus $K$ is a {\it satellite knot}\index{satellite knot}\index{knot!satellite}
 of the {\it companion knot}\index{companion knot}\index{knot!companion} 
 $k\subset S^3$
with {\it pattern}\index{pattern} $(V,K)$ (cf. 
Figure \ref{fig:satellite-knot} and \cite[Section 4.D]{Rolfsen}).

It should be noted that composition of knots is a special case of the satellite operation.
In fact, the composite knot $K_1\#K_2$ is a satellite of $K_1$ with pattern
$(S^3-\interior N(\mu_2),K_2)$ where $\mu_2$ is a meridian of $K_2$.
It is also a satellite of $K_2$ with pattern $(S^3-\interior N(\mu_1),K_1)$.

\medskip

It turns out that JSJ pieces of knot exteriors are expressed as link exteriors
(Theorem \ref{thm:JSJ-Knot2}).
A {\it link}\index{link} $L$ is a smoothly (or piecewise-linearly) embedded
disjoint union of circles in $S^3$,
namely $L$ is a closed $1$-submanifold of $S^3$.
Thus a knot is a link of $1$ component.
A link of $\mu$-components is called a 
{\it $\mu$-component trivial link}
if it bounds $\mu$ disjoint disks in $S^3$,
and we denote it by $O_{\mu}$. 
The {\it exterior} of a link $L$ is defined by
$E(L):=S^3-\interior N(L)$,
where $N(L)$ is a regular neighborhood of $L$.
The links in the example below play a key role in torus decompositions of knot exteriors.

\begin{example}
{\rm 
(1)
The $\mu+1$-component 
{\it key chain link}\index{key chain link}
$H_{\mu+1}=K_0\cup O_{\mu}$ 
is a union of the $\mu$-component trivial link $O_{\mu}$
and the trivial knot $K_0$ which intersects each member of 
$\mu$ disjoint links bounded by $O_{\mu}$.
Then $E(H_{\mu+1})$ is homeomorphic to $\mbox{($\mu$ holed disk)}\times S^1$,
and is called a {\it composing space}\index{composing space}.
If $E(H_{\mu+1})$ is the root JSJ piece of a knot exterior $E(K)$,
then $K$ is a connected sum of $\mu$ prime knots.

(2)
For a pair $(p,q)$ of relatively prime integers,
the {\em $(p,q)$-torus knot}\index{torus knot}
$K_{p,q}$ is defined by 
\[
K_{p,q}:=\{(z_1,z_2)\in S^3 \ | \ z_1^p+z_2^q = 0\}.
\]
$K_{p,q}$ is a regular orbit
of the circle action on $S^3$ given by
\[
\omega\cdot (z_1,z_2)=(\omega^q z_1,\omega^p z_2)
\quad (\omega\in S^1\subset \CC).
\]
Thus $E(K_{p,q})$ is a Seifert fibered space.
$K_{p,q}$ is contained in the torus
\[
T:=\{(z_1,z_2)\in S^3 \ |\ |z_1|^p=|z_2|^q\},
\]
and it wraps $q$ times in the $z_1$ direction and $p$ times in the $z_2$ direction.
The annulus $A:=T\cap E(K)$ divides $E(K_{p,q})$ into two solid tori.
By van-Kampen's theorem in this setting, we see that
\[
G(K_{p,q})=\langle a,b \ | \ a^p=b^q \rangle.
\]
The cyclic subgroup generated by $a^p=b^q$
forms the infinite cyclic center of $G(K_{p,q})$.
Moreover, a knot $K$ is a torus knot if and only if
$G(K)$ has a nontrivial center.
$K_{p,q}$ is nontrivial if and only if
both $p$ and $q$ have absolute value $\ge 2$.
If $E(K_{p,q})$ is a JSJ piece of a knot exterior $E(K)$,
then $K$ is a satellite of the torus knot $K_{p,q}$.

(3)
For a pair $(p,q)$ of relatively prime integers with $p\ge 2$,
the {\it $(p,q)$-Seifert link}\index{Seifert link} $C_{p,q}$ is defined by
\[
C_{p,q}:=K_0\cup K_{p,q} \quad \mbox{with $K_0=\{(z_1,z_2)\in S^3 \ | \ z_2 = 0\}$}.
\]
If $C_{p,q}$ is the root JSJ piece of a knot exterior $E(K)$,
then $K$ is the {\it $(p,q)$-cable}\index{cable}
of some nontrivial knot.
}
\end{example}

We have the following characterization of the torus decompositions
of knot exteriors.

\begin{theorem}
\label{thm:JSJ-Knot2}
{\rm (1)}
A compact orientable $3$-manifold $M$ is a JSJ piece of $E(K)$
for some nontrivial knot $K$ in $S^3$
if and only if $M\cong E(L)$ for some link $L$ in $S^3$,
which is the union of a knot $K_0$
and a trivial link $O_{\mu}$ (with $\mu$ possibly $0$),
such that $E(L)$ is either 
(i) hyperbolic or (ii) a Seifert fibered space homeomorphic to
a composing space, a nontrivial torus knot exterior, or a cable space.

{\rm (2)}
Let $K$ be a nontrivial knot in $S^3$, and 
let $\mathcal{T}$ be a union of disjoint essential tori in
$E(K)$, satisfying the following conditions.
\begin{enumerate}
\item[{\rm (i)}]
Each closed up component of $E(K)-\mathcal{T}$ is homeomorphic to 
a link exterior $E(L)$ which satisfies the condition in {\rm (1)}.
\item[{\rm (ii)}]
There does not exist a pair of 
adjacent closed up components of $E(K)-\mathcal{T}$, 
both of which are composing spaces.
\end{enumerate}
Then $\mathcal{T}$ is the characteristic toric family of $E(K)$.
\end{theorem}

The way JSJ pieces fit together in $E(K)$ is recorded by
the {\it companionship tree}\index{companionship tree}, defined as follows:
The vertices correspond to the JSJ pieces,
and the edges correspond to the components of $\mathcal{T}$,
where if an edge corresponds to a component $T$ of $\mathcal{T}$, it
joins the vertices corresponding to the two JSJ pieces
containing $T$ as a boundary component.
Since $H_1(S^3)=0$, this graph is a tree.
For a more detailed description of torus decompositions, 
see \cite[Chapter 2]{Bonahon-Siebenmann} and \cite{Budney}.

\section{The orbifold theorem and the Bonahon--Siebenmann decomposition of links}
\label{sec:Bonahon-Siebenmann_orbifold-thm}

In \cite[Chapter 13]{Thurston0},
Thurston initiated the systematic study of orbifolds,
namely quotients of spaces 
by properly discontinuous group actions
which are not necessarily free.
In 1978, he announced the orbifold theorem,
the geometrization theorem of $3$-orbifolds
which have non-empty $1$-dimensional singular set.
Every link $L=\cup_j K_j$ determines an infinite family of orbifolds,
by regarding each component $K_j$ as the singular locus
of cone angle $2\pi/n_j$ for some $n_j\ge 2$. 
The case when $n_j=2$ for every $j$ is particularly important,
and the Bonahon--Siebenmann decomposition theory of links 
is essentially the decomposition theory of such orbifolds.
Their theory is intimately related with
Conway's ingeneous analysis of link diagrams,
and gives us a nice method for understanding links directly from their diagrams.
In particular, it gives a complete classification
of the \lq\lq algebraic links'',
which implies, for example,
that the Kinoshita--Terasaka knot and the Conway knot are different
and that they admit no symmetry.

The purpose of this section is to recall the orbifold theorem
and its impact on knot theory.
To be precise, we will give surveys of
(i) the Bonahon--Siebenmann decomposition theory,
(ii) the classification of $2$-bridge links,
(iii) the orbifold theorem,
and (iv) application of the orbifold theorem to the study of branched cyclic coverings.

\subsection{The Bonahon--Siebenmann decompositions for simple links}
\label{subsec:Bonahon-Siebenmann}
By the geometric decomposition Theorems \ref{thm:JSJ-Knot1} and \ref{thm:JSJ-Knot2}
of knot exteriors,
the classification of knots is reduced to that of the links 
whose exteriors have trivial torus decompositions.
Deriving from Montesinos' work \cite{Montesionos1, Montesionos1b}
on double branched coverings of links
and Thurston's work on $3$-dimensional orbifolds,
Bonahon and Siebenmann 
established a new decomposition theorem for such links.
This is essentially a $\ZZ/2\ZZ$-equivariant JSJ decomposition theory,
applied to the double branched coverings of links.

To explain their results, we introduce a few definitions.
A link $L$ in $S^3$ is {\it splittable}\index{splittable} if there is a $2$-sphere $S$ in $S^3$
which separates the components of $L$.
$L$ is {\it unsplittable}\index{unsplittable} if it is not splittable.
This is equivalent to the condition that $E(L)$ is irreducible.
$L$ is {\it simple for Schubert}\index{simple for Schubert} if 
$E(L)$ is irreducible and atoroidal.
If $L$ is simple for Schubert, then
the JSJ decomposition of $E(L)$ is trivial.
The converse also holds for knots, but not
for links.
For example, the key-chain link $H_{\mu+1}$ is not simple for Schubert,
but the torus decomposition of $E(H_{\mu+1})$ is trivial.

Let $(M,L)$ be a pair consisting of a compact orientable $3$-manifold
and a proper $1$-submanifold $L$ in $M$.
A {\it Conway sphere}\index{Conway sphere} in $(M,L)$ is a $2$-sphere $\Sigma$ in $\interior M$
or $\partial M$ which meets $L$ transversely in $4$ points.
A Conway sphere $\Sigma$ is said to be {\it pairwise compressible}\index{pairwise compressible} 
if there is a disk $D$ in $M-L$
such that $D\cap \Sigma=\partial D$ does not bound a disk in $\Sigma-L$.
Otherwise, $\Sigma$ is said to be {\it pairwise incompressible}.
Two Conway spheres $\Sigma$ and $\Sigma'$ in $M$ are said to be {\em pairwise parallel},
if there is a closed up component $N$ of $M-(\Sigma\cup \Sigma')$
bounded by $\Sigma$ and $\Sigma'$ such that 
$(N, N\cap L)\cong (\Sigma, \Sigma\cap L)\times [0,1]$.
A Conway sphere is {\it essential}\index{essential} if it is pairwise incompressible
and is not pairwise parallel to a boundary component.
$(M,L)$ is {\em simple for Conway}\index{simple for Conway} if it does not contain an essential Conway sphere.

\begin{figure}[ht]
\begin{center}
 {
\includegraphics[height=4cm]{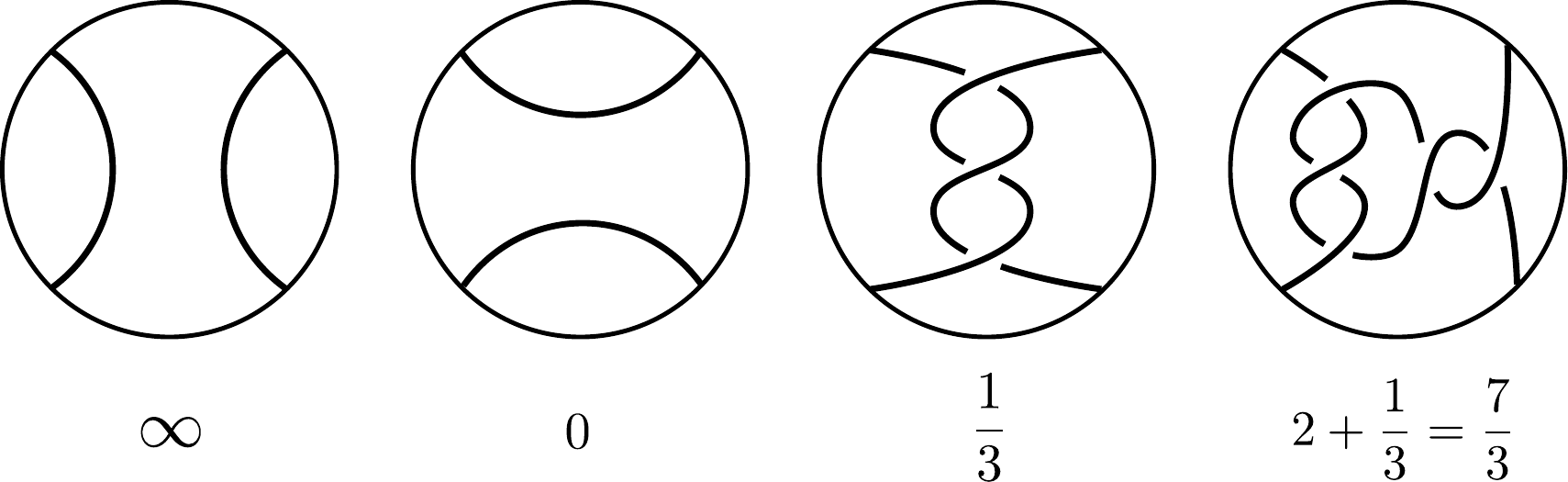}
 }
\end{center}
\caption
{
The pair of arcs forming a trivial tangle is parallel to 
a pair of arcs on the boundary Conway sphere.
If we identify the Conway sphere with the quotient of $\RR^2$ by the group generated by 
$\pi$-rotations around the lattice points,
then the inverse image of the pair of arcs in $\RR^2$ 
forms a family of mutually disjoint lines 
of rational slope passing through the lattice points.
This slope is the {\it slope} of the rational tangle.
}
\label{fig:rational-tangle}
\end{figure}

A {\it trivial tangle}\index{trivial tangle} is a pair $(B^3, t)$,
where $t$ is a union of two arcs properly embedded arcs in $B^3$
which is parallel to a pair of disjoint arc in $\partial B^3$.
A {\it rational tangle}\index{rational tangle} is a trivial tangle $(B^3, t)$
which is endowed with an identification of $\partial (B^3,t)$ with the
Conway sphere standardly 
embedded in $\RR^3\subset S^3$.
A rational tangle, up to the natural equivalence relation, 
is determined uniquely by its {\it slope} as illustrated in Figure \ref{fig:rational-tangle}.
(See \cite{Conway} for the original definition, and see
\cite[Chapter 8.1]{Bonahon-Siebenmann} or 
\cite[Section 8.6]{Cromwell2004} for detailed exposition.)

\begin{figure}[ht]
\begin{center}
 {
\includegraphics[height=3cm]{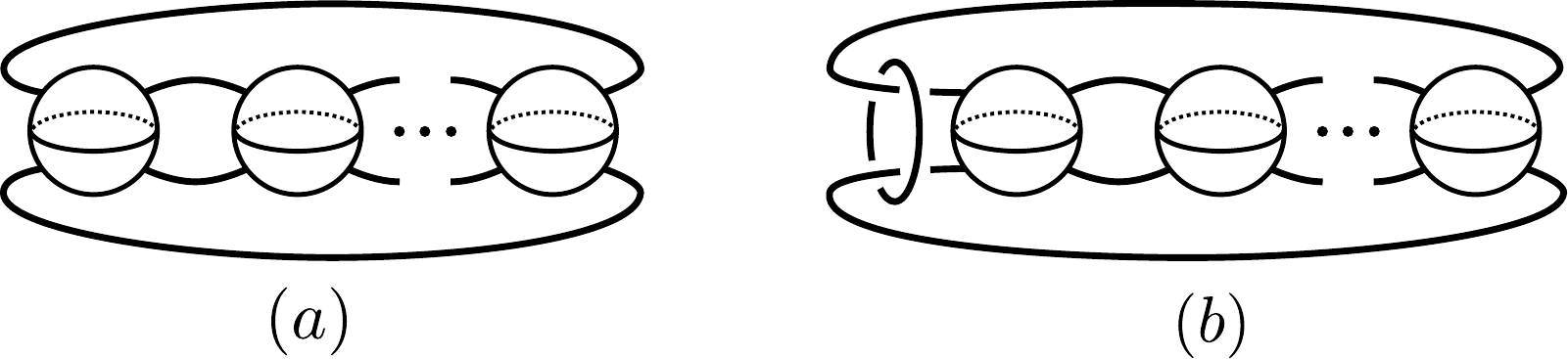}
 }
\end{center}
\caption
{
(a) a hollow Montesinos pair, (b) a hollow Montesinos pair with a ring
}
\label{fig:Montesinos-pair}
\end{figure}

A {\it Montesinos pair}\index{Montesinos pair} is a pair $(M,L)$ which is built from 
a {\it hollow Montesinos pair} or a {\it hollow Montesinos pair with a ring}
in Figure \ref{fig:Montesinos-pair} by plugging some of the holes 
with rational tangles of finite slope.

Bonahon and Siebenmann established the following decomposition theorem
\cite[Theorem 3.4]{Bonahon-Siebenmann}.

\begin{theorem}
\label{thm:BS-decomosition}
For a link $L$ in $S^3$ that is simple for Schubert,
there is a unique (up to isotopy respecting $L$)
compact subsurface $\mathcal{G}\subset S^3$
satisfying the following property.
\begin{enumerate}[(a)]
\item
Each component of $\mathcal{G}$ is a pairwise incompressible Conway sphere.
\item
Each closed up component $N$ of $S^3-\mathcal{G}$ gives a pair $(N,N\cap L)$
that either is simple for Conway, or else is a Montesinos pair.
\item
If any component of $\mathcal{G}$ is deleted, 
Property (b) fails.
\end{enumerate}
\end{theorem}

The above decomposition is called the 
{\it characteristic decomposition}\index{characteristic decomposition}
(or the {\it Bonahon--Siebenmann decomposition}\index{Bonahon--Siebenmann decomposition})
of $(S^3,L)$.
The union of the Montesinos pairs is called the 
{\it algebraic part}\index{algebraic part} or 
{\it arborescent part}\index{arborescent part} of $(S^3,L)$.
The link $(S^3,L)$ is called an 
{\it arborescent link}\index{arborescent link}
if its arborescent part is equal to the whole pair $(S^3,L)$.
This terminology comes from the fact that
arborescent parts can be represented by weight planar trees.
The classification the arborescent parts and links
is given by \cite[Part V]{Bonahon-Siebenmann}.

\begin{example}
{\rm
The Bonahon--Siebenmann decomposition of the
Kinoshita--Terasaka knot and the Conway knot
are given by the spheres in Figure \ref{fig:KT-Conway}
(cf. \cite[Sections 4.1 and 4.2]{Montesionos1b}).
This fact gives an intuitive proof of the inequivalence of these two knots.
It also shows that both knots are arborescent.
}
\end{example}

For a link $L$ in $S^3$,
let $p:M_2(L)\to S^3$ be the double branched covering of $S^3$ branched over $L$,
and let $\tau$ be the covering involution. 
Then the Bonahon--Siebenmann decomposition of $L$ can be regarded as a $\ZZ_2$-equivariant version
of the torus decomposition of $M_2(L)$
for the following reasons:

\begin{itemize}
\item[$\circ$]
The inverse image of an essential Conway sphere of $(S^3,L)$ 
is an essential torus in $M_2(L)$.
\item[$\circ$]
Let $(N,N\cap L)$ be a piece of the Bonahon--Siebenmann decomposition of $(S^3,L)$
which is a Montesinos pair.
Then the inverse image $p^{-1}(N)$ 
 is a Seifert fibered space,
where the base orbifold is orientable or non-orientable according to whether
$(N,N\cap L)$ is obtained from 
a hollow Montesinos pair or 
that with a ring
(see \cite{Montesionos1, Montesionos2}).
Moreover, the covering involution $\tau$ preserves the Seifert fibration of $p^{-1}(N)$.
The image of its fiber in $S^3$
is either a circle disjoint from $L$
or an interval with endpoints in $L$.
\item[$\circ$]
The above fact implies that the inverse image in $M_2(L)$ of the arborescent part of 
$(S^3,L)$ is a graph manifold (cf. Waldhausen \cite{Waldhausen_graph}).
In particular,
if $(S^3,L)$ is an arborescent link
then $M_2(L)$ is a graph manifold.
\item[$\circ$]
For each piece $(N,N\cap L)$ of the Bonahon--Siebenmann decomposition of $(S^3,L)$
which is not a Montesinos pair,
the inverse image $p^{-1}(N)$ is 
irreducible and atoridal.
Moreover, by the orbifold theorem (Theorem \ref{orbifold-theorem}) 
explained later in this section,
$p^{-1}(N)$ admits 
a complete hyperbolic structure of finite volume,
with respect to which $\tau|_{p^{-1}(N)}$ is an isometry.
\end{itemize}

We note that
the Bonahan-Siebenmann decomposition 
is intimately related with Conway's ingenious analysis of knot diagrams,
which in turn is based on Kirkmann's idea from the 19th century
(see \cite{HTM} for the history).
In fact, it reveals that Conway's notation for a link diagram
is not merely a convenient tool for describing diagrams
but also contains geometric information of the link represented by a diagram.
This is certainly the case for algebraic parts of the link.
As shown in \cite[Theorems 1.4 and 6.11]{Bonahon-Siebenmann},
Conway's notation for non-algebraic parts also 
has geometric information under certain conditions.

\subsection{$2$-bridge links}

In this subsection, 
we introduce $2$-bridge links,
which form a very special but important class of links.
For a rational number $r\in\QQ\cup\{1/0\}$,
the {\it $2$-bridge link, $K(r)$, of slope $r$}\index{$2$-bridge link}
is defined as the \lq\lq sum''
of the rational tangles of slopes $r$ and $1/0$.
To be precise, it is obtained from the rational tangles, $(B^3, t(r))$
and $(B^3, t(1/0))$, of slopes $r$ and $1/0$, respectively,
by glueing $(B^3, t(r))$ and $(-B^3, t(1/0))$
along the boundaries via the identity map.
(Note that the boundaries of rational tangles are 
identified with the Conway sphere standardly embedded in $\RR^3$.)
If $r=q/p$ where $p\ge 0$ and $q$ are relatively prime integers,
then $K(r)$ is a knot or a two-component link according to whether
$p$ is odd or even.
The following classification theorem was proved by
Schubert \cite{Schubert1956}, by establishing the uniqueness up to isotopy of 
{\it $2$-bridge spheres}\index{$2$-bridge sphere}
($2$-spheres which divide $K(r)$ into two trivial tangles). 
\begin{itemize}
\item[$\circ$]
Two $2$-bridge links $K(q/p)$ and $K(q'/p')$ are isotopic
if and only if $p=p'$ and either $q\equiv q' \pmod{p}$
or $qq'\equiv 1 \pmod{p}$.
They are homeomorphic if and only if
$p=p'$ and either $q\equiv \pm q' \pmod{p}$
or $qq'\equiv \pm1 \pmod{p}$.
\end{itemize}
The double branched covering, $M_2(K(q/p))$, of $S^3$
branched over $K(q/p)$ is the lens space $L(p,q)$,
and the above classification of $2$-bridge links
can be also deduced from the classification of lens spaces,
which in turn was established by Reidemeister \cite{Reidemeister},
using the Reidemeister torsion.
Moreover, the following characterization of $2$-bridge link
was obtained by Hodgson and Rubinstein \cite{Hodgson-Rubinstein-1985},
by classifying involutions on lens spaces with $1$-dimensional fixed point sets.

\begin{itemize}
\item[$\circ$]
A link $L$ in $S^3$ is a $2$-bridge link
if and only if the double branched covering $M_2(L)$
is a lens space.
\end{itemize}

\noindent
The result of \cite{Hodgson-Rubinstein-1985} is a special but important case
of the orbifold theorem (Theorem \ref{orbifold-theorem}) explained later in this section. 
They also proved the uniqueness up to isotopy of genus $1$ Heegaard surfaces
of lens spaces, which in turn gives a purely topological proof 
of the classification of lens spaces.

Thurston's uniformization theorem for Haken manifold 
(cf. Theorems \ref{thm:geometrization-mfd}),
together with an analysis of incompressible surfaces in the exterior
of $2$-bridge links, imply 
the following 
(cf. \cite[p.102]{Riley1979}, \cite[Lemmas 4.4]{Kobayashi1984}, \cite{Hatcher-Thurston}).
\begin{itemize}
\item[$\circ$]
The $2$-bridge link $K(q/p)$ is hyperbolic
if and only if $q\not\equiv \pm 1 \pmod {p}$.
\end{itemize}

\subsection{Bonahon--Siebenmann decompositions and $\pi$-orbifolds}
\label{subsec:Bonahon-Siebenmann-decomposition}

For a link $L$ in $S^3$,
the pair $(S^3,L)$ is homeomorphic to the quotient
$(M_2(L),\fix(\tau))/\tau$,
where $M_2(L)$ is the double branched covering of $S^3$ branched over $L$
and $\tau$ is the covering involution.
This means that the good $3$-orbifold $\OO(L):=M_2(L)/\tau$
has $S^3$ as underlying space and $L$ as singular set,
and each component of the singular set $L$ has cone angle $\pi$.
The Bonahon--Siebenmann decomposition is regarded as the 
torus decomposition of this orbifold.

Recall that an {\it $n$-orbifold}\index{orbifold} is a metrizable topological space $\OO$ 
locally modeled on the quotient of 
$\RR^n$ by a finite subgroup $G$ of the orthogonal group $\mathrm{O}(n)$.
If a point $x\in\OO$ corresponds to the image of the origin of $\RR^n$,
then the finite group $G$ is called the {\it local group} at $x$, and is denoted by $G_x$.
If $G_x$ is trivial, $x$ is {\it regular}, otherwise $x$ is {\it singular}.
The {\it singular locus}\index{singular locus} is the subset, $\Sigma_{\OO}$, of $\OO$
consisting of the singular points. 
When $G_x$ is the cyclic group generated by a $2\pi/m$-rotation around 
the codimension $2$ subspace $\RR^{n-2}\times\{0\}$,
we say that the point $x$ (and the stratum of the singular set containing $x$)
has {\it cone angle}\index{cone angle} $2\pi/m$ or {\it index}\index{index} $m$.
  
A quotient space $\OO:=X/\Gamma$, where $X$ is a smooth $n$-manifold and $\Gamma$ is a smooth 
properly discontinuous action, is an $n$-orbifold,
and its singular locus is the image of the subspace of $X$
consisting of points with nontrivial stabilizer.
If $\Gamma$ is a finite group,
such an orbifold is called a 
{\it good orbifold}\index{good orbifold}.
The {\it orbifold fundamental group}
$\pi_1^{\mathrm{orb}}(\OO)$ of $\OO$
is defined as the group consisting of all lifts of $\Gamma$
to the universal covering space $\tilde X$ of $X$.
Thus we have the following exact sequence.
\[
1 \to \pi_1(X) \to \pi_1^{\mathrm{orb}}(\OO) \to \Gamma \to 1
\]

For a link $(S^3,L)$, the orbifold $\OO(L):=(M_2(L),\fix(\tau))/\tau$
is called the {\it $\pi$-orbifold}\index{$\pi$-orbifold} 
associated with $L$.
The orbifold fundamental group  
$\pi_1^{\mathrm{orb}}(\OO(L))$ is called the 
{\it $\pi$-orbifold group}\index{$\pi$-orbifold group} of $L$.
It is calculated from 
the link group $G(L)=\pi_1(S^3-L)$
and a set of meridians $\{\mu_1,\dots, \mu_m\}$ as follows:
\[
\pi_1^{\mathrm{orb}}(\OO(L)) = G(L)/\langle\langle \mu_1^2,\dots, \mu_m^2 \rangle\rangle.
\]
Here $m$ is the number of components of $L$,
and $\mu_j$ is a meridian of the $j$-th component of $L$.
By using the orbifold theorem explained in the next subsection,
Boileau and Zimmermann \cite{Boileau-Zimmermann1} proved the following theorem which
shows that $\pi_1^{\mathrm{orb}}(\OO(L))$ is a very strong invariant for links.

\begin{theorem}
\label{thm:pi-orbifold}
Let $L$ be a prime unsplittable link in $S^3$ such that $\pi_1^{\mathrm{orb}}(\OO(L))$ is infinite.
Then the following hold.

(1) For any link $L'$ in $S^3$, 
the pairs $(S^3,L)$ and $(S^3,L')$ are homeomorphic if and only if
their $\pi$-orbifold groups $\pi_1^{\mathrm{orb}}(\OO(L))$ and $\pi_1^{\mathrm{orb}}(\OO(L'))$ are isomorphic.

(2) The natural homomorphism from the symmetry group $\Sym(S^3,L)$ to 
the outer-automorphism group $\Out(\pi_1^{\mathrm{orb}}(\OO(L)))$ is an isomorphism.
\end{theorem} 

Here the {\it symmetry group} $\Sym(S^3,L)$ is the group of 
diffeomorphisms of the pair $(S^3,L)$ up to isotopy.
It should be noted that the problem of determining the symmetry group of a knot 
is a refinement of the fundamental problem of determining 
whether the knot is chiral/invertible. 

Using the above theorem,
Boileau and Zimmermann \cite{Boileau-Zimmermann2} determined the symmetry groups 
of all non-elliptic Montesinos links,
i.e., the Montesinos links with infinite $\pi$-orbifold groups.
(The symmetry groups of elliptic Montesinos links were determined by 
\cite{Sakuma2} using the orbifold theorem.)
This result may be regarded as a broad extension of Trotter's proof \cite{Trotter}
of non-invertibility of the pretzel knot $P(p,q,r)$
with $|p|, |q|, |r|$ distinct odd integers $\ge 3$.
Trotter's proof is based on the fact that $\pi_1^{\mathrm{orb}}(\OO(P(p,q,r)))$
is an extension of the hyperbolic triangular reflection group
\[
[p,q,r]=\langle x,y,z \ | \ x^2, y^2, z^2, (xy)^p, (yz)^q, (zx)^r \rangle
\]
by the infinite cyclic group, which in turn is a consequence
of the fact that the $\pi$-orbifold $\OO(P(p,q,r))$
is a {\it Seifert fibered orbifold}\index{Seifert fibered orbifold}
over the $2$-dimensional hyperbolic orbifold
$\HH^2/[p,q,r]$.

The symmetry groups of the arborescent links are 
completely determined by Bonahon and Siebenmann in \cite{Bonahon-Siebenmann}.
In particular, this implies that the symmetry groups of the Kinoshita--Terasaka knot
and the Conway knot are trivial,
and so they are chiral and noninvertible.
The knot $8_{17}$ is also arborescent, and its symmetry group is 
the order $2$ cyclic group, generated by an orientation-reversing involution
representing the negative-amphicheirality of the knot.
This is the Bonahon-Siebenmann's proof of the non-invertibility of $8_{17}$
(cf. Subsection \ref{subsec:Alexander-polynomial}).

\subsection{The orbifold theorem and the Smith conjecture}
\label{subsec:orbifold-theorem}

Many of the concepts for $3$-manifolds,
such as irreduciblity, atoroidality and Seifert fibrations,
have natural generalization for $3$-orbifolds,
and a characteristic splitting (torus decomposition) theorem was  
establised by Bonahon and Siebenmann \cite{Bonahon-Siebenmann0}
(cf. \cite{BMP, Bonahon}). 
The characteristic splitting Theorem \ref{thm:BS-decomosition} for links
is a special case of the general splitting theorem,
though the detailed analysis for the algebraic parts and 
the application to knot theory in \cite{Bonahon-Siebenmann} 
cannot be found in \cite{Bonahon-Siebenmann0}.

Bonahon-Siebenmann's theory forms the first step towards the proof of the following geometrization theorem
for orbifolds, which was announced by Thurston \cite{Thurston1},
and finally proved by Boileau, Leeb and Porti \cite{BLP} 
(see also Cooper-Hodgson-Kerckhoff \cite{CHK2000} and
Boileau-Porti \cite{Boileau-Porti} for an earlier account,
and Dinkelbach-Leeb \cite{DL} for the generalization to non-orientable orbifolds
using the equivariant Ricci flow).

\begin{theorem}[Orbifold Theorem]\index{orbifold theorem}
\label{orbifold-theorem}
Every compact orientable good $3$-orbifold with nonempty singular set
has a canonical splitting 
by spherical $2$-dimensional suborbifolds and toric $2$-dimensional suborbifolds
into geometric $3$-orbifolds.
\end{theorem}

Here a $3$-orbifold $\OO$ is {\it geometric}\index{geometric}\index{orbifold!geometric}
if either it is the quotient of a ball by an orthogonal action,
or its interior has one of the eight Thurston geometries,
namely $\OO=X/\Gamma$,
where $X$ is one of the eight Thurston's geometries
and $\Gamma$ is a discrete subgroup of $\Isom(X)$.
(If $X$ is different from the constant curvature spaces $\HH^3$, $\EE^3$ and $\Sphere^3$,
then there is no canonical metric on $X$, however,
it admits a family of natural metrics for which $\Isom (X)$ are identical.
See the beautiful surveys \cite{Scott1983, Bonahon}.)

The orbifold theorem was first announced as the following symmetry theorem
concerning finite group actions on $3$-manifolds. 

\begin{theorem}[Symmetry Theorem]\index{symmetry theorem}
\label{thm:symmetry}
Let $M$ be a compact orientable irreducible $3$-manifold. 
Suppose $M$ admits an action by a finite group $G$ of orientation-preserving diffeomorphisms 
such that some non-trivial element has a fixed point set of dimension one. 
Then $M$ admits a geometric decomposition preserved by the group action.
\end{theorem}

This theorem poses a very strong restriction on finite group actions on knots
(see \cite{Boileau-Flapan1987, Luo, Sakuma1}).
In particular, 
it includes, as a special case, the following positive answer to the Smith conjecture.
 
\begin{theorem}[The Smith Conjecture]\index{Smith conjecture}
\label{thm:SmithConj}
If $h:S^3\to S^3$ is an orientation-preserving periodic diffeomorphism 
with non-empty fixed point set,
then $h$ is smoothly conjugate to an orthogonal diffeomorphism.
In particular, $\fix(h)$ is the trivial knot.
\end{theorem}

The proof of this conjecture recorded in \cite{Morgan-Bass}
may be regarded as
the first major impact of Thurston's unifromization theorem for Haken manifolds,
and it was established using the uniformization theorem,
the equivariant loop theorem by Meeks-Yau \cite{Meeks-Yau}, 
and a refinement of Bass-Serre theory \cite{Serre}.

In Theorems \ref{thm:symmetry} and \ref{thm:SmithConj},
the smoothness of the action is essential.
In fact there is an orientation-preserving periodic homeomorphism $h$ of $S^3$
which has a wild knot as the fixed point set;  in particular,
the cyclic action generated by $h$ is not topologically conjugate to an orthogonal action.
It is this phenomena that lead Shin'ichi Kinoshita and Hidetaka Terasaka, 
the founders of knot theory in Japan, into knot theory.
It is an amazing coincidence that Terasaka published 
an introductory book \cite{Terasaka} to non-Euclidean geometry for the general public
in 1977, around the time Thurston started the series of lectures
on the geometry and topology of $3$-manifolds.

\subsection{Branched cyclic coverings of knots}
\label{subsec:branched-covering}
In Subsection \ref{subsec:Bonahon-Siebenmann-decomposition},
we explained the important role of the double branched coverings
of knots and links.
Not only the double branched covering but also
the cyclic branched covering
has attracted keen attention 
of various mathematicians,
because the latter gives a bridge between knot theory and $3$-manifold theory
and because of its special beauty.
In this subsection, we review the impact of Thurston's work,
in particular the orbifold theorem, 
on the study of branched cyclic coverings of knots.

For a knot $K$ in $S^3$, let $M_n(K)$ be the $n$-fold cyclic branched covering
of $S^3$ branched over $K$.
We also call $M_n(K)$ the {\it $n$-fold cyclic branched covering of $K$}. 
Then we have the following natural question.

\begin{problem}
\label{prob:branched-covering}
{\rm
To what extent does the topological type of $M_n(K)$ determine $K$?
}
\end{problem}

It should be noted that
$M_n(K)$ inherits the orientation of the ambient space $S^3$,
but it is independent of the orientation of the circle $K$.
Namely $M_n(K)\cong M_n(-K)$ as oriented manifolds.
Thus the precise meaning of the above question is as follows.
To what extent does the topological type of the {\it oriented manifold} 
$M_n(K)$ determine the {\it isotopy type} of the {\it unoriented} knot $K$?
 
The positive solution of the Smith conjecture is essentially equivalent to
the following partial answer to the above problem
(see \cite{Morgan-Bass}).

\begin{theorem}[Branched covering theorem]
\label{thm:branched-covering-thm}
A knot $K$ in $S^3$ is trivial if and only if $M_n(K)\cong S^3$ for some $n\ge 2$. 
\end{theorem} 

The orbifold theorem gives a very strong tool for the study of 
Problem \ref{prob:branched-covering}.
Before describing its influence, let us recall two classical constructions
of pairs of knots sharing the same cyclic branched covering.

\medskip

{\bf Construction 1.}
Let $K$ be a non-invertible prime oriented knot.
Then, by the unique prime factorization theorem, 
the knots $K\# K$ and  $K\#(-K)$ are not isotopic
as unoriented knots.
However, they share the same $n$-fold cyclic branched covering for all $n\ge 2$,
because:
\[
M_n(K\# K) \cong M_n(K)\# M_n(K) \cong M_n(K)\# M_n(-K) \cong M_n(K\#(-K))
\]

\medskip

{\bf Construction 2.}
Let $L=K_1\cup K_2$ be a $2$-component link consisting of two trivial knots.
For integers $n_1, n_2\ge 2$ which are relatively prime to the linking number
$lk(K_1,K_2)$, the inverse image $\tilde K_1$ of $K_1$ in $M_{n_1}(K_2)\cong S^3$ 
is a knot, and so is
the inverse image $\tilde K_2$ of $K_2$ in $M_{n_2}(K_1)\cong S^3$.
Moreover, both $M_{n_2}(\tilde K_1)$ and $M_{n_1}(\tilde K_2)$ are homeomorphic to
the $(\ZZ/n_1\ZZ)\oplus (\ZZ/n_2\ZZ)$-covering of $S^3$ branched over $L$,
and hence they are homeomorphic.
(There is an analogous construction by using a three-component link
such that any $2$-component sublink is a Hopf link
(see \cite[0.2]{Reni-Zimmermann2001}).)

\medskip
Now, we state an important consequence of the orbifold theorem
(see \cite{CHK2000}).

\begin{theorem}
\label{thm:hyp-branched-cover}
Let $K$ be a hyperbolic knot in $S^3$,
i.e., $K$ is a knot which is neither a torus knot nor a satellite knot. 
Then $M_n(K)$ is hyperbolic for all $n\ge 3$, except for the $3$-fold covering of the figure eight knot (which is a Euclidean manifold).
Moreover, the covering transformation group acts on $M_n(K)$ by isometries. 
\end{theorem}

\begin{remark}
\label{rem:nono-pi-hyperbolic}
{\rm
In the above theorem the assumption $n\ge 3$ is essential.
In fact, if a hyperbolic knot contains an essential Conway sphere, $\Sigma$,
then the inverse image, $\tilde \Sigma$, of $\Sigma$ in $M_2(K)$ is an essential torus 
and hence $M_2(K)$ is non-hyperbolic even though $K$ itself is hyperbolic.
Moreover, every arborescent link
has a graph manifold as double branched covering.
}
\end{remark}

The hyperbolic Dehn surgery theorem implies that
if $n$ is sufficiently large, then the branch line forms the unique shortest closed geodesic in $M_n(K)$ (cf. Subsection \ref{subsec:Geometry-Dehn-filling}).
Using this fact, we can see that $M_n(K)$ for sufficiently large $n$ determines the knot $K$.
More generally, Kojima \cite{Kojima1986} proved the following theorem,
which gives a positive answer to a question of Goldsmith 
\cite[Problem 1.27]{Kirby1978}.

\begin{theorem}
For each prime knot $K$ there exists a constant $n_K$, 
such that two prime knots $K$ and $K'$ are equivalent 
if their $n$-fold cyclic branched covers are homeomorphic for some $n>\max(n_K,n_{K'})$.
\end{theorem}

We can reformulate Problem \ref{prob:branched-covering} as follows:
For a given connected closed orientable $3$-manifold $M$,
in how many different ways can $M$ occur as a cyclic branched covering
of a knot in $S^3$? 
There are two basic cases: the case when $M$ is a Seifert fibered space
and the case when $M$ is a hyperbolic manifold.
(The general case can be treated by using the equivariant sphere theorem
and torus decomposition \cite{Meeks-Scott1986} into
Seifert fibered space and hyperbolic manifolds.)

When $M$ is a Seifert fibered space,
the covering transformation group, $H$, is fiber-preserving
by \cite{Meeks-Scott1986} (when $\pi_1(M)$ is infinite)
and by the orbifold theorem (when $\pi_1(M)$ is finite).
If $H$ reverses the fiber-orientation,
then the quotient knot is a Montesinos knot
whereas if $H$ preserves the fiber-orientation
then the quotient knot is a torus knot.

In the case where $M$ is hyperbolic, we may assume, by the orbifold theorem,
that $H$ is a cyclic subgroup of
the finite group $\Isom^+(M)$.
The group $H$ must be a {\it hyper-elliptic group}\index{hyper-elliptic group},
namely $H$ is a finite cyclic group
such that $\fix h$ is a circle for every non-trivial element $h\in H$,
and $M/H\cong S^3$
(cf. \cite[Definition 1]{BMFPZ2018}).
Thus there is a one-to-one correspondence
\begin{align*}
\{\mbox{knots $K$ such that $M_n(K)\cong M$ for some $n\ge 2$}\}/\mbox{isotopy}\\
\leftrightarrow
\{\mbox{hyper-elliptic subgroups of $\Isom^+(M)$}\}/\mbox{conjugacy}.
\end{align*}
By Kojima's theorem \cite{Kojima1988},
any finite group can be the full isometry group of a closed orientable hyperbolic $3$-manifold.
However, the geometric condition 
for a hyper-elliptic group, $H$,
implies purely group theoretical conditions on $H$.
For example, we can see by using
the Smith conjecture (Theorem \ref{thm:SmithConj}) that
the normalizer of $H$ in $\Isom^+(M)$ is a finite subgroup of
the semi-direct product
$(\ZZ/2\ZZ)\ltimes(\QQ/\ZZ\oplus \QQ/\ZZ)$,
where $\ZZ/2\ZZ$ acts on $\QQ/\ZZ\oplus \QQ/\ZZ$ 
as multiplication by $-1$
(see \cite[Remark 3]{BMFPZ2018}).
Thus we have a chance to apply finite group theory to the study of cyclic branched coverings.
For example, if we are interested in the case when the degree $n$ is a prime number $p$,
then by Sylow's theorem,
every hyper-elliptic subgroup of order $p$
is conjugate to a cyclic sugroup of a single Sylow $p$-subgroup $S_p$ 
of $\Isom^+(M)$.
This interplay between the study of cyclic branched coverings
and finite group theory was initiated by Reni and Zimmermann,
and various interesting results were obtained,
including the following.

\begin{itemize}
\item[$\circ$]
Reni-Zimmermann \cite{Reni-Zimmermann2001}:
Let $K$ and $K'$ be two hyperbolic knots 
such that $M_n(K)\cong M_{n'}(K')$ for some $n, n'\ge 3$. 
Suppose further that $n$ and $n'$ have a common prime divisor $p>2$.
Then $K$ and $K'$ are related by Construction 2.
In particular, if $n=n'$ is not a power of $2$,
then the same conclusion holds (cf. \cite{Zimmermann1998}).
\item[$\circ$]
Paoluzzi \cite{Paoluzzi2005}:
A hyperbolic knot is determined by any three of its cyclic branched coverings
of order $\ge 2$.
Indeed, two coverings suffice if their orders are not coprime.
\item[$\circ$]
Boileau-Franchi-Mecchia-Paoluzzi-Zimmermann \cite{BMFPZ2018}:
A closed hyperbolic $3$-manifold is a cyclic branched covering of at most fifteen inequivalent knots in $S^3$.
\end{itemize}
A noteworthy aspect of the proof of the last result is the substantial use of finite group theory, in particular of the classification of finite simple groups. 

\medskip
For the double branched coverings,
we have the following additional construction.

\medskip
{\bf Construction 3.}
Let $\theta$ be a $\theta$-curve in $S^3$,
namely a spatial graph consiting of two vertices and three edges $\alpha_1$, $\alpha_2$ and $\alpha_3$,
each of which connects the two vertices.
For $\{i,j,k\}=\{1,2,3\}$,
suppose $A_{k}:=\alpha_i\cup\alpha_j$ forms a trivial knot.
Then the inverse image, $K_k$, of the arc $\alpha_k$ in $M_2(A_{k})\cong S^3$
forms a (strongly invertible) knot,
and $M_2(K_k)$ is identified with the $(\ZZ/2\ZZ)^2$-covering of $S^3$
branched over $\theta$.
If $A_{k}$ is a trivial knot for more than one $k\in \{1,2,3\}$,
we obtain knots in $S^3$ sharing the same double branched coverings.
(A similar construction is applied to embeddings of the $1$-skeleton
of the tetrahedron and the Kuratwski graph in $S^3$,
which produce potentially distinct $4$ and $9$ knots, respectively,
sharing the same double branched coverings
(see \cite{Mecchia-Reni2002}).)

\medskip

A link $L$ is said to be {\it $\pi$-hyperbolic}\index{$\pi$-hyperbolic}\index{link!$\pi$-hyperbolic} if $M_2(L)$ is hyperbolic.
For double coverings of $\pi$-hyperbolic knots, 
the following results were obtained.

\begin{itemize}
\item[$\circ$]
Boileau-Flapan \cite{Boileau-Flapan1995}:
If $K$ is a $\pi$-hyperbolic knot, then
every knot $K'$ which shares the same double branched covering with $K$
is constructed by repeatedly applying Constructions 2 and 3.
\item[$\circ$]
Reni \cite{Reni2000}:
There are at most nine different $\pi$-hyperbolic knots
with the same double branched coverings.
Mecchia-Reni \cite{Mecchia-Reni2002} gave a more geometric proof to this estimate
and proved that the same estimate holds for $\pi$-hyperbolic links.
\item[$\circ$]
Kawauchi \cite{Kawauchi2006}:
Reni's estimate is the best possible,
i.e., there are nine mutually inequivalent $\pi$-hyperbolic 
knots $K_i$ ($i=1,\cdots, 9$), in $S^3$
with the same double branched coverings. 
\end{itemize}
In the proof of the second result,
the study of the Sylow $2$-subgroup of $\Isom^+(M)$
of a closed orientable hyperbolic $3$-manifold holds a key.
The third result was obtained by using Kawauchi's imitation theory,
which yields, for a given $(3,1)$-manifold pair $(M,L)$,
a family of $(3,1)$-manifold pairs $(M^*,L^*)$
which is \lq\lq topologically similar'' to $(M,L)$.
A key example in the theory is 
the Kinoshita--Terasaka knot, which is an {\it imitation}\index{imitation} of the trivial knot.
(This fact was first found by Nakanishi \cite{Nakanishi1981}
and a beautiful generalization of this fact was given by Kanenobu \cite{Kanenobu1988}.)

\medskip
For the double branched covering of links
which are not $\pi$-hyperbolic,
we have the following additional construction.
(See Paoluzzi \cite{Paoluzzi2001} for further construction.)

\medskip
{\bf Construction 4.}
(Mutation)
Let $\Sigma$ be an essential Conway sphere of a link $L$ in $S^3$.
Cut $(S^3,L)$ along $\Sigma$ and reglue by an orientation-preserving involution 
of $(\Sigma,\Sigma\cap L)$ whose fixed point set is disjoint from $\Sigma\cap L$.
This process, called a {\it mutation}\index{mutation}, results in a new link $L'$ in $S^3$,
called a {\it mutant}\index{mutant} of $L$.
A pair of links are called {\it mutants}
if they are related by a sequence of mutations.
It was proved by Viro \cite[Theorem 1]{Viro} that
if $L$ and $L'$ are mutants then they share the same double branched coverings
(cf. \cite[Proposition 3.8.2]{Kawauchi1996}).

\medskip
In \cite{Greene2013}, Greene studied the Heegaard Floer homology of 
the double branched coverings of alternating links,
and proved that a reduced alternating link diagram is determined up to mutation by
the Heegaard Floer homology of the double branched covering of the link.
In particular, the following result follows.  

\begin{itemize}
\item[$\circ$]
Two reduced alternating links $L$ and $L'$ share the same double branched covering
if and only if $L$ and $L'$ are mutants.
\end{itemize}

He also proposes the mysterious conjecture:
{\it if a pair of links have homeomorphic double branched coverings, then either both are alternating or both are non-alternating}.

\section{Hyperbolic manifolds and the rigidity theorem}
\label{sec:hyperbolic-mfd}
In this section,
we recall basic facts concerning hyperbolic manifolds
and the Mostow-Prasad rigidity theorem
for hyperbolic manifolds of finite volume and of dimension $\ge 3$.
The rigidity theorem has had tremendous influence on knot theory,
because it guarantees that any geometric invariant
of the hyperbolic structure of a hyperbolic knot complement
is automatically a topological invariant of the knot complement.

For further information on hyperbolic geometry,
see the textbooks 
Benedetti-Petronio \cite{Benedetti-Petronio},
Ratcliffe \cite{Ratcliffe},
Matsuzaki-Taniguchi \cite{Matsuzaki-Taniguchi},
Anderson \cite{Anderson}
and
Marden \cite{Marden2}.

\subsection{Hyperbolic space}
\label{subsec:hyperbolic-space}
Let $\HH^n$ be the hyperbolic $n$-space,
i.e., the upper-half space 
\[
\HH^n:=\{(x_1,\dots, x_n)\in \RR^n \ | \ x_n>0\}
\]
in $\RR^n$ 
equipped with the Riemannian metric
\[
ds^2=\frac{1}{x_n^2}(dx_1^2+\cdots +dx_n^2).
\]
$\HH^n$ is the unique connected, simply connected, complete Riemannian manifold
of constant sectional curvature $-1$.
The isometry group $\Isom(\HH^n)$ is a real Lie group and
acts transitively on $\HH^n$,
and the stablizer of each point is identified with the 
orthogonal group $\Orth(n)$.
If $n\ge 3$, the ideal boundary $\partial \HH^n=(\RR^{n-1}\times\{0\})\cup\{\infty\}$
has a natural conformal structure,
and the orientation-preserving isometry group $\Isom^+(\HH^n)$
is identified with the group of conformal maps of $\partial \HH^n$.

Let $\gamma$ be a nontrivial element of $\Isom^+(\HH^n)$.
Then precisely one of the following holds.

\begin{enumerate}
\item
$\gamma$ is {\it elliptic}\index{elliptic}, i.e., $\gamma$ has a fixed point in $\HH^n$.
\item
$\gamma$ is {\it parabolic}\index{parabolic}, i.e.,
$\gamma$ has a unique fixed point, $x$, in $\partial \HH^n$,
called the {\it parabolic fixed point}\index{parabolic fixed point}.
Then $\gamma$ preserves every {\it horoball}\index{horoball}, $H_x$, centered at $x$.
Here, if $x\ne \infty$, then $H_x$ is the intersection of a (closed) Euclidean ball
with $\HH^n$ which touches $\partial\HH^n$ at $x$, and
if $x=\infty$ then $H_x$ is the closed upper-half space
\[
H_{\infty,c}:=\{(x_1,\dots, x_n)\in\HH^n \ | \ x_n\ge c\}
\quad \mbox{for some $c>0$,}
\]
called the horoball centered at $\infty$ with hight $c$.
The {\it horosphere}\index{horosphere} $\partial H_x$ inherits a Euclidean metric from the hyperbolic metric,
which is invariant by $\gamma$.
\item
$\gamma$ is {\it hyperbolic}\index{hyperbolic}, i.e.,
$\gamma$ has precisely two fixed points in $\partial \HH^n$,
one of which is repelling and the other is attracting.
The geodesic in $\HH^n$ joining the two fixed points 
is the unique geodesic which is preserved by $\gamma$;
it is called the {\it axis} of $\gamma$, and denoted by $\axis \gamma$.
\end{enumerate}

For low dimensions $n=2$ and $3$, we have:
\[
\Isom^+(\HH^2)\cong\PSL(2,\RR), \quad
\Isom^+(\HH^3)\cong\PSL(2,\CC).
\]
We identify the upper-half space $\HH^3=\RR^2\times \RR_+$ with $\CC\times \RR_+$ and 
identify the ideal boundary $\partial \HH^3$ with the Riemann sphere $\CC\cup\{\infty\}$.
Then the action of 
$A=\begin{pmatrix}
a & b\\
c & d
\end{pmatrix}
\in \PSL(2,\CC)$ on 
$\partial\HH^3=\CC\cup\{\infty\}$
is given by the linear fractional transformation
\[
A(z)=\frac{az+b}{cz+d}.
\]
Assume that $A\ne \pm E$, where $E$ is the identity matrix.
Then, as we see in the following, 
the orientation-preserving isometry of $\HH^3$
corresponding to
$A\in\PSL(2,\CC)$
is elliptic, parabolic, or hyperbolic
according as the trace
$\tr A$ (which is defined up to sign change) belongs to $(-2,2)$, $\{\pm 2\}$, or $\CC-[-2,2]$.

Case 1. $\tr A\ne \pm 2$.
Then $A$ has precisely two fixed points in $\partial\HH^3$.
After conjugation in $\PSL(2,\CC)$, we may assume that
they are $0$ and $\infty$.
Thus $A(z)=az$ for some $a=re^{\theta\sqrt{-1}}\in \CC^*-\{1\}$.
The action of the isometry $A$ on $\HH^3$ is given by:
\[
A(z,t)=(az, |a|t)=(re^{\theta\sqrt{-1}}z, rt).
\]
This is a skrew motion along the geodesic {\it axis} $0\times \RR_+$
with (signed) translation length $\log r$ and rotation angle $\theta$.
The quantity
\[
\LL_A:=\log r+\theta\sqrt{-1} =\log a \in \CC/2\pi \sqrt{-1}\ZZ
\]
is called the {\it complex translation length}\index{complex translation length} of the isometry $A$.
If we interchange $0$ and $\infty$ by conjugation in $\PSL(2,\CC)$,
then the complex translation length changes into $-(\log r+\theta\sqrt{-1})$.
Thus the complex translation length is defined only modulo $2\pi \sqrt{-1}\ZZ$
and up to multiplication by $\pm 1$.
In fact, a simple calculation implies
$\LL_A$ is characterized by the following identity:
(Note that $\tr A$ for $A\in\PSL(2,\CC)$ is defined only up to sign.)
\[
\pm \tr A=2\cosh \frac{\LL_A}{2}
\]
If we fix an orientation of the axis,
then the complex translation length is defined as an element in $\CC/2\pi \sqrt{-1}\ZZ$.
Note that $A$ is elliptic if and only if $r=1$,
which equivalent to the condition that $\tr A\in (-2,2)$
(under the assumption that $A\ne \pm E$).
Thus $A$ is hyperbolic or elliptic
according to whether $\tr A$ is contained in $(-2,2)$ or $\CC-[-2,2]$.

Case 2. $\tr A=\pm 2$.
Then $A$ has a unique fixed point in $\partial\HH^3$
and hence parabolic.
After conjugation in $\PSL(2,\CC)$, we may assume that
it is $\infty$.
Thus $A(z)=z+\tau$ for some $\tau\in \CC^*$.
The action of the isometry $A$ on $\HH^3$ is given by:
\[
A(z,t)=(z+\tau, t).
\]
In this case, the complex translation length $\LL_A$ is defined to be $0$.

\medskip

The following lemma can be easily proved. 

\begin{lemma}
\label{lem:centralizer}
Let $A$ be a nontrivial element in $\PSL(2,\CC)$ 
which is not an elliptic element of order $2$.
Then the centralizer $C(A)$ in $\PSL(2,\CC)$ 
is as follows.
\begin{enumerate}
\item 
If $A$ is elliptic or hyperbolic,
then $C(A)-\{E\}$ consists of elliptic/hyperbolic elements
which share the same axis with $A$.
Thus $C(A)$ is isomorphic to the multiplicative group $\CC^*$
\item
If $A$ is parabolic,
then $C(A)-\{E\}$ consists of parabolic elements
which share the same parabolic fixed point with $A$.
Thus $C(A)$ is isomorphic to the additive group $\CC$
\end{enumerate}
\end{lemma}

\subsection{Basic facts for hyperbolic manifolds}
By a {\em hyperbolic structure}\index{hyperbolic structure} on an $n$-manifold $M$,
we mean a Riemannian metric on $M$
of constant sectional curvature $-1$:
the curvature condition means that
every point in $M$ has a neighborhood isometric to an open set 
of $\HH^n$.
A hyperbolic structure on $M$ induces a hyperbolic structure 
on the universal cover $\tilde M$ of $M$
which is invariant by the action of the covering transformation group.
Thus we obtain a local isometry 
$D:\tilde M\to \HH^n$, called the {\em developing map}\index{developing map},
and a homomorphism $\rho:\pi_1(M)\to \Isom\HH^n$,
called the {\em holonomy representation}\index{holonomy representation}\index{representation!holonomy},
such that $D$ is $\rho$-equivariant,
i.e., $D\circ \gamma = \rho(\gamma)\circ D:\tilde M \to \HH^n$. 

A hyperbolic structure on $M$ is {\em complete}\index{complete}\index{hyperbolic structure!complete}
if the induced metric on $M$ is complete.
This condition is equivalent to the condition that
the induced metric on $\tilde M$ is complete,
which in turn is equivalent to the condition that
the developing map $D:\tilde M\to \HH^n$ is an isometry.
Then the holonomy representation $\rho:\pi_1(M)\to \Isom\HH^n$
is faithful and discrete,
nameley $\rho$ gives an isomorphism from $\pi_1(M)$ to 
a discrete torsion-free subgroup, $\Gamma$, of $\Isom \HH^n$.
Thus the complete hyperbolic manifold $M$ is
identified with $\HH^n/\Gamma$.

By a {\it Kleinian group}\index{Kleinian group} we mean a discrete subgroup of $\Isom^+\HH^3$,
and by a {\it Fuchsian group}\index{Fuchsian group} we mean a discrete 
subgroup of $\Isom^+\HH^2$.
By Lemma \ref{lem:centralizer},
any commutative torsion-free Kleinian group
is conjugate to one of the three groups in the following example.

\begin{example}[Commutative torsion-free Kleinian groups]
\label{example:elementary-group}

{\rm
(1) The infinite cyclic group $J_0=J_0(re^{\theta\sqrt{-1}})$ 
generated by the hyperbolic element 
$A(z)=re^{\theta\sqrt{-1}}z$ with $r>1$.
The hyperbolic manifold $\HH^3/J_0$ is homeomorphic to the interior of the solid torus,
and it has the unique closed geodesic with length $\Re\LL_A=\log r$.
For any $r>0$, the closed $r$-neighborhood of $\axis A$ is invariant by $J_0$,
and its
quotient by $J_0$ is called a {\it tube}\index{tube} around the closed geodesic.

(2) The infinite cyclic group $J_1$
generated by the parabolic transformation $A(z)=z+1$.
The hyperbolic manifold $\HH^3/J_1$ is homeomorphic to the product
$\interior D^*\times \RR_+$,
where $D^*=D^2-\{0\}$ is a once-punctured disk.
This hyperbolic manifold does not contain a closed geodesic.
For any $c>0$, the horoball $H_{\infty,c}$ is invariant by $J_1$
and and its quotient by $J_1$ is called
an {\it annulus cusp}\index{annulus cusp}\index{cusp!annulus}.

(3) The rank $2$ free abelian group 
$J_2=J_2(\tau)$ generated by the two parabolic transformations
$A(z)=z+1$ and $B(z)=z+\tau$ with $\tau\in\CC-\RR$.
The hyperbolic manifold $\HH^3/J_2$ is homeomorphic to the product
$T^2\times \RR_+$,
and it does not contain a closed geodesic.
For any $c>0$, the horoball $H_{\infty,c}$ is invariant by $J_2$
and its quotient by $J_2$ is called a 
{\it torus cusp}\index{torus cusp}\index{cusp!torus}.
The boundary torus $\partial H_{\infty,c}/J_2$
admits a Euclidean structure which is conformally equivalent to
the Euclidean torus $\CC/\langle 1, \tau\rangle$. 
Though the cusp neighborhood $H_{\infty,c}/J_2$ is noncompact,
its volume 
$\vol(H_{\infty,c}/J_2)=
\frac{1}{2}\operatorname{area}(\partial H_{\infty,c}/J_2)$
is finite.
The complex number $\tau$ is called the {\it modulus} of the cusp torus
with respect to the basis $\{A,B\}$.
}
\end{example}

For an orientable complete hyperbolic $3$-manifold $M=\HH^3/\Gamma$ 
and a point $x\in M$ 
the {\it injectivity radius}\index{injectivity radius} 
$r(x,M)$ of $M$ at $x$ is defined by
\[
r(x,M)=\sup\{r>0 \ | \ 
\mbox{the $r$-neigborhood of $x$ in $M$ is isometric to an $r$-ball in $\HH^3$}\}.
\]
For a given constant $\epsilon>0$, we can decompose $M$ into
the {\it $\epsilon$-thick part}\index{thick part}
\[
M_{\ge \epsilon}=\{x\in M \ | \ r(x,M)\ge \frac{1}{2}\epsilon\}
\]
and its complement
\[
M_{< \epsilon}=\{x\in M \ | \ r(x,M)< \frac{1}{2}\epsilon\}.
\]
The closure of $M_{< \epsilon}$ is denoted by $M_{\leq \epsilon}$
and is called the {\it $\epsilon$-thin part}\index{thin part} of $M$.
(This complicated definition eliminates the trouble
which occurs when there is a closed geodesic of length $\epsilon$
\cite[p.254]{Thurston_book}.)
The following is a consequence of the Margulis lemma
(see \cite[Theorem 5.10.1 and Corollary 5.10.2]{Thurston0}).

\begin{theorem}
There is a universal constant $\epsilon_0>0$,
such that for any positive constant $\epsilon < \epsilon_0$
and for any orientable complete hyperbolic manifold $M=\HH^3/\Gamma$,
the $\epsilon$-thin part $M_{\leq \epsilon}$
is a disjoint union of 
tubes around (short) simple closed geodesics,
annulus cusps, and torus cusps.
\end{theorem}

The following proposition can be proved 
by using the above theorem
and the concept of convex core introduced in Subsection \ref{subsec:core-conformal-boundary}.
(See \cite[Proposition 5.11.1]{Thurston0}).

\begin{proposition}
\label{prop:finite-volume}
If an orientable complete hyperbolic manifold $M=\HH^3/\Gamma$
has finite volume,
then $M$ is the union of a compact submanifold (bounded by tori)  and 
finitely many torus cusps $C_1,\dots, C_m$ for some $m\ge 0$.
In particular, $M$ is identified with the interior 
of a compact $3$-manifold $\bar M$ with (possibly empty) toral boundary.
\end{proposition}

To end this subsection, 
we recall an important consequence of Thurston's 
hyperbolization theorem for Haken manifolds.

\begin{definition}
\label{def:hyperbolic-link}
{\rm
A knot or link $L$ in $S^3$ is {\it hyperbolic}\index{hyperbolic!link}\index{hyperbolic!knot}\index{link!hyperbolic}\index{knot!hyperbolic}
if its complement $S^3-L\cong \interior E(L)$
admits a complete hyperbolic structure of finite volume.
}
\end{definition}

The following theorem is a special case of Theorem \ref{thm:JSJ-Knot1},
which in turn is a special case of the geometrization theorem.

\begin{theorem}
\label{thm:hyperbolic-link}
A prime knot in $S^3$ is hyperbolic
if and only if it is neither a torus knot nor a satellite knot.
More generally,
an unsplittable prime link $L$ is hyperbolic
if and only if $E(L)$ is atoroidal and is not a Seifert fibered space.
\end{theorem}

\subsection{Rigidity theorem for complete hyperbolic manifolds of finite volume}
\label{subsec:rigidity}

For complete hyperbolic structures of finite volume of dimension $\ge 3$,
the following strong rigidity theorem is established 
by Mostow \cite{Mostow1968} and Prasad \cite{Prasad1973}
(cf. \cite[Theorem 5.7.2]{Thurston0}).

\begin{theorem}[The Mostow-Prasad rigidity theorem]
\label{thm:Mostow-Prasad}
If an orientable $n$-manifold with $n\ge 3$ admits a complete hyperbolic structure of finite volume,
then this structure is unique.
To be precise, the following holds.
Let $\Gamma_i$ ($i=1,2$) be discrete torsion free subgroups
of $\Isom^+ \HH^n$ with $n\ge 3$ of cofinite volume, i.e., 
$\vol(\HH^n/\Gamma_i)<\infty$.
Then any isomorphism $\phi:\Gamma_1\to \Gamma_2$ is realized 
by a unique isometry $f: \HH^n/\Gamma_1 \to \HH^n/\Gamma_2$.
\end{theorem}

This theorem together with Thurston's hyperbolization theorem
had tremendous impact in knot theory.
Because Theorem \ref{thm:hyperbolic-link} says that almost all knots are hyperbolic
(moreover, the Geometrization Theorem \ref{thm:JSJ-Knot1} reduces the study of knots 
to the study of hyperbolic links)
and the above theorem imply that 
geometric invariants, such as volumes, cusp shapes, 
and lengths of shortest closed geodesics,
of the complete hyperbolic structures
on knot/link complements are topological invariants of the knots/links.

\section{Computation of hyperbolic structures and canonical decompositions of cusped hyperbolic manifolds}
\label{sec:computation-canonical}
Epstein and Penner proved that
every cusped hyperbolic manifold of finite volume
admits a natural ideal polyhedral decomposition, called the canonical decomposition.
This fact 
(together with the rigidity theorem
and the Gordon-Luecke knot complement theorem) 
has the following striking consequence in knot theory.
{\it The combinatorial structure of the canonical decomposition
of a hyperbolic knot complement is a complete knot invariant.}
Moreover the marvelous computer program {\it SnapPea}
developed by Jeffrey Weeks
enabled us to compute the canonical decompositions of
knot complements.
For example, SnapPea immediately tells us that
the Kinoshita--Terasaka knot and the Conway's knot are different
and that they admit no symmetry. 

In this section, we recall the Epstein--Penner canonical decomposition 
and its impact on knot theory.
We also recall a method for constructing hyperbolic structures
by using ideal triangulation, which was first explained in 
Thurston's lecture notes \cite[Chapter 4]{Thurston0}, 
and explain a method for finding the canonical decomposition.
In the final subsection,
we give a list of geometric invariants of hyperbolic knots,
which are guaranteed to be knot invariants
by the rigidity theorem,
and introduce their study 
from the viewpoint of {\it effective geometrization}.

\subsection{The canonical decompositions of cusped hyperbolic manifolds}
\label{subsec:canonical}
Let $M=\HH^n/\Gamma$ be an orientable complete hyperbolic $n$-manifold of finite volume
with $m\ge 1$ cusps.
Pick mutually disjoint cusps
$C_1,\dots, C_m$ of $M$, and set $C=\cup_{i=1}^m C_i$.
Then we can canonically construct a spine $\Ford$ 
and a canonical ideal polyhedral decomposition $\DD$ of $M$ as 
follows.

Observe that a generic point in $M-C$
has a unique shortest geodesic path to $C$
but that there are exceptional points which have more than one shortest geodesic paths to $C$.
Let $\Ford$ be the subset of $M-C$  
consisting of these exceptional points.
Namely, $\Ford$ is the {\it cut locus} in $M$ with respect to
the cusps $C=\cup_{i=1}^m C_i$.
Then $\Ford$ is a locally finite totally geodesic cell complex in $M$,
and there is a deformation retraction of $M$ onto $\Ford$.
We call it the {\em Ford complex}\index{Ford complex} 
or {\em Ford spine}\index{Ford spine} of $M$,
with respect to the choice of cusps $C_1, C_2,\dots, C_m$.

By taking the geometric dual to $\Ford$ as follows, 
we obtain an ideal polyhedral decomposition $\DD$ of $M$.
Let $\tilde\Ford$ and $\tilde C$ be the inverse images of $\Ford$ and $C$
in the universal covering $\HH^n$ of $M$.
Pick a vertex $x$ of $\tilde\Ford$.
Then there are finitely many shortest geodesic paths from $x$ to $\tilde C$.
Let $\{v_i\}$ be the ideal points in $\partial \HH^n$ 
forming the centers of the horoball components of $\tilde C$ 
which are joined to $x$ by a shortest geodesic path.
The convex hull of the ideal points $\{v_i\}$ forms an 
$n$-dimensional {\it ideal polyhedron}\index{ideal polyhedron} of $\HH^n$,
and the collection of all such ideal polyhedra, 
where $x$ runs over the vertices of $\tilde\Ford$,
determines a $\Gamma$-invariant tessellation of $\HH^n$.
The tessellation descends to an ideal polyhedral decomposition $\DD$ of
$M=\HH^n/\Gamma$.

Epstein and Penner \cite{Epstein-Penner} gave a description of 
the decomposition $\DD$
by using a convex hull construction in Minkowski space.
Their description shows that each cell of $\DD$
admits a natural (incomplete) Euclidean structure:
so, these decompositions are called {\em Euclidean decompositions}\index{Euclidean decomposition}.
In \cite{Akiyoshi-Sakuma}, 
a generalization of the Epstein--Penner construction to 
cusped hyperbolic manifolds of infinite volume is given,
and their relationship to the convex cores are discussed.

The Ford complex $\Ford$ and its geometric dual $\DD$ depend only on the ratio of the volumes
$\vol (C_1):\vol(C_2):\cdots:\vol(C_m)$.
Moreover, it is proved by Akiyoshi \cite{Akiyoshi}
that the combinatorial structures of $\Ford$, 
when the ratio varies, are finite.
The ideal polyhedral decomposition $\DD$,
for the case when $C_1, C_2,\dots, C_m$ have the same volume,
is uniquely determined by the hyperbolic manifold $M$,
and is called the {\it canonical decomposition}\index{canonical decomposition} of the cusped hyperbolic manifold $M$.

\begin{example}
{\rm
(1) Let $M$ be the hyperbolic thrice-punctured sphere,
obtained by gluing two ideal hyperbolic triangles
through identification of their boundaries via the identity map.
Then this decomposition of $M$ into the two copies of ideal triangles
is the canonical decomposition of $M$.
The corresponding Ford complex of $M$ is a $\theta$-shaped geodesic spine of $M$
consisting of two vertices and three edges.

(2) As shown by Thurston \cite[Chapter 4]{Thurston1},
the complete hyperbolic structure of the
figure-eight knot complement $M$ is obtained by glueing
two copies of the regular ideal tetrahedron.
The decomposition of $M$ into the 
two copies of the regular ideal tetrahedron
is the canonical decomposition of $M$.
}
\end{example}

Since the complete hyperbolic structure of a given knot complement is unique
by the Mostow-Prasad rigidity theorem
(Theorem \ref{thm:Mostow-Prasad}),
and since 
by the knot complement theorem
(Theorem \ref{thm:knot-complement-thm})
a knot is determined by its complement,
it follows that the combinatorial structure of the canonical decomposition
of a hyperbolic knot complement is a complete topological invariant of the knot.

\begin{theorem}
\label{thm:canonical}
(1) Two hyperbolic knots are equivalent,
if and only if the canonical decompositions of their complements are
combinatorially equivalent.

(2) Let $K$ be a hyperbolic knot and $\DD$ the canonical decomposition of 
$S^3-K$.
Then 
\[
\Sym(S^3,K)\cong \Isom (S^3-K) \cong \Aut (\DD).
\] 
\end{theorem}

In the above theorem, $\Sym(S^3,K):=\pi_0 \Diff(S^3,K)$ denotes 
the {\em symmetry group}\index{symmetry group}\index{knot!symmetry group} of the knot $K$,
and $\Aut(\DD)$ denotes the combinatorial automorphism group of $\DD$.

\begin{example}
\label{example:KT-Conway_canonical}
{\rm
(1) It is a simple exercise to see that the automorphism group
of the canonical decomposition $\DD$ of the complement figure-eight knot $K$
is isomorphic to the order $8$ dihedral group $D_8$.
Thus we have $\Sym(S^3,K)\cong \Isom (S^3-K) \cong D_8$.

(2) The canonical decompositions of the complements of
the Kinoshita--Terasaka knot and the Conway knot
consist of 
$12$ and $14$ ideal tetrahedra, respectively.
Hence they are inequivalent,
even though they are mutants of each other and so
they share the same Alexander polynomial, the Jones polynomial, the hyperbolic volume
and the same double branched coverings.
Moreover, the automorphism groups of both canonical decompositions are trivial.
Thus the symmetry groups of these two knots are trivial.
In particular, both of them are neither amphicheiral nor invertible.
The noninvertibility of $8_{17}$ can be also proved by
using the canonical decomposition of the knot complement.
}
\end{example}

As explained in the next subsection,
the canonical decompositions are amenable to computer calculation,
and wonderful computer programs were developed:
{\it SnapPea}\index{SnapPea} by Weeks \cite{Weeks1}, 
{\it Snap}\index{Snap} by Coulson-Goodman-Hodgson-Neumann \cite{CGHW},
{\it SnapPy}\index{SnapPy} by Culler-Dunfield-Goerner\cite{CDG}, and 
a computer verified program 
{\it HIKMOT}\index{HIKMOT} 
by Hoffman-Ichihara-Kashiwagi-Masai-Oishi-Takayasu \cite{HIKMOT}.
The results in Example \ref{example:KT-Conway_canonical}(2) are, 
of course, obtained
by any of these programs.

This enabled
Hoste, Thistlethwaite and Weeks \cite{HTM}
to extend (and correct) Conway's enumeration of all $11$ crossing knots
to include all prime knots up $16$ crossings.
There are $1,701,936$ such knots, and all except for $32$ knots are hyperbolic!
To be precise, Hoste and Weeks used the canonical decomposition,
and Thistelethwaite used the \lq\lq universal method'' described at the end of 
Subsection \ref{subsec:universal-method}.
Thus their table is double checked,
and this fact shows the strength of both methods.

This is something like a magic wand for knot theorists
as long as finitely many knots of reasonable crossing numbers are concerned.
However, to understand the canonical decompositions of infinite families of knots 
or cusped hyperbolic manifolds is not easy.
For the Farey manifolds, namely punctured torus bundles and 2-bridge knot complements,
the combinatorial structures of the canonical decompositions are determined
by Jorgensen \cite{Jorgensen} and Gu\'eritaud \cite{Gueritaud2006b}
(cf. \cite{ASWY, Gueritaud2006a, Sakuma-Weeks}).

To end this subsection,
we remark that it is still an open problem
to see whether every orientable cusped hyperbolic $3$-manifold of finite volume
admits a {\it ideal triangulation}\index{ideal triangulation},
namely an ideal polyhedral decomposition
consisting of ideal tetrahedra.
Here an {\it ideal tetrahedron}\index{ideal tetrahedron} is a closed convex hull in $\HH^3$
of $4$ ideal points in $\partial\HH^3$,
called the {\it ideal vertices}\index{ideal vertex}.
(Since any such manifold $M$ admits an ideal polyhedral decomposition
by \cite{Epstein-Penner}  
and since every ideal polyhedron is decomposed into ideal tetrahedra,
$M$ admits a {\it partially flat} ideal triangulation, 
namely one in which some of the tetrahedra degenerate into 
flat quadrilaterals with distinct vertices
(see \cite{Petronio-Porti2000}).
But this does not necessarily lead to a genuine ideal triangulation of $M$.)
Wada, Yamashita and Yoshida \cite{Wada-Yamashita-Yoshida1966} and Yoshida \cite{Yoshida1966}  
proved the existence of such triangulations 
under certain combinatorial conditions on the polyhedral decomposition, and 
Luo, Schleimer and Tillman
\cite{Luo-Schleimer-Tillman} proved that every 
such manifold virtually admits an ideal triangulation,
namely some finite cover has an ideal triangulation.
Hodgson, Rubinstein and Segerman \cite{Hodgson-Rubinstein-Segerman2012}
considered a relaxed version of the problem,
and proved, in particular, that every hyperbolic link complement in $S^3$
admits a topological ideal triangulation with a \lq\lq strict angled structure''.

\subsection{Ideal triangulations and computations of the hyperbolic structures}
\label{sebsec:ideal-triangulation}

Let $M=\HH^3/\Gamma$ be an orientable complete hyperbolic $3$-manifold of finite volume
with $m\ge 1$ cusps,
and let $\rho:\pi_1(M) \to \Gamma < \PSL(2,\CC)$ be the holonomy representation.
Then, as we have observed in the previous section,
$M$ admits an ideal polyhedral decomposition $\DD$.
We now assume that $\DD$ is an ideal triangulation,
namely $\DD$ consists of ideal tetrahedra.
Any ideal tetrahedron $\Delta$ (up to isometry) is represented by a complex number $z$
with positive imaginary part,
such that the Euclidean triangle cut out of any vertex of $\Delta$
by a horosphere is similar to the triangle in $\CC$
with vertices $0$, $1$, and $z$.
In fact, $\Delta$ is isometric to the ideal tetrahedron $\Delta(z)$
spanned by $0$, $1$, $\infty$ and $z$ in the upper half-space model $\CC\times \RR_+$
of $\HH^3$.
We call $z$ the {\it shape parameter} of the ideal tetrahedron $\Delta(z)$.
(If $z$ has negative imaginary part,
then $\Delta(z)$ is regarded as negatively oriented.
If $z$ is a real number different from $0$ and $1$,
then $\Delta(z)$ is regarded as a degenerate ideal tetrahedron.)
The complex numbers $z$, $(z-1)/z$, and $1/(1-z)$ give isometric ideal tetrahedra,
and we give each edge $e$ of $\Delta=\Delta(z)$ 
one of the three complex numbers by the following rule,
and call it the {\it edge parameter}\index{edge parameter}\index{ideal tetrahedron!edge parameter} of $\Delta$ associated with $e$.
\begin{itemize}
\item[$\circ$]
Edges $[0,\infty]$ and $[1,z]$ have edge parameter $z$.
\item[$\circ$]
Edges $[1,\infty]$ and $[0,z]$ have edge parameter $1/(1-z)$.
\item[$\circ$]
Edges $[z,\infty]$ and $[0,1]$ have edge parameter $(z-1)/z$.
\end{itemize}

Let $e$ be an edge of the ideal triangulation $\DD$ of $M$,
and let $z_1,\dots, z_k$ be the edge parameter of the edges of
ideal tetrahedra glued to $e$.
Since these ideal tetrahedra close up as one goes around $e$,
the parameters satisfies the following equation.
\[
\prod_{i=1}^k z_j=1 \quad{\mbox{\rm and}} \quad \sum_{j=1}^k \arg(z_j)=2\pi
\]
This condition is identical to the following equation,
which is called the {\it gluing equation}\index{gluing equation} around $e$.
\[
\sum_{j=1}^k \log(z_j)=2\pi \sqrt{-1},
\]
where $\log:\CC-\RR_{\leq 0}\to \CC$ is the branch of the logarithm function
whose imaginary part lies in $(-\pi,\pi)$.

Let $T$ be a torus boundary component of the compact $3$-manifold $\bar M$
whose interior is homeomorphic to the hyperbolic manifold $M=\HH^3/\Gamma$,
and let $\mu$ be an oriented essential simple loop on $T$.
(A simple loop on $T$ is {\it essential}\index{essential} if it does not bound a disk in $T$.)
By identifying $T$ with a cusp torus,
and considering the intersection with the cusp torus with 
the ideal triangulation $\DD$,
we obtain a triangulation of $T$,
whose vertices correspond to the edges of $\DD$
and whose triangles correspond to truncations of ideal tetrahedra
around ideal vertices.
We may assume that $\mu$ intersects the edges of the triangulation transversely 
and does not intersect the vertices of the triangulation.
Each segment of $\mu$ in a triangle cuts off a single vertex of the triangle,
and so has an associated edge parameter $z_j$.
Define $\epsilon_j=+1$ or $-1$
according to whether the vertex lies to the {\it left} of $\mu$ or not.
(Here we assume that 
$\infty$ is a parabolic fixed point of $\Gamma$,
$\pi_1(T)$ is identified with the stabilizer $\Gamma_{\infty}$ of $\infty$,
and $T$ is identified with the Euclidean torus $\CC/\Gamma_{\infty}$
via the projection from a horosphere centered at $\infty$ to $\CC$.
The left/right convention is determined by the standard orientation of $\CC$.)
Then we can see that the complex translation length of 
the image $\rho(\mu)$ of $\mu$ 
by the holonomy representation $\rho$
of the complete hyperbolic manifold $M$
 is represented by the complex number
\[
\LL_{\mu}:=\sum_{j} \epsilon_j\log(z_j).
\]
Since $\rho(\mu)$ is parabolic, we have $\LL_{\mu}=0$.
Thus we have
the following {\it completeness equation}
\[
\sum_{i} \epsilon_i\log(z_i)=0.
\]

Conversely, let $\bar M$ be an orientable compact manifold whose boundary is non-empty
and consists of tori,
and let $\DD$ be a {\it topological ideal triangulation}\index{topological ideal triangulation}\index{ideal triangulation!topological} of $M=\interior \bar M$.
Namely $\DD$ is a topological trianglulation 
(a cell decomposition whose cells are identified with simplices)
of the space $\hat M=\bar M/\sim$,
where $\sim$ is the equivalence relation which identifies 
all points of each boundary component of $\bar M$,
such that the vertex set of $\DD$ is equal to the finite set
consisting of the image of $\partial \bar M$. 
By a simple argument using the Euler characteristic,
we see that the number of edges in $\DD$
is equal to the number, $t$, of tetrahedra in $\DD$.
Now let $\HH_+=\{z\in\CC \ | \ \Im z>0\}$ be the upper-half space of the complex plane.
Pick a $t$-tuple of complex numbers 
$\zb=(z_1,\dots,z_t)\in (\HH_+)^t\subset\CC^t$ 
with positive imaginary parts, and 
identify the topological ideal tetrahedra $\{\Delta_1,\dots,\Delta_t\}$ 
with hyperbolic ideal tetrahedra $\{\Delta(z_1),\dots,\Delta(z_t)\}$.
Since all hyperbolic ideal triangles are isometric,
we can realize the topological gluing maps among the
faces of the topological ideal tetrahedra by hyperbolic isometries.
Thus we obtain a hyperbolic structure on the complement of the $1$-skeleton
of $\DD$.
We have the following theorem
(see \cite{Thurston0, Neumann-Zagier}).

\begin{theorem}
\label{thm:gluing}
Under the above setting, the following hold for each $\zb=(z_1,\dots,z_t)\in (\HH_+)^t\subset\CC^t$.
\begin{enumerate}
\item
The hyperbolic structure on the complement of the $1$-skeleton of $\DD$
extends a hyperbolic structure on the whole $M$
if and only if $\zb$ satisfies the gluing equation 
at every edge of $\DD$.
\item
When condition (1) is satisfied,
the resulting hyperbolic structure on $M$ is complete
if and only if $\zb$ also satisfies the completeness equation at 
every boundary component of $\bar M$
(for a single choice of an oriented essential simple loop $\mu$ 
for each boundary component).
\end{enumerate}
\end{theorem}

\begin{remark}
\label{rem:dim-deformation}
{\rm
Let $\XX$ be the variety of $\zb=(z_1,\dots,z_t)\in\CC^t$
consisting of the solutions of the gluing equations.
Then, by a combinatorial argument,
we can see that $\XX$ has dimension $m$ over $\CC$, 
where $m$ is the number of boundary components of $\bar M$
(\cite[Theorem 5.6]{Thurston0}, \cite[Proposition 2.3]{Neumann-Zagier}).
By the rigidity theorem,
there is a unique point $\zb^0\in \XX\cap (\HH_+)^t$
which satisfies the completeness equation.
It is proved by \cite[Section 4]{Neumann-Zagier}
that $\zb^0$ is a smooth point of $\XX\cap (\HH_+)^t$,
namely  
there is a neighborhood of $\zb^0$ in $\XX\cap (\HH_+)^t$ 
which is biholomorphically equivalent to an open set in $\CC^m$.
(Moreover, it was proved by Choi \cite{Choi2004} that
$\XX\cap (\HH_+)^t$ is a smooth complex manifold.)
This fact plays a crucial role in a proof of 
the hyperbolic Dehn filling theorem (see Subsection \ref{subsec:proof-HDFT}).
}
\end{remark}

On the other hand, 
there is a convenient method for obtaining topological ideal triangulations 
of knot/link complements from diagrams
(see \cite{Menasco, Takahashi, Weeks3}).
Thus we have a good chance to construct a complete hyperbolic structure 
on a given knot/link complement by applying Theorem \ref{thm:gluing}.
In fact, this works extremely well,
though the proof of Thurston's uniformization theorem is very difficult.

Moreover, if a given ideal triangulation $\DD$ of $M$
satisfies a certain inequality at each codimension $1$ face of $\DD$,
then $\DD$ is the canonical decomposition (see \cite{Weeks2}).
If the inequality was not satisfied at some face of $\DD$,
then apply the Pachner $3-2$ move to $\DD$ at the face, if it is geometrically realizable,
and check if the conditions for the faces hold.
If this does not lead to the canonical decomposition,
then retriangulate $\DD$ randomly, and repeat the above procedure.
This is the way SnapPea finds the canonical decompositions.
Though there is no theoretical guarantee,
this method is extremely efficient (see \cite{Weeks2, Weeks3}).
For the treatment of the case 
when the canonical decomposition is not an ideal triangulation,
see the work of Hodgson and Weeks \cite{Hodgson-Weeks}.

\subsection{Other geometric invariants for hyperbolic knots and effective geometrization}
\label{subsec:other-invariants}

In addition to the canonical decomposition,
there are various important geometric invariants of hyperbolic knots and links.

\begin{itemize}
\item[$\circ$]
The volumes and the Chern-Simons invariants 
of the hyperbolic link complements.
\item[$\circ$]
The volumes of the maximal cusps.
\item[$\circ$]
The moduli of the Euclidean cusp tori.
\item[$\circ$]
Length spectrum, i.e., the multi-set of lengths
of closed geodesics,
in particular the length of the shortest closed geodesic.
\item[$\circ$]
Lengths of the vertical geodesic paths,
joining maximal cusps to themselves.
\item[$\circ$]
Euclidean length spectrum of the maximal cusp torus.
\end{itemize}

Volumes of hyperbolic manifolds are
treated in Section \ref{sec:volume}.

In the recent beautiful survey \cite{FKP_survey},
Futer, Kalfagianni, and Purcell
discuss these invariants from the viewpoint
of {\it effective geometrization}
or {\it WYSIWYG topology},
where WYSIWYG stands for \lq\lq what you see is what you get'',
which aims to determine the geometry of link complements directly
from the link diagrams.
A typical example in this direction is the following estimate by Lackenby \cite{Lackenby2}
of the volume of alternating link complements
in terms of the twist number.

\begin{theorem}
Let $D$ be a reduced alternating diagram of a hyperbolic link $L$ in $S^3$,
and let $t(D)$ be the twist number of the diagram $D$. Then
\[
\frac{1}{2}V_{\mathrm{tet}}(t(D)-2) \le \vol(S^3-L) \le  V_{\mathrm{tet}}(16 t(D) - 16)
\]
where $V_{\mathrm{tet}}=1.0149416...$ is the volume of 
the regular ideal tetrahedron.
\end{theorem}

\noindent
Here the {\it twist number}\index{twist number}\index{diagram!twist number} $t(D)$ of a link digram $D$ is the number of twists of $D$,
where a {\it twist}\index{twist}\index{diagram!twist} of $D$ is 
either a connected collection of bigon regions in $D$ arranged in a row 
which is maximal in the sense that it is not part of a longer row of bigons,
or a single crossing adjacent to no bigon regions.  

The article \cite{FKP_survey} presents 
a nice survey on the recent great progress towards
effective geometrization,
including a refinement of the above result.

\section{Flexibility of incomplete hyperbolic structures and the hyperbolic Dehn filling theorem}
\label{sec:Hyperbolic-Dehn-filling}

By the Mostow-Prasad rigidity theorem,
the complete hyperbolic structure on a $3$-manifold $M$ of finite volume is rigid.
However, when $M$ has a cusp, the complete hyperbolic structure
admits nontrivial continuous
deformations into incomplete hyperbolic structures
(see Remark \ref{rem:dim-deformation}).
In the generic case,
the metric completion yields
a pathological topological space which is not even Hausdorff.
However, in certain special isolated cases, 
the metric completion produces a complete hyperbolic manifold.
This is a rough idea of 
Thurston's hyperbolic Dehn filling Theorem.
This theorem has stimulated keen attention of many mathematicians 
and enormous amount of research grew out of this result.
In this section, we give an outline of a proof of this theorem
and a brief survey of its influence on knot theory.

\subsection{Hyperbolic Dehn filling theorem}
\label{subsec:HDF-theorem}
We begin by recalling the topological operation, Dehn filling.
By an  {\it oriented slope}\index{oriented slope}\index{slope!oriented} on a torus $T$, we mean the isotopy class of an oriented
essential simple loop on $T$.
Each oriented slope  
represents a primitive element of $H_1(T;\ZZ)$, and conversely
any primitive element of $H_1(T;\ZZ)$ is represented by a unique oriented slope.
If we fix a basis $\{\mu,\lambda\}$ of $H_1(T;\ZZ)$,
then a primitive element of $H_1(T;\ZZ)$
is expressed as $p\mu + q\lambda$
where $(p,q)$ is a pair of relatively prime integers.
Thus we can identify the set of
oriented slopes on $T$ 
with the set of pairs of relatively prime integers $(p,q)\in \ZZ^2\subset \RR^2\cup \{\infty\}\cong S^2$.

Let $M$ be a connected compact orientable $3$-manifold whose boundary consists of 
$m$ tori $T_1,\dots,T_m$.
Pick an oriented slope $\nu_j$ on $T_j$ for each $j$,
and attach a solid torus $V_j=D_j^2\times S^1$ to $M$ along $T_j$,
so that the meridian $\partial D_j^2 \times \{*\}$ is identified with the slope $\nu_j$.
The resulting manifold is denoted by $M(\nub)=M(\nu_1,\dots,\nu_{m})$
and called the result of {\it Dehn filling}\index{Dehn filling} of $M$ along the 
tuple $\nub=(\nu_1,\dots,\nu_{m})$ of oriented slopes.
We extend this operation to the case where some $\nu_j$ is the symbol $\infty$,
by the rule that
if $\nu_j=\infty$ then we leave the boundary $T_j$ as it is.
In particular, $M(\infty,\dots,\infty)=M$. 

The following theorem is proved by Thurston \cite[Chapters 4 and Section 5.8]{Thurston0}.

\begin{theorem}
[Hyperbolic Dehn filling Theorem]\index{hyperbolic Dehn filling theorem}
\label{thm:hyp-Dehn}
Let $M$ be a connected compact orientable $3$-manifold whose boundary consists of 
$m$ tori, and
suppose that $\interior M$ admits a complete hyperbolic structure of finite volume.
Then, except for finitely many choices of the slopes of $\nu_j$ for each $1\le j\le m$,
the manifold $M(\nu_1,\dots,\nu_{m})$ admits a complete hyperbolic structure.
To be more precise, there exists a neighborhood $V$
of $(\infty,\dots,\infty)$ in $(\RR^2\cup \{\infty\})^m$
such that $M(\nu_1,\dots,\nu_{m})$ admits a complete hyperbolic structure
for every slope $(\nu_1,\dots,\nu_{m})$ contained in $V$.
\end{theorem}

\begin{remark}
{\rm
(1) The operation at $T_j$ is actually determined by the {\it slope}\index{slope}
(the isotopy class of an unoriented essential simple loop on a torus) 
obtained from
$\nu_j$ by forgetting the orientation.

(2)
When $M$ is the exterior of an $m$-component link $L=\cup_{j=1}^{m} K_j$ in $S^3$,
we fix an orientation of each component $K_j$ of $L$, and
choose the meridian-longitude systems $\{\mu_j,\lambda_j\}$ as  
a preferred basis for $H_1(T_j;\ZZ)$,
and represent an oriented slope, $\nu_j$, on $T_j$ 
by a pair of relatively prime integers $(p_j,q_j)$
with $\nu_j=p_j\mu_j+q_j\lambda_j$.
The slope obtained from $\nu_j$ by forgetting the orientation
is uniquely determined by  the rational number $p_j/q_j\in \QQ\cup \{1/0\}$.
(It should be noted that slope $1/0$ and the symbol $\infty$ have different meanings.)
Moreover, this does not depend on the choice of the orientation of $K_j$.
We denote the manifold $M(\nu_1,\dots,\nu_{m})$ by $M(p_1/q_1,\dots, p_m/q_m)$,
and call it the result of {\it Dehn surgery}\index{Dehn surgery} on $L$
with slope $(p_1/q_1,\dots, p_m/q_m)$.
}
\end{remark}

In Theorem \ref{subsec:HDF-theorem}, 
a slope (or a tuple of slopes) which does not produce a hyperbolic manifold
is called an {\it exceptional slope}\index{exceptional slope}\index{slope!exceptional}.

\begin{example}
\label{ex:exceptional-slope}
{\rm
(1) The exceptional slopes of the figure-eight knot $K$
are the slopes $p/q$ with $-4\le p\le 4$ and $-1\le q\le 1$.
Thus the set of exceptional slopes is
$\{1/0, 0, \pm1,\pm2,\pm 3, \pm 4\}$
(see \cite[Section 4.6]{Thurston0}).

(2) Let $M$ be the exterior of the Whitehead link $L=K_1\cup K_2$ in $S^3$.
Consider the Dehn filling only along $T_1=\partial N(K_1)$.
Then the exceptional slopes for this Dehn filling
are those slopes contained in the parallelogram
with vertices $\pm(4,-1)$ and $\pm(0,1)$
(see \cite[Section 6]{Neumann-Reid1992}).
}
\end{example}

\subsection{Outline of a proof 
and generalized Dehn filling coefficients}
\label{subsec:proof-HDFT}
We give an outline of the proof of Theorem \ref{thm:hyp-Dehn}
by Neumann-Zagier \cite{Neumann-Zagier}
(cf. \cite[Section E.6]{Benedetti-Petronio}),
when the hyperbolic manifold $\interior M$ admits an ideal triangulation $\DD$. 
(See Petronio-Porti \cite{Petronio-Porti2000} for a proof without assuming the existence
of an ideal triangulation,
and using a partially flat ideal triangulation of $M$.)
Let $\Delta_1,\dots,\Delta_t$ be the ideal tetrahedra in $\DD$,
and let $\zb^0=(z_1^0,\dots,z_t^0)$ be their shape parameters.
By the rigidity theorem and Theorem \ref{thm:gluing},
$\zb^0$ is the unique solution of the gluing and the completeness equations.  
Let $\XX$ be the variety of $\zb=(z_1,\dots,z_t)\in\CC^t$
consisting of the solutions of the gluing equations.
For $\zb\in\XX\cap (\HH_+)^t$,
let $M_{\zb}$ be the (almost certainly incomplete) hyperbolic manifold
determined by the parameter $\zb$,
and let $\rho_{\zb}:\pi_1(M) \to \PSL(2,\CC)$ be the holonomy representation
of $M_{\zb}$. 
For each boundary component $T_j$ of $M$ ($1\le j\le m$),
fix an oriented slope $\mu_j$.
For $\zb\in\XX\cap (\HH_+)^t$,  
let $u_j(\zb)$ be the complex number $\LL_{\mu_j}$,
defined as in Subsection \ref{sebsec:ideal-triangulation},
which represents the complex translation length of $\rho_{\zb}(\mu_j)$.
(Though the complex translation length is defined only modulo $2\pi \sqrt{-1}\ZZ$
and up to multiplication by $\pm 1$, 
the construction in Subsection \ref{sebsec:ideal-triangulation}
gives a well-defined continuous lift to $\CC$.)

For each boundary component $T_j$, 
pick an oriented slope $\lambda_j$
which intersects $\mu_j$ transversely in a single point
and so $\{\mu_j,\lambda_j\}$ forms a generator system of $H_1(T_j;\ZZ)$.
Let $\vb:=(v_1,\dots, v_m)$ be the map from $\XX$ to $\CC^m$, where
$v_j(\zb)$ is the complex number $\LL_{\lambda_j}$,
defined as in Subsection \ref{sebsec:ideal-triangulation},
which represents the complex translation length of $\rho_{\zb}(\lambda_j)$.

Recall the key Remark \ref{rem:dim-deformation} 
that there is a neighborhood of $\zb^0$ in $\XX\cap (\HH_+)^t$ 
which is biholomorphically equivalent to an open set in $\CC^m$.
By using this fact, we can see that 
$\ub:=(u_1,\dots,u_m)$ maps a neighborhood of $\zb^0\in\XX\cap (\HH_+)^t$
biholomorphically onto a neighborhood, $\XX_0$, of $0\in\CC^m$
(cf. \cite[Section 4]{Neumann-Zagier}).
 
We now change notation as follows.
For $\ub\in \XX_0$, we denote the corresponding hyperbolic manifold and
the holonomy representation by $M_{\ub}$ and $\rho_{\ub}$, respectively,
and we regard $\vb$ as a map from $\XX_0$ to $\CC^m$.

By replacing $\XX_0$ with a smaller neighboorhood of $0$, if necessary,
we can assume that
$\ub$ and $\vb$ are independent over $\RR$, for all $\ub\in\XX_0-\{0\}$.
In fact, there is an analytic function 
$\taub=(\tau_1,\dots,\tau_m):\XX_0\to \CC^m$,
satisfying the following conditions \cite[Lemma 4.1]{Neumann-Zagier}:
\begin{enumerate}
\item
$v_j(\ub)=\tau_j(\ub) u_j$ for every $\ub=(u_1,\dots, u_j,\dots,u_m)\in \XX_0$
and  $j=1,\dots,m$. 
\item
$\tau_j(0,\dots,0)$ is equal to the modulus
of the cusp torus $T_j$ of the complete hyperbolic manifold $M$
with respect to $\{\mu_j,\lambda_j\}$.
\end{enumerate}
In particular, we may assume that $\tau_j(\ub)$ is non-real for 
every $\ub\in\XX_0$, and so $\ub$ and $\vb$ are independent over $\RR$
for every $\ub\in \XX_0-\{0\}$.

Now we define the 
{\em generalized Dehn filling coefficients}\index{generalized Dehn filling coefficient} of the  
$j$-th boundary torus component 
$\nu_j\in \RR^2\cup \{\infty\}\cong S^2$ by the formula:
\[
\begin{cases}
\nu_j=\infty \quad \mbox{if $u_j=0$}\\
\nu_j=(p_j, q_j) \quad \mbox{where $p_ju_j + q_jv_j= 2\pi \sqrt{-1}$} \quad \mbox{if $u_j\ne 0$}
\end{cases}
\]

The hyperbolic Dehn filling Theorem \ref{thm:hyp-Dehn} is a consequence
of the following theorem.

\begin{theorem}
\label{thm:hyp-Dehn2}
Under the above setting,
the \lq\lq generalized Dehn filling coefficients map'' $\ub\mapsto \nub=(\nu_1,\dots, \nu_m)$
gives a homeomorphism 
from a neighborhood $U\subset \XX_0$ of $0$ in $\CC^m$ 
onto a neigborhood $V$ of $(\infty,\dots, \infty)$ in
$(\RR^2\cup \{\infty\})^m$.
Moreover, the following hold.

\begin{itemize}
\item[$\circ$]
If $\nu_j=\infty$,
the hyperbolic structure at the $j$-th end is complete.
\item[$\circ$]
If $\nu_j=(p_j,q_j)$ where $p_j, q_j\in\ZZ$ are coprime,
then the completion of the $j$-th end is a hyperbolic $3$-manifold,
which is topologically the Dehn filling such that
the simple loop $p_j\mu_j+q_j\lambda_j$ on $T_j$ bounds a disk.
\item[$\circ$]
When $p_j/q_j\in\QQ\cup \{\infty\}$,
let $m_j,n_j\in\ZZ$ be coprime integers such that
$(p_j,q_j)=d(m_j,n_j)$ for some $d>0$.
The completion is a hyperbolic cone $3$-manifold 
obtained by gluing a solid torus with singular core,
such that 
the simple loop $m_j\mu_j+n_j\lambda_j$ on $T_j$ bounds a disk
which has a singularity at the center,
and that the cone angle of the singular locus is $2\pi/d$.
\item[$\circ$]
If $p_j/q_j\in\RR-\QQ$,
then the metric completion of the $j$-th end is not even topologically a manifold.
\end{itemize}
\end{theorem}

In the above, a {\em hyperbolic cone $3$-manifold}\index{hyperbolic cone manifold} is 
a smooth $3$-manifold $C$ equipped with a complete metric (distance function)
which is locally isometric to $\HH^3$
or to the space $\HH^3(\alpha)$ obtained from a geodesic cheese-cake-shaped polyhedron
of angle $\alpha>0$
by identifying two sides.
The singular locus $\Sigma\subset C$ is the set of points 
modeled on the singular line of some $\HH^3(\alpha)$,
and $\alpha$ is called the {\em cone angle}\index{cone angle} at a singular point
modeled on this singular line
(for precise definition, see 
\cite[Section 1]{Hodgson-Kerckhoff1998},
\cite[Chapter 3]{CHK2000},
\cite[Chapter 1]{Boileau-Porti}, \cite[Section 3]{BLP}).
Hyperbolic $3$-cone manifolds play a key role in the proof of the orbifold theorem
(Theorem \ref{orbifold-theorem}).

\begin{remark}
\label{rem.parameter-space}
{\rm
(1) 
Assume that a tuple of oriented slopes
$\nub=(\nu_1,\dots,\nu_{m})$ is the image of
a paremeter $\ub=(u_1,\cdots,u_m)\in U$ in Theorem \ref{thm:hyp-Dehn2},
namely the metric completion of the hyperbolic manifold $M_{\ub}$
is homeomorphic to the manifold $M(\nub)$
obtained from $M$ by Dehn filling along $\nub$.
Let $\nub'$ be the tuple of oriented slopes
obtained from $\nub$ by replacing some component 
$\nu_j=(p_j,q_j)$ with $-\nu_j=(-p_j,-q_j)$.
Then $\nub'$ is the image of the parameter $\ub'$
obtained from $\ub$ by replacing the component $u_j$ with $-u_j$. 
Since $M(\nub')$ is homeomorphic to $M(\nub)$
by a homeomorphism preserving the subspace $M$,
the rigidity theorem implies that
$M_{\ub'}$ is isometric to $M_{\ub}$.
In fact, such an isometry exists whenever two parameters $\ub$ and $\ub'$
are related by the involution
$(u_1,\dots, u_j,\dots,u_m)\mapsto (u_1,\dots, -u_j,\dots,u_m)$.
Thus deformations of the complete hyperbolic manifold $M$
are parametrized by the quotient of $U$ by the $(\ZZ/2\ZZ)^m$-action,
generated by the above involutions with $j=1,\dots, m$.
In other words, the space $U$ is identified with a $(\ZZ/2\ZZ)^m$-branched covering
of a deformation space of $M$.
The space $U$ actually parametrizes the incomplete hyperbolic manifolds $M_{\ub}$
endowed with an ideal triangulation
(see \cite[p.323]{Neumann-Zagier}).   

(2) 
In Theorem \ref{thm:hyp-Dehn2},
the complete hyperbolic manifolds $\{M(\nub)\}$
are regarded as discrete points in the deformation space 
$U/(\ZZ/2\ZZ)^m \cong V/(\ZZ/2\ZZ)^m$.
Thus the discrete set of complete hyperbolic manifolds 
$\{M(\nub)\}$ are linked together 
in the connected space $V/(\ZZ/2\ZZ)^m$.
}
\end{remark}

\subsection{Geometry of the hyperbolic manifolds obtained by Dehn fillings}
\label{subsec:Geometry-Dehn-filling}

In the hyperbolic Dehn filling Theorem \ref{thm:hyp-Dehn},
the complete hyperbolic manifolds 
$M(\nub)=M(\nu_1,\dots,\nu_{m})$ 
{\it geometrically converge}\index{geometrically converge} to
the original complete hyperbolic manifold
$\interior M$
as $\nub=(\nu_1,\dots,\nu_{m})\to \inftyb=(\infty,\dots,\infty)$
\cite[Section 5.11]{Thurston0}.
Namely, there are positive numbers $\epsilon(\nub)$ 
converging to $0$ as $\nub\to\inftyb$,
and numbers $k(\nub)>1$
converging to $1$ as $\nub\to\inftyb$,
such that there is a
$k(\nub)$-bi-Lipschitz diffeomorphism
\[
\phi_{\nub}:M(\nub)_{\ge \epsilon(\nub)}\to (\interior M)_{\ge \epsilon(\nub)}
\]
between that $\epsilon(\nub)$-thick parts.
This in particular implies that
the lengths of core loops of the attached solid tori in $M(\nub)$ 
converge to $0$ as $\nub\to\inftyb$.
This fact plays an essential role in various researches,
including \cite{Kojima1986, Hodgson-Weeks, Bleiler-Hodgson-Weeks1999, 
Rieck-Yamashita2016, Lackenby*}.

This also implies that the volumes $\vol(M(\nub))$ of the Dehn filled manifolds
converge to the volume $\vol(\interior M)$ of the original hyperbolic manifold
as $\nub\to\inftyb$.
Moreover,
Thurston \cite{Thurston0} proved, by using the Gromov norm
(cf. Subsection \ref{subsec:Gromov-norm}), that
$\vol(M(\nub))$ is strictly smaller than  $\vol(M)$ if $\nub\ne\inftyb$.
This is refined to quantitative estimates of $\vol(M(\nub))$ 
by Neumann-Zagier \cite{Neumann-Zagier}, Hodgson-Kerckhoff \cite{Hodgson-Kerckhoff2005} 
and Futer-Kalfagianni-Purcell \cite{FKP1}.

\medskip

Gromov and Thurston obtained the following result,
by constructing a 
Riemannian metric of negative curvature on $M(\nub)$,
when each surgery curve is \lq\lq sufficiently long'',
by modifying the complete hyperbolic metric of $\interior M$
(see \cite{Bleiler-Hodgson1996} for a detailed proof).

\begin{theorem}
[The 2$\pi$-Theorem]\index{2$\pi$-theorem}
Let $M$ be an orientable complete hyperbolic $3$-manifold of finite volume,
and let $C_1,\dots, C_m$ be disjoint torus cusps of $M$.
Suppose $\nu_i$ is a slope on $\partial C_i$
represented by a geodesic with length $> 2\pi$
with respect to the Euclidean metric.
Then $M(\nu_1,\dots, \nu_m)$ has a Riemannian metric of negative curvature.
\end{theorem}

The metric on 
$M(\nu_1,\dots, \nu_m)$ outside the filling solid tori
is identical to the hyperbolic metric on $M-\cup_{j=1}^m C_j$.
The geometrization theorem 
(Theorem \ref{thm:geometrization-mfd})
established by Perelman guarantees that
the resulting manifold $M(\nu_1,\dots, \nu_m)$
is actually hyperbolic.

The 2$\pi$-theorem was refined to
the $6$-theorem by Agol \cite{Agol2000} and Lackenby \cite{Lackenby2},
and it plays
a key role in the study of exceptional surgeries
(see the next subsection).

\subsection{Exceptional surgeries}
For a given hyperbolic knot $K$ in $S^3$,
or more generally an orientable complete hyperbolic manifold with one cusp,
there are only finitely many exceptional slopes $\nu$
which produce non-hyperbolic manifolds.
For example, the figure-eight knot has $10$ exceptional slopes
(Example \ref{ex:exceptional-slope}(1)).
In the survey \cite{Gordon1998},
Gordon proposed various interesting conjectures,
including one which says that $10$ is the largest possible number of exceptional slopes
of a hyperbolic knot complement.

The natural and important problem of determining exceptional surgery slopes
has attracted attention of many mathematicians,
and an enormous amount of research grew out of this problem,
including:

\begin{itemize}
\item[$\circ$]
the $2\pi$-theorem of Gromov-Thurston \cite{Gromov-Thurston1987}
and its improvement to the $6$-theorem by Agol \cite{Agol2000} and Lackenby \cite{Lackenby2},
\item[$\circ$]
the cyclic surgery theorem by Culler-Gordon-Luecke-Shalen \cite{CGLS},
obtained by combining two different kinds of arguments, namely
(i) arguments using the $\SL(2,\CC)$-character varieties 
(cf. Subsection \ref{subsec:Culler-Shalen-theory}) and 
(ii) combinatorial, graph-theoretic analysis of 
the intersection of two incompressible, planar surfaces in knot exteriors,
\item[$\circ$]
the proof of the Property P conjecture by
Kronheimer--Mrowka \cite{Kronheimer-Mrowka2004},
by using Seiberg-Whitten theory.
\item[$\circ$]
study of finite surgery by Boyer-Zhang \cite{Boyer-Zhang1, Boyer-Zhang2}
and Ni-Zhang \cite{Ni-Zhang2018}, by mainly
using the $\SL(2,\CC)$-character varieties
(Heegaard Floer homology and the Casson-Walker invariant
are also used in \cite{Ni-Zhang2018}),
\item[$\circ$]
the proof of the Property R conjecture by Gabai \cite{Gabai1987}, by using taut foliations,
\item[$\circ$]
a universal upper bound of the number of exceptional slopes by
Hodgson-Kerckhoff \cite{Hodgson-Kerckhoff2005, Hodgson-Kerckhoff2008},
by developing deformation theory of hyperbolic structures
(cf. \cite{Hodgson-Kerckhoff1998}).
\item[$\circ$]
the optimal universal upper bound, $10$, on the number of exceptional slopes
of a one-cusped hyperbolic manifold by Lackenby-Meyerhoff
\cite{Lackenby-Meyerhoff2013}
(see Agol \cite{Agol2010a} for related work),
\item[$\circ$]
the optimal universal upper bound, $8$, on the geometric intersection numbers
of pairs of exceptional slopes of one-cusped hyperbolic manifolds by Lackenby-Meyerhoff
\cite{Lackenby-Meyerhoff2013},
\item[$\circ$]
the complete classification of exceptional surgeries on hyperbolic alternating knots
by Ichihara-Masai \cite{Ichihara-Masai2016},
building on a result of \cite{Lackenby1}
and through computer-aided verified computation \cite{HIKMOT}
using a super-computer.
\end{itemize}
The last three results give affirmative answers to some conjectures in \cite{Gordon1998}.
See the survey articles \cite{Boyer, Gordon1998, Gordon2012} for 
background and further information.

As for {\it Seifert surgeries}\index{Seifert surgery} of knots,
namely surgeries which produce Seifert fibered spaces,
Deruelle, Miyazaki and Motegi \cite{DKM2012}
embarked on the project to understand 
the whole shape of relationships among all such surgeries,
and various interesting results are obtained in this direction.

Among Seifert surgeries,
{\it lens space surgeries}\index{lens space surgery} are particularly interesting.
Berge \cite{Berge*} presented a conjecturally complete list 
of lens space surgery on knots in $S^3$.
Based on Berge's conjecture, 
Goda and Teragaito \cite{Goda-Teragaito2000} conjectured that
if a $p$-surgery on a hyperbolic knot $K$ produces a lens space
then $K$ is a fibered knot and its genus $g$ satisfies the inequality
$2g+8 \le |p|\le 4g-1$.
(Note that by the cyclic surgery theorem $p$ is an integer.)
Rasmussen \cite{Rasmussen2004} attacked this problem by using 
the Heegaard Floer homology,
and obtained the estimate $|p|\le 4g+3$.
This in fact relies on the fact that 
lens spaces belong to larger class of spaces, known as {\em L-spaces}\index{L-space},
which are rational homology $3$-spheres 
with the \lq\lq simplest Heegaard-Floer homology'' 
(see Ozsv\'ath-Szab\'o \cite{Ozsvath-Szabo2005}).
See Greene \cite{Greene2015} and references therein for further information
on L-space surgery,
and see the reviews \cite{Ozsvath-Szabo2006, Juhasz2015} for the background. 

A nice overall survey (in Japanese) on surgery was recently written by Motegi \cite{Motegi*},
and its English translation will appear soon.
This survey is strongly recommended.

\section{Volumes of hyperbolic $3$-manifolds}
\label{sec:volume}
The volume is the most basic invariant
of hyperbolic manifolds.
After quickly recalling a method for calculating hyperbolic volumes,
we explain (i) the Jorgensen-Thurston theory concerning the volume spectrum 
of hyperbolic $3$-manifolds,
(ii) results concerning small volume hyperbolic manifolds,
(iii) relation to the Gromov norm,
and finally (iv) the volume conjecture,
which lies in the two innovations,
hyperbolic geometry and quantum topology, in knot theory.

\subsection{Calculation of hyperbolic volumes}
We explain a method for calculating the volumes of hyperbolic $3$-manifolds,
which is implemented in SnapPea.
The method depends on the fact that
every hyperbolic $3$-manifold $M$ is obtained by hyperbolic Dehn filling
on a cusped hyperbolic manifold, say $M_0$.
This follows from the facts that 
the complement of a simple closed geodesic is a cusped hyperbolic manifold 
(see \cite{Sakai2})
and that the shortest closed geodesic in $M$ is simple.
SnapPea usually succeeds in finding an ideal triangulation
of the complete hyperbolic manifold $M_0$,
which can be deformed into an ideal triangulation 
of the incomplete hyperbolic structure on $M_0$
whose completion yields the complete hyperbolic structure of $M$
(cf. Subsection \ref{subsec:proof-HDFT}).
Thus the calculation of $\vol(M)$ is reduced to that
of the volumes of ideal tetrahedra.

Recall that the isometry type of an ideal tetrahedron
is determined by its shape parameter $z\in\HH\subset \CC$,
which in turn represent the similarity class of the
Euclidean triangle with vertex set $\{0,1,z\}$.  
Let $\alpha,\beta,\gamma$ be the inner angles of this triangle.
Then the volume of the ideal tetrahedron $\Delta(z)$ of shape parameter $z$
is given by the following formula:
\[
\vol(\Delta(z))=\Lambda(\alpha)+\Lambda(\beta)+\Lambda(\gamma),
\]
where $\Lambda(\theta)$ is the {\it Lobachevsky function}\index{Lobachevsky function} defined by
\[
\Lambda(\theta)
=
-\int_0^{\theta}\log|2\sin t|dt
=
\frac{1}{2}\sum_{n=1}^{\infty}\frac{\sin(2n\theta)}{n^2}.
\]
The volume function $\vol(\Delta(z))$ takes the maximal value
$V_{\mathrm{tet}}=3\Lambda(\pi/3)=1.0149416...$
precisely at $z=\exp(\pi\sqrt{-1}/3)$,
i.e., exactly when $\Delta(z)$ is a regular ideal tetrahedron.
See \cite[Chapters 6 and 7]{Thurston0} for details.

\label{subsec:volume}
\subsection{J\o rgensen--Thurston theory for the volumes of hyperbolic $3$-manifolds}
Let $V_n\subset \RR_+$ be the ordered set 
consisting of the volumes of complete hyperbolic $n$-manifolds.
If $n\ne 3$ then $V_n$ is isomorphic to $\NN$,
by Gauss--Bonnet theorem for $n=2$ and by
Wang's theorem \cite{Wang1972} for $n\ge 4$.
For dimension $n=3$, we have the following surprising theorem
due to J\o rgensen and Thurston (see \cite{Thurston0}),
which forms a sharp contrast to Wang's theorem.

\begin{theorem}
[J\o rgensen--Thurston Theorem]\index{J\o rgensen--Thurston theorem}
\label{thm:Jorgensen-Thurston}
$V_3$ is a well-orderd closed set which is isomorphic to $\omega^{\omega}$. 
Moreover, the map 
\[
\vol:\{\mbox{\rm complete hyperbolic $3$-manifolds of finite volume}\}/({\rm isometry})
\to V_3
\]
is finite to one.
\end{theorem}
This means that there is 
a smallest volume $v_1$, a next smallest volume $v_2$,
and so on,
and these are the volumes of closed hyperbolic $3$-manifolds.
The increasing sequence $v_1<v_2<\cdots< v_k < \cdots$ has a limit $v_{\omega}$,
and this is the volume of a complete hyperbolic $3$-manifold with one cusp
(cf. Subsection \ref{subsec:Geometry-Dehn-filling}).
There is a smallest volume $v_{\omega+1}$ bigger than $v_{\omega}$,
a second smallest volume $v_{\omega+2}$ bigger than $v_{\omega+1}$, and so on,
and these are the volumes of closed hyperbolic $3$-manifolds,
and their limit $v_{2\omega}$ 
is the second smallest volume of a complete hyperbolic $3$-manifold with one cusp.
The increasing sequence $v_{\omega}<v_{2\omega}<\cdots< v_{k\omega} < \cdots$ has a limit $v_{\omega^2}$,
and this is the volume of a complete hyperbolic $3$-manifold with two cusps, and so on.

The second statement of Theorem \ref{thm:Jorgensen-Thurston}
says that the volume is \lq\lq almost'' a complete invariant of complete hyperbolic manifolds.

Of course, the volume is not a complete invariant.
For example, the complements of the Kinoshita--Terasaka knot and the Conway knot
have the same volume 11.21911773....
In fact, Ruberman \cite{Ruberman} proved
that the hyperbolic volume, 
more generally the Gromov invariant
(cf. Subsection \ref{subsec:Gromov-norm} below),
is unchanged by mutation.
Hodgson and Masai \cite{Hodgson-Masai2013} studied 
the number $N(v)$ of
orientable hyperbolic $3$-manifolds with given volume $v\in V_3$:
they constructed infinitely many $v\in V_3$
for which $N(v)=1$,
and proved the exponential growth of $N(v)$
by showing $N(4nV_{\mathrm{oct}})\ge 2^n/(2n)$.
See Chesebro-DeBlois \cite{Chesebra-Debois2014} and 
Millichap \cite{Millichap2015} for related results.

\subsection{Small volume hyperbolic manifolds}
\label{subsec:small-volume}
It is a natural and important problem to determine the small volumes,
such as $v_1$, $v_{\omega}$, $v_{\omega^2}$, etc.
For the minimal volume $v_{\omega^n}$
of orientable complete hyperbolic $3$-manifolds with $n$-cusps,
the following results are established.
\begin{itemize}
\item[$\circ$]
Gabai-Meyerhoff-Milley (2009) \cite{Gabai-Meyerhoff-Milley_2009}:
The {\it Fomenko-Matveev-Weeks manifold}\index{Fomenko-Matveev-Weeks manifold},
which is obtained by $(5,2)$ and $(5,1)$ Dehn surgery
on the Whitehaed link,
has the smallest volume
$v_1=0.94270736...$.
\item[$\circ$]
Cao-Meyerhoff (2001) \cite{Cao-Meyerhoff}: 
The figure-eight knot complement and its sister,
namely $(5,1)$-filling on
one component of the Whitehead link complement, 
have the volume
$v_{\omega}=2V_{\mathrm{tet}}=2.02988...$,
where $V_{\mathrm{tet}}=1.0149416...$ is the volume of 
the regular ideal tetrahedron.
The figure-eight knot is the orientation double cover
of the Gieseking manifold,
the non-orientable hyperbolic $3$-manifold,
which has the smallest volume among the all 
(orientable or not) complete non-compact hyperbolic $3$-manifolds
(see Adams \cite{Adams1987}).
\item[$\circ$]
Agol (2010) \cite{Agol2010b}:
The Whitehaed link complement and the complement of the
pretzel link $P(-2,3,8)$ have the volume 
$v_{\omega^2}=V_{\mathrm{oct}}=3.66386...$,
where $V_{\mathrm{oct}}$ is the volume of regular ideal octahedron.
\item[$\circ$]
Yoshida (2013) \cite{Yoshida2013}:
The complement of the minimally twisted hyperbolic $4$-chain link has the volume
$v_{\omega^4}=2V_{\mathrm{oct}}=7.32772...$.
\end{itemize}

See the review \cite{Gabai-Meyerhoff-Milley_2014},
for further information.
It should be noted that all of the above small volume hyperbolic manifolds
are arithmetic (cf. \cite{Hatcher1983}, \cite[Theorem 5.1]{Neumann-Reid1992}
and Subection \ref{subsec:arithmetic-group}).

As is noted in \cite[Introduction]{Gabai-Meyerhoff-Milley_2011},
Thurston had long promoted the idea that volume is a good measure of the complexity of a hyperbolic $3$-manifold.
In fact, in \cite[the end of Chapter 6]{Thurston0}, 
he writes the following:
{\it One gets a feeling that volume is a very good measure of the complexity of a link complement, and that the ordinal structure is really inherent in three-manifolds.}
The following conjecture,
due to Thurston, Weeks, Matveev-Fomenko and Mednykh-Vesnin,
states the idea more rigorously,
and the results presented above can be regarded as
partial answers to this conjecture.

\begin{conjecture}
The complete low-volume hyperbolic $3$-manifolds can be obtained by filling cusped hyperbolic $3$-manifolds of small topological complexity.
\end{conjecture}

To end this subsection, we explain another approach 
to Thurston's idea above,
by using the notions of {\it shadows}\index{shadow} of $3$ and $4$-manifolds introduced 
by Turaev \cite{Turaev1, Turaev2}.
Costantino and Thurston \cite{Costantino-Thurston} 
introduced the {\it shadow complexity}\index{shadow complexity} $\mathrm{sc}(M)$ of a compact orientable
$3$-manifold $M$ with (possibly empty) toral boundary,
and proved the following estimate of the Gromov norm $||M||$
(cf. Subsection \ref{subsec:Gromov-norm}, below):
\[
\frac{V_{\mathrm{tet}}}{2V_{\mathrm{oct}}}||M||
\le \mathrm{sc}(M)
\le C||M||^2
\quad
\mbox{for some universal constant $C$.}
\] 
In the same paper, they implicitly introduced the notion of 
{\it stable map complexity} and
studied its relation between (branched) shadow complexity as well.
Ishikawa and Koda \cite{Ishikawa-Koda} showed the two complexities are actually equal, and 
moreover, using the result of \cite{FKP1},
they gave an elaborate refinement of the above (left) inequality when $M$ is hyperbolic.
They also defined the {\it branched shadow complexity}\index{shadow complexity!branched} $\mathrm{bsc}(M,L)$
for a link $L$ in a compact orientable
$3$-manifold $M$ with (possibly empty) toral boundary,
and gave a complete characterization of hyperbolic links $L$ in $S^3$
with $\mathrm{bsc}(S^3,L)=1$.

\subsection{Gromov norm}
\label{subsec:Gromov-norm}
In \cite{Gromov_1982}, Gromov introduced the notion of {\it simplicial volume} 
$||M||$ of a closed manifold $M$ as follows,
using real singular homology:
\[
||M||:=\inf\{||z|| \ | \ \mbox{$z$ is a singular cycle 
representing the fundamental class $[M]$}\}
\]
Here, for a (real) singular chain $z=\sum_j a_j\sigma_j$,
its norm $||z||$ is defined as the sum  $\sum_j |a_j|$ of the absolute values of its coefficients.
He used it to estimate the \lq\lq{\it minimal volume}''
of closed smooth manifold (see \cite{Gromov_1982}).
Building on this work,
Thurston \cite[Chapter 6]{Thurston0} defined the {\it Gromov norm}\index{Gromov norm} $||M||$
of a compact orientable $3$-manifold $M$ with (possibly empty) toral boundary 
as follows:
\[
||M||:=
\lim_{\epsilon\to 0}
\inf\{||z|| \ | \ \mbox{$z$ is a singular chain 
representing $[M,\partial M]$ and $||\partial z||\leq\epsilon$}\}
\]
He then proved the following.
\begin{enumerate}
\item
If $M$ is hyperbolic (and hence $\interior M$ admits a complete hyperbolic structure of finite volume), 
then 
\[
||M||=\frac{1}{V_{\mathrm{tet}}}\vol(\interior M).
\]
\item
If $M$ is a Seifert fibered space, then $||M||=0$.
\item
Let $T$ be a torus embedded in $\interior M$ and let $M_T$ be the
manifold obtained by cutting $M$ along $T$. Then $||M||\leq||M_T||$.
\end{enumerate}
Soma \cite{Soma} proved that 
when $T$ is incompressible, equality holds in (3)
and that similarly equality holds 
for an incompressible annulus properly embedded in $M$. 
He then defined, for a link $L$ in $S^3$,
the {\it Gromov invariant}\index{Gromov invariant}\index{link!Gromov invariant} $||L||$ of $L$  by
$
||L||=||E(L)||
$,
and obtained the following theorem.

\begin{theorem}[Soma]
\label{thm:Soma}
For a link $L$ in $S^3$, the following hold.
\begin{enumerate}
\item
If $L$ is a split sum of two links $L_1$ and $L_2$,
then $||L||=||L_1||+||L_2||$.
\item
If $L$ is a connected sum of two links $L_1$ and $L_2$,
then $||L||=||L_1||+||L_2||$.
\item
Suppose $L$ is a non-splittable link,
and let $\{M_j\}$ be the hyperbolic pieces of
the JSJ decomposition of $E(L)$.
Then
\[
||L||=\sum_j ||M_j||=\frac{1}{V_{\mathrm{tet}}}\sum_j\vol(\interior M_j).
\]
\end{enumerate}
\end{theorem}

\subsection{The volume conjecture}
In addition to the revolution caused by William Thurston,
knot theory has experienced yet another revolution
through the discovery of the Jones polynomial by Vaughan Jones \cite{Jones}.
The Volume Conjecture, first stated by Rinat Kashaev \cite{Kashaev} 
and then reformulated and expanded by Hitoshi Murakami and Jun Murakami \cite{Murakami-Murakami},
provoked deep interaction between the two innovations,
hyperbolic geometry and quantum topology.

The conjecture says that the hyperbolic volume of a hyperbolic knot in $S^3$
(more generally, the Gromov norm of a knot in $S^3$)
is determined by the asymptotic behavior of {\it Kashaev's invariant} $\langle K\rangle_N$,
which is shown by \cite{Murakami-Murakami} to coincide with
the evaluation, $J_N(K)$, of the {\it $N$-colored Jones polynomial} 
(with a certain normalization)
at the primitive $N$-th root of unity $\exp(2\pi i/N)$.

\begin{conjecture}[Volume Conjecture]\index{volume conjecture}
For any knot $K$ in $S^3$, the following holds:
\[
||K||=\frac{2\pi}{V_{\mathrm{tet}}}\lim_{N\to\infty}\frac{\log|J_N(K)|}{N}.
\]
In particular, if $K$ is a hyperbolic knot,  the following holds:
\[
\vol(S^3-K)=2\pi\lim_{N\to\infty}\frac{\log|J_N(K)|}{N}.
\]

\end{conjecture}

Moreover, H. Murakami and J. Murakami proved that
Kashaev's invariant also coincides with an evaluation 
of the {\it generalized Alexander polynomial}
defined by Y. Akutsu, T. Deguchi and T. Ohtsuki \cite{ADO}.
They say in \cite[page 86]{Murakami-Murakami} that
{\it the set of colored Jones polynomials and the set of generalized Alexander polynomials of Akutsu-Deguchi-Ohtsuki intersect at Kashaev's invariants.}

Furthermore, 
H. Murakami, J. Murakami, M. Okamoto, T. Takata and Y. Yokota \cite{MMOTY}
proposed the following complexification of Kashaev's conjecture:

\begin{conjecture}
[Complexification of the Volume Conjecture]
For any hyperbolic knot $K$ in $S^3$, the following holds:
\[
\vol(S^3-K)+\sqrt{-1}\mathrm{CS}(S^3-K)=2\pi\lim_{N\to\infty}\frac{\log{J_N(K)}}{N}.
\]
\end{conjecture}

In the above conjecture $\mathrm{CS}(S^3-K)$ denotes the {\it Chern-Simons invariant} of 
$S^3-K$
(see \cite{Chern-Simons, Meyerhoff, YoshidaT}).
For further information, see the surveys \cite{MurakamiH_1, MurakamiH_2}
and the recently published book \cite{Murakami-Yokota2018}.

\section{Commensurability and arithmetic invariants of hyperbolic manifolds}
\label{sec:arithmetic-invariant}
In \cite[Sections 6.7 and 6.8]{Thurston0},
Thurston studied the commensurability 
relation among hyperbolic knot/link complements,
and gave various commensurable and incommesurable examples.
This work has promoted intimate interaction between knot theory and number theory.
In this section, we recall basic arithmetic invariants of 
commensurability classes of Kleinian groups, 
and describe application to knot theory.
We also describe the dichotomy between arithmetic groups and non-arithmetic groups
found by Margulis and Borel.
In the final subsection, we 
recall the solution due to Gehring, Marshal and Martin
of the $3$-dimensional Siegel problem 
to determine the minimal volume of hyperbolic orbifolds,
lying emphasis on the role of arithmetic groups.
For further information on 
the topic of this section,
see the textbook Maclachlan-Reid \cite{Maclachlan-Reid}.

\subsection{Commensurability classes and invariant trace fields}
Two Kleinian groups $\Gamma_1$ and $\Gamma_2$ are said to 
be {\it commensurable}\index{commensurable}
if there is a conjugate, $\Gamma_2^{g}:=g^{-1}\Gamma_2g$ ($g\in\PSL(2,\CC)$)
such that $\Gamma_1\cap\Gamma_2^{g}$ has finite index both in
$\Gamma_1$ and $\Gamma_2^{g}$.
This is equivalent to the condition that
the two hyperbolic manifolds
$M_1=\HH^3/\Gamma_1$ and $M_2=\HH^3/\Gamma_2$ are {\it commensurable},
i.e., 
there is a hyperbolic manifold which is a finite covering of both
$M_1$ and $M_2$.
As is explained in Subsection \ref{subsec:canonical},
the canonical decomposition provides us an efficient method for checking if
two (cusped) hyperbolic manifolds are isometric.
But, the method is not directly applicable for checking commensurability,
though there is a nice application of the canonical decomposition 
for the commensurability problem (see \cite{GHH2008} and Subsection \ref{subsec:arithmetic-group}).

Number theory enables us to define a very useful invariant of the commensurability
classes of Kleinian groups of cofinite volume.
Let $M=\HH^3/\Gamma$ be an orientable complete hyperbolic manifold of finite volume.
Consider the set $\tr\Gamma = \{\pm \tr(\gamma) \ | \ \gamma\in\Gamma\}\subset \CC$
and the field $\QQ(\tr\Gamma)$ generated by the set.
(Note that the trace $\tr \gamma$ for $\gamma\in \PSL(2,\CC)\cong \Isom^+(\HH^3)$ is well-defined up to sign.)
This is called the {\it trace field}\index{trace field}\index{Kleinian group!trace field} of the Kleinian group $\Gamma$.
It follows from the rigidity theorem
that the trace field $\QQ(\tr\Gamma)$ has finite degree over $\QQ$,
i.e., it is a {\it number field}.
By the rigidity theorem again,
this is an invariant of the topological space $M$.

Though the trace field $\QQ(\tr\Gamma)$ itself is not, in general, an invariant of 
the commensurability class,
it provides us with a very useful commensurability invariant as follows.
Let $\Gamma^{(2)}$ be the subgroup of $\Gamma$ generated by
$\{\gamma^2 \ | \ \gamma\in\Gamma\}$.
Then $\Gamma^{(2)}$ is normal in $\Gamma$
and $\Gamma/\Gamma^{(2)}$ is a finite abelian group which is a direct sum of
order $2$ cyclic groups.
The following theorem was proved by Reid \cite{Reid1990}.

\begin{theorem}
\label{thm:invariant-trace-field}
Let $\Gamma$ be a Kleinian group of finite covolume.
Then $\QQ(\tr\Gamma^{(2)})$ is an invariant of the commensurability class of $\Gamma$.
Moreover
\[
\QQ(\tr\Gamma^{(2)})=\QQ(\{(\tr \gamma)^2 \ | \ \gamma\in\Gamma\}).
\]
\end{theorem}
\noindent
The field $\QQ(\tr\Gamma^{(2)})$ is denoted by $k(\Gamma)$ and is called 
the {\it invariant trace field}\index{invariant trace field}\index{Kleinian group!invariant trace field}  of $\Gamma$.
By \cite[Corollary 2.3]{Neumann-Reid1992},
if $M=\HH^3/\Gamma$ is a knot complement
(or more generally, the complement of a link in a $\ZZ/2\ZZ$-homology sphere)
then $k(\Gamma)=\QQ(\tr\Gamma)$:
thus in this case the trace field itself is an invariant of 
the commensurability class.

If $M$ is a cusped hyperbolic manifold which admits an ideal triangulation 
into the hyperbolic ideal tetrahedra $\{\Delta(z_1),\dots,\Delta(z_t)\}$,
then the following holds \cite[Theorem 2.4]{Neumann-Reid1992}:
\[
k(\Gamma)=\QQ(z_1,\dots, z_t).
\]

The {\it invariant quaternion algebra}\index{invariant quaternion algebra}\index{Kleinian group!invariant quaternion algebra}
of $\Gamma$ is the $k(\Gamma)$-algebra of the $2\times 2$ matrix algebra $M_2(\CC)$
generated over $k(\Gamma)$ by the elements of $\Gamma^{(2)}$.
It is denoted by $A(\Gamma)$.
This algebra is also an invariant of the commensurability class of $\Gamma$.
Both $k(\Gamma)$ and $A(\Gamma)$ are preserved by mutation 
(see \cite{Neumann-Reid1991}).

The computer program \lq\lq Snap'' 
calculates various arithmetic invariants 
including the invariant trace field
and the invariant quaternion algebra (see \cite{CGHW}).

\subsection{Commensurators and hidden symmetries}

For a Kleinian group $\Gamma$ of cofinite volume,
the {\it commensurator}\index{commensurator}\index{Kleinian group!commensurator} of $\Gamma$ is defined by
\[
\Comm(\Gamma)=\{g\in \Isom\HH^3 \ | \ [\Gamma; \Gamma\cap\Gamma^{g}]<\infty\},
\]
and its orientation-preserving subgroup is denoted by $\Comm^+(\Gamma)$.
The commensurator $\Comm(\Gamma)$ is identified with the group
of equivalence classes of virtual automorphisms of $\Gamma$.
A {\it virtual automorphism}\index{virtual automorphism}\index{Kleinian group!virtual automorphism}
of $\Gamma$ is an isomorphism 
$\phi:\Gamma_1\to\Gamma_2$ between subgroups of finite index in 
$\Gamma$, and two virtual automorphisms are defined to be equivalent 
if they agree on some subgroup of $\Gamma$ of finite index.
A virtual automorphism represents an isometry 
between two finite coverings $\HH^3/\Gamma_1$ and $\HH^3/\Gamma_2$
of the hyperbolic manifold $M=\HH^3/\Gamma$.
It is called a {\it hidden symmetry}\index{hidden symmetry} of $M$
if it is not a lift of an isometry of $M$.
By a {\it hidden symmetry} of a hyperbolic knot in $S^3$,
we mean a hidden symmetry of the knot complement.
We can see  as follows
that the figure-eight knot $K$ has a hidden symmetry.
Recall that $S^3-K=\HH^3/\Gamma$
has an ideal triangulation consisting of two copies
of the regular ideal tetrahedron $\Delta(\omega)$
with $\omega =\exp(\frac{\pi\sqrt{-1}}{3})$.
This implies that the invariant trace field $k(\Gamma)$
is equal to $\QQ(\omega)=\QQ(\sqrt{-3})$.
Moreover, we see that $\Gamma$ is a subgroup $\PSL(2,\OO_3)$ 
of finite index (actually equal to $24$),
where $\OO_3$ is the ring of integers of the number field of $\QQ(\sqrt{-3})$.
This implies that $\PGL(2,\QQ(\sqrt{-3}))$ belongs to the commensurator 
subgroup of $\Gamma$.
In fact, we have
$\Comm^+(\Gamma)=\PGL(2,\QQ(\sqrt{-3}))$.
Since $\PGL(2,\QQ(\sqrt{-3}))$ is dense in $\PSL(2,\CC)$,
the normalizer of $\Gamma$ must be a proper subgroup of
$\Comm^+(\Gamma)=\PGL(2,\QQ(\sqrt{-3}))$.
Hence $\Gamma$ (and so the figure-eight knot) has 
a hidden symmetry.

In addition to the figure-eight knot,
the two dodecahedral knots of Aitchison and Rubinstein \cite{AR1992}
admit hidden symmetries, and these three are the only known such knots.
Neumann and Reid \cite[Question 1]{Neumann-Reid1992}
conjecture that they are all.
For results related to the conjecture, 
see \cite{Reid-Walsh2008, BBCW2012, BBCW2015, Millichap-Worden}
and references therein.

\subsection{Arithmetic versus non-arithmetic}
\label{subsec:arithmetic-group}
The above explanation for the existence of hidden symmetries
of the figure-eight knot
is based on the fact that the figure-eight knot group 
belongs to the particularly nice family of Kleinian groups,
called {\it arithmetic groups}\index{arithmetic group}.
For the definition of arithmetic groups, 
see the textbook 
\cite{Reid-Maclachlan2003} or the course notes \cite[Chapter 3, Section 3]{Neumann1999b}.
If we restrict our attention to a cofinite volume Kleinian group $\Gamma$ 
such that $M=\HH^3/\Gamma$ has a cusp,
then $\Gamma$ is arithmetic
if and only if $\Gamma$ is 
conjugate to a subgroup of 
$\PGL(2,\OO_d)$ for some positive integer $d$.
Here $\OO_d$ is the ring of integers of the number field of $\QQ(\sqrt{-d})$.
In this case, we have $k(\Gamma)=\QQ(\sqrt{-d})$
and $A(\Gamma)=M_2(\QQ(\sqrt{-d}))$,
and the invariant trace field $k(\Gamma)$ is the complete commensurability invariant of 
the arithmetic group $\Gamma$.
However, most cusped hyperbolic manifolds are non-arithmetic;
in particular,
the figure-eight knot is the unique hyperbolic knot in $S^3$
whose complement is arithmetic (see Reid \cite{Reid1991}).

Margulis \cite{Margulis1974} (see also Borel \cite{Borel})
establised the following striking dichotomy between the arithemtic Kleinian groups
and non-arithmetic Kleinian groups.

\begin{theorem}
\label{thm:arithemtic-nonarithemtic}
Let $\Gamma$ be a cofinite volume Kleinian group.
Then the following hold.
\begin{enumerate}
\item
$\Gamma$ is non-arithmetic if and only if $\Gamma$ has finite index in
$\Comm^+(\Gamma)$.
In this case, $\Comm^+(\Gamma)$ is the unique maximal element
in the commensurability class of $\Gamma$.
\item
$\Gamma$ is arithmetic if and only if 
$\Comm^+(\Gamma)$ is dense in $\PSL(2,\CC)$.
In this case, there are infinitely many maximal elements
in the commensurability class of $\Gamma$.
\end{enumerate}
\end{theorem}

The first assertion of Theorem \ref{thm:arithemtic-nonarithemtic} 
shows that the commensurability class of
a non-arithmetic 
cofinite volume Kleinian group $\Gamma$
is particularly simple,
namely it consists only of conjugates
of finite index subgroups 
of the Kleinian group $\Comm^+(\Gamma)$.
In terms of orbifolds,
this means that two non-arithmetic orientable hyperbolic $3$-manifolds
$M_1$ and $M_2$ are commensurable if and only if
they cover a common orbifold.
Based on this fact and by using the Epstein--Penner decomposition 
\cite{Epstein-Penner}
and finiteness of Epstein--Penner decompositions of a given
cusped hyperbolic manifolds (see Akiyoshi \cite{Akiyoshi}),
Goodman-Heard-Hodgson \cite{GHH2008} gave a practical algorithm
for determining when two cusped hyperbolic 
non-arithmetic $3$-manifolds are commensurable.
Their algorithm is based on the fact that
two cusped hyperbolic $n$-manifolds $M$ and $M'$ cover a common orbifold if and only if
they admit Epstein--Penner decompositions lifting to isometric tilings of $\HH^n$
(see \cite[Theorem 2.4]{GHH2008}).
Their algorithm is implemented in a computer program,
which enabled them to determine the commensurability classes
of the complements of all hyperbolic knots and links up to
$12$ crossings.
In particular, they have shown that
the complements of
the Kinoshita--Terasaka knot and the Conway knot
belong to different commensurability classes,
even though they share the same invariant trace fields
and invariant quaternion algebras.
See Chesebro-DeBlois \cite{Chesebra-Debois2014}
and Millichap-Worden \cite{Millichap-Worden}
for related works.

The second assertion of Theorem \ref{thm:arithemtic-nonarithemtic} 
shows that the commensurability class of arithmetic Kleinian groups 
is very complicated. 
Walter Neumann 
describes a geometric way of thinking of this situation as follows,
in his course notes \cite[Chapter 3, Section 6]{Neumann1999b}.

{\it
A Kleinian (or Fuchsian) group is the symmetry group of some \lq\lq pattern'' in $\HH^3$
(respectively $\HH^2$).
This pattern might just be a tessellation --- for instance,
a tesselation by fundamental domains,
or it might be an Escher-style drawing.
If one superposes two copies of this pattern,
displaced with respect to each other,
one will usually get a pattern which no longer has a Kleinian (or Fuchsian) symmetry group
in our sense --- the symmetry group has become too small to have finite volume quotient.
But in the arithmetic case --- and {\rm only} in this arithmetic case ---
one can always change the displacement very slightly to make the superposed pattern
have a symmetry group
that is of finite index in the original group.}

In the course notes \cite{Neumann1999b},
we can also find a beautiful introduction to 
the idea of Scissor congruence,
with a historical background which goes back to Euclid, Dehn and Hilbert.
For more details of this important topic, see \cite{Neumann1999a}.

\subsection{Siegel's problem and arithmetic manifolds}

In Subsection \ref{subsec:small-volume},
we surveyed various important results concerning small volume hyperbolic $3$-manifolds.
It is equally natural and important to study small volume hyperbolic orbifolds.
In 1943, Siegel \cite{Siegel1943, Siegel1945} 
posed the problem of identifing the infimum
\[
\mu(n)=\inf_{\Gamma}\vol(\HH^n/\Gamma)
\]
where the infimum is taken over the {\it lattices}\index{lattice}
$\Gamma < \Isom^+ \HH^n$,
i.e., discrete subgroups of cofinite volume. 
Siegel solved the problem in dimension $2$,
by showing that
the $(2,3,7)$-triangle group
is the unique Fuchsian group of minimal coarea
\[
\mu(2)=2\pi\left|\frac{1}{2}+\frac{1}{3}+\frac{1}{7}-1\right|=\frac{\pi}{21}.
\]

In 1986, Kazhdan and Margulis \cite{Kazhdan-Marguli1968} made an important contribution
to the Siegel problem, by proving that
$\mu(n)$ is positive and attained for each $n$.

Arithmetic groups play a crucial role in the study of the Siegel problem.
One big reason is that, 
due to formulas of Borel \cite{Borel},
various explicit calculations can be made
for arithmetic Kleinian groups.
According to Gaven Martin \cite{Martin2015}, another reason is that
{\it it turns out that nearly all the extremal problems one might formulate are realised by arithmetic groups, perhaps the number theory forcing additional symmetries in a group and therefore making it \lq\lq smaller'' or \lq\lq tighter''.}

After a long term collaboration,
Gehring, Marshal and Martin \cite{Gehring-Marshall-Martin2009, Gehring-Marshall-Martin2012}
finally solved the $3$-dimensional {\it Siegel problem}\index{Siegel problem}.

\begin{theorem}
The minimum $\mu(3)$ of the volumes of hyperbolic $3$-orbifolds
is 
\[
\mu(3)=\vol (\HH^3/\Gamma_0) = 275^{3/2}2^{-7}\pi^{-6}\zeta_k(2) \sim 0.03905...,
\]
where $\zeta_k$ is the Dedekind zeta function of the underlying number field 
$\QQ(\gamma_0)$, with $\gamma_0$ a complex root of
$\gamma^4+6\gamma^3+12\gamma^2+9\gamma+1$, of discriminant $-275$.
Here $\Gamma_0$ is 
an arithmetic Kleinian group obtained as a $\ZZ/2\ZZ$-extension of the 
index $2$ orientation-preserving subgroup of the group generated by reflection in the faces of the $3$-$5$-$3$-hyperbolic Coxeter tetrahedron.
The group $\Gamma_0$ is generated by two elliptic elements,
one of order $2$ and the other of order $5$.
\end{theorem}

\begin{figure}[ht]
\begin{center}
 {
\includegraphics[height=4cm]{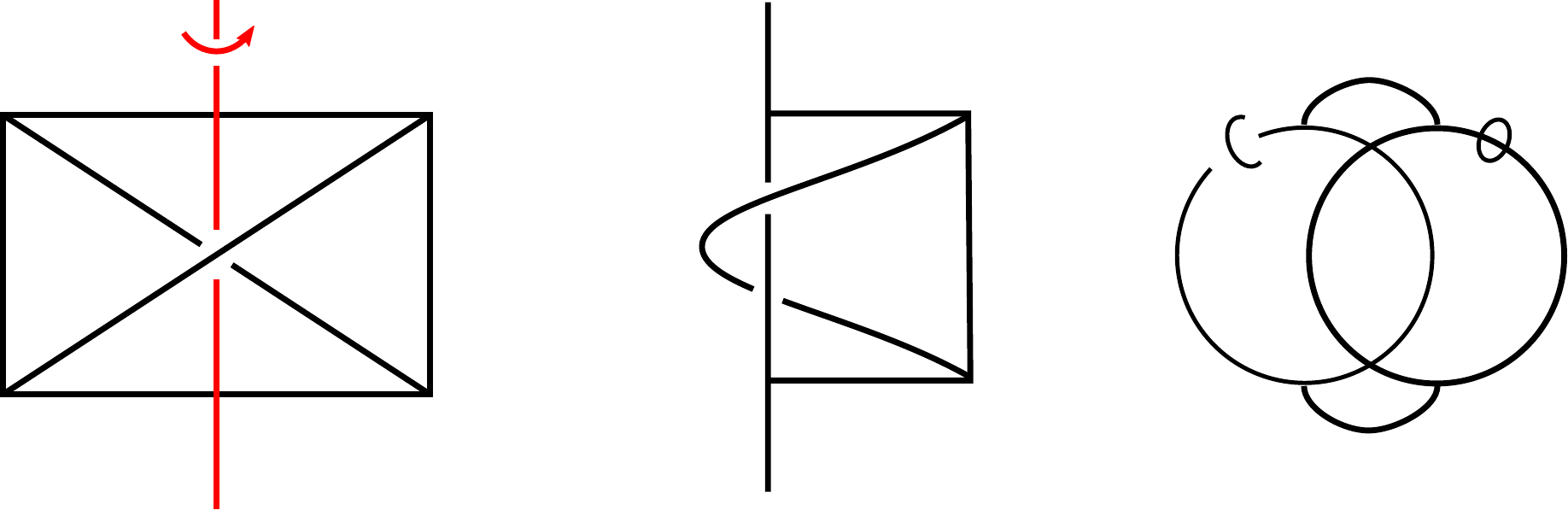}
 }
\end{center}
\caption
{
The minimal volume $3$-orbifold $\HH^3/\Gamma_0$.
The blue eyeglass frame represents the generating pair of $\Gamma_0$
consisting of elliptic elements.
}
\label{fig:minimal-orb}
\end{figure}

\begin{figure}[ht]
\begin{center}
 {
\includegraphics[height=4cm]{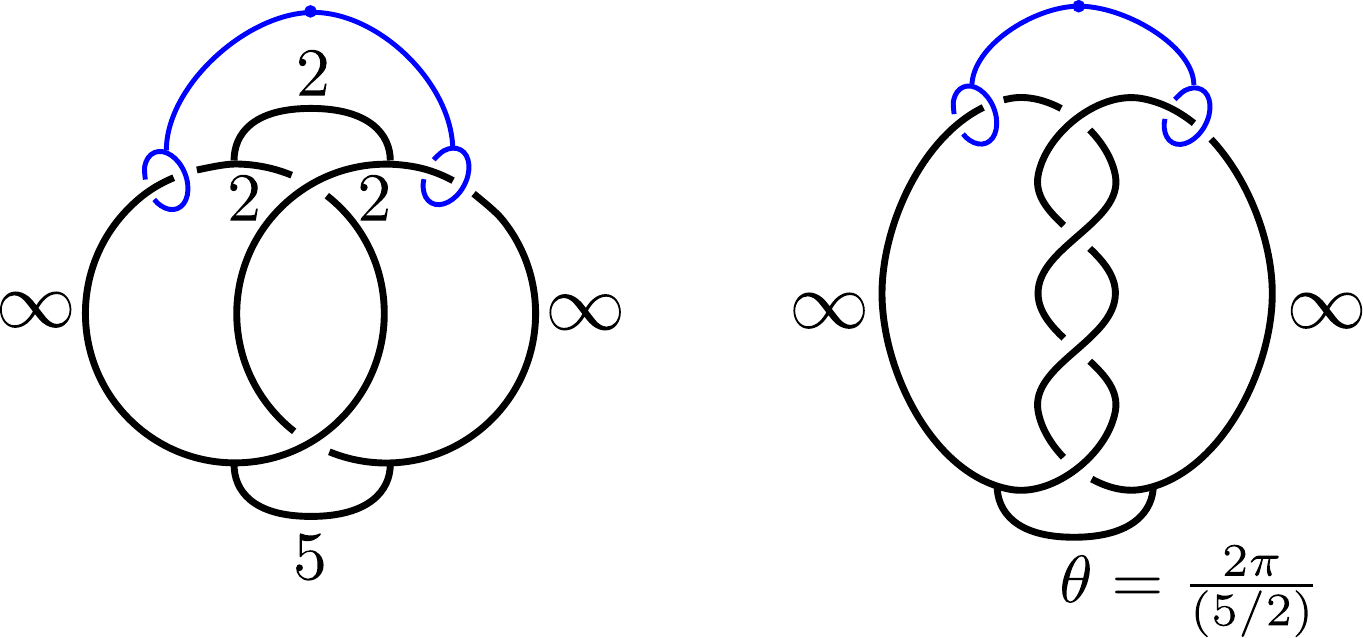}
 }
\end{center}
\caption
{
The left picture illustrates the Heckoid orbifold $H(1/4;5/2)$
and the parabolic generating pair of the Heckoid group $\pi_1^{\mathrm{orb}}(H(1/4;5/2))$.
The Heckoid group is identified with the
image of the holonomy representation of the hyperbolic cone manifold 
$C\left( 1/4;\frac{2\pi}{(5/2)} \right)$ depicted by the right picture.
}
\label{fig:Heckoid}
\end{figure}

\begin{remark}
{\rm
The quotient orbifold $\OO_0=\HH^3/\Gamma_0$ is as illustrated in 
Figure \ref{fig:minimal-orb},
where the blue eyeglasses represent the generating pair.
This orbifold is obtained from the  \lq\lq Heckoid orbifold
$H(1/4;5/2)$'' in Figure \ref{fig:Heckoid}
by an orbifold surgery.
Here a Heckoid orbifold is a hyperbolic $3$-orbifold 
whose orbifold fundamental group is a {\it Heckoid group},
which is a Kleinian group generated by two parabolic transformations
introduced by Riley \cite{Riley1992} as an analogy of Hecke groups
and formulated by \cite{Lee-Sakuma2013}.
Heckoid orbifolds are also intimately related to $2$-bridge links.
As noted by Martin \cite{Martin2015, Martin2016},
most of small volume $3$-orbifolds arise from $2$-bridge links.
}
\end{remark}

See the survey by Martin \cite{Martin2015} for backgrounds and details 
concerning the $3$-dimensional Siegel problem,
and the surveys by 
Belolipetsky \cite{Belolipetsky} and
Kellerhals \cite{Kellerhals2014}
for the higher dimensional Siegel problem.

The above theorem has the following application to finite group actions
on hyperbolic $3$-manifolds.

\begin{corollary}
Let $M$ be an orientable complete hyperbolic $3$-manifold of finite volume,
and let $G$ be a finite group acting on $M$ 
effectively and orientation-preservingly.
Then
\[
|G|\le \frac{\vol M}{\mu(3)}.
\] 
\end{corollary}

A refinement of this corollary for hyperbolic knot complements
can be found in \cite[Theorem 4.14]{FKP_survey}.

\section{Flexibility of complete hyperbolic manifolds - deformation theory of hyperbolic structures -}
\label{sec:flexibility}
Let $M$ be a complete hyperbolic manifold
homeomorphic to the interior of a compact orientable $3$-manifold $\bar M$.
If $\partial\bar M$ is a (possibly empty) union of tori,
then $\vol(M)<\infty$ and so by the Mostow-Prasad rigidity theorem
the complete hyperbolic structure on $M$ is unique. 
However, when $\partial \bar M$ contains a component different from a torus,
the complete hyperbolic structure of $M$ admits a nontrivial deformation,
and there is a rich and deep deformation theory.
This deformation theory is one of the central themes in Kleinian group theory and it
plays a crucial role in the proof of 
the geometrization theorem of Haken manifolds.
In particular, the existence of complete hyperbolic structures
on surface bundles over the circle,
e.g. the complements of hyperbolic fibered knots,
was established as a consequence the double limit theorem
\cite[Theorem 4.1]{Thurston5}
concerning the deformation space of hyperbolic structures 
on $\Sigma\times\RR$ where $\Sigma$ is a (fiber) surface.
The idea of a Cannon--Thurston map,
a $\pi_1(\Sigma)$-equivariant sphere filling curve,
grew out of this construction.

On the other hand, Agol \cite{Agol2011} proved that
a hyperbolic punctured surface bundle over the circle
admits a very special topological ideal triangulation,
called a {\it veering triangulation},
which is canonical in the sense that 
it is determined by the fiber structure.
It was revealed by Gu\'eritaud \cite{Gueritaud2016} that 
the veering triangulation
is intimately related to the Cannon--Thurston map.

The purpose of this section is (i) to give an introduction to
the deformation theory of Kleinian groups
and its relation to the hyperbolic structures of surface bundles over the circle,
and (ii) to explain Cannon--Thurston maps and veering triangulations.
For further information on deformation theory,
see 
Otal \cite{Otal1, Otal2},
Matsuzaki-Taniguchi \cite{Matsuzaki-Taniguchi},
Kapovich \cite{Kapovich2000},
Ohshika \cite{Ohshika2002} and
Marden \cite{Marden1, Marden2}.

\subsection{Convex cores and conformal boundaries of hyperbolic manifolds}
\label{subsec:core-conformal-boundary}
In this subsection, we recall the basic concepts of convex cores and conformal boundaries
of hyperbolic manifolds.

Though the action of a Kleinian group $\Gamma$ on $\HH^3$
is properly discontinuous,
the action of $\Gamma$ on $\partial \HH^3$ does not have this property.
To see this, pick a point $x\in \HH^3$ and consider its orbit $\Gamma x$.
Of course the orbit is discrete in $\HH^3$. But,
it has nonempty accumulation points in the $3$-ball $\HH^3\cup \partial \HH^3$
(provided that $\Gamma$ is not a finite group).
The set of all accumulation points is independent of the choice of $x$
and forms a $\Gamma$-invariant closed set in $\partial \HH^3$.
This set is denoted by $\Lambda(\Gamma)$ and is called the 
{\it limit set}\index{limit set}\index{Kleinian group!limit set} 
of $\Gamma$.
The action of $\Gamma$ on $\Lambda(\Gamma)$ is not properly discontinuous and 
is chaotic.
The complement $\Omega(\Gamma):=\partial \HH^3-\Lambda(\Gamma)$
is called the 
{\it domain of discontinuity}\index{domain of discontinuity}\index{Kleinian group!domain of discontinuity} 
of $\Gamma$,
and it is a (possibly empty) maximal open domain in $\partial \HH^3$
on which $\Gamma$ acts properly discontinuously.

The {\it convex core}\index{convex core} $C_M$ of a complete hyperbolic manifold $M=\HH^3/\Gamma$
is defined as the quotient $C(\Lambda(\Gamma))/\Gamma$,
where $C(\Lambda(\Gamma))$ is the convex hull in $\HH^3$ of the limit set $\Lambda(\Gamma)$.
Note that any closed geodesic in $M$ corresponds to a conjugacy class 
of a hyperbolic element of $\Gamma$ and that the endpoints of its axis 
are contained in $\Lambda(\Gamma)$:
this implies that the axis is contained in $C(\Lambda(\Gamma))$
and so the closed geodesic is contained in $C_M$.
In fact, $C_M$ is the smallest locally convex closed subset of $M$
which contains all closed geodesics of $M$.
The convex core $C_M$ is also characterized as the smallest locally convex submanifold of $M$
whose inclusion is a homotopy equivalence.

On the other hand, since the action of $\Gamma$ on $\partial \HH^3$ (and hence on $\Omega(\Gamma)$)
is conformal, 
the quotient space $\partial_{\infty}M:=\Omega(\Gamma)/\Gamma$ has a natural conformal structure
and forms the boundary of the {\it Klein manifold}\index{Klein manifold}
$(\HH^3\cup \Omega(\Gamma)/\Gamma)$.
The Riemann surface $\partial_{\infty}M=\Omega(\Gamma)/\Gamma$ is called the 
{\it conformal boundary}\index{conformal boundary} 
of $M$.

\begin{example}[Infinite cyclic Kleinian group]
\label{example:limit-set1}
{\rm
For the infinite cyclic Kleinian group $\Gamma$ generated by a hyperbolic transformation
$A(z)=az$ with $|a|\ne 1$
in Example \ref{example:elementary-group}(1),
the convex core of the quotient hyperbolic manifold $\HH^3/\Gamma\cong \interior(D^2 \times S^1)$ 
is equal to the core closed geodesic $(0\times \RR_+)/\Gamma$,
and the conformal boundary is the torus $(\CC-\{0\})/(z\sim az)$.}
\end{example}

In the remainder of this section
$\Sigma\cong \interior\Sigma_{g,b}$
denotes the closed orientable surface of genus $g$ with $b$ punctures,
and with negative Euler characteristic. 

\begin{definition}[Type-preserving representation]
\label{def:type-preserving}
{\rm
A representation $\rho:\pi_1(\Sigma)\cong \pi_1(\Sigma_{g,b}) \to \Isom^+\HH^3$
is {\it type-preserving}\index{type-preserving}\index{representation!type-preserving} if it satisfies the following conditions.
\begin{enumerate}
\item
$\rho$ maps peripheral elements (elements represented by boundary loops
of $\Sigma_{g,b}$)
to parabolic elements.
\item
$\rho$ is 
{\it irreducible}\index{irreducible}\index{representation!irreducible}, 
i.e., $\rho(\pi_1(\Sigma))$ does not have a common fixed point 
on $\partial \HH^3$.
\end{enumerate}
}
\end{definition}

\begin{example}[Fuchsian group]
\label{example:limit-set2}
{\rm
The surface $\Sigma$ admits a complete hyperbolic structure of finite area $\pi|\chi(\Sigma)|$.
Pick a complete hyperbolic metric on $\Sigma$ and
let $\rho_0:\pi_1(\Sigma)\to \Isom^+\HH^2$ be the holonomy representation.
Then it is discrete, faithful and type-preserving, and
its image $\Gamma_0=\rho_0(\pi_1(\Sigma))$ is a Fuchsian group.
The limit set of the Fuchsian group $\Gamma_0$
is equal to $\partial\HH^2$.
Regard $\Gamma_0$ as a Kleinian group,
i.e., a discrete subgroup of $\Isom^+\HH^3$.
Then the limit set $\Lambda(\Gamma_0)$
is the round circle $\partial\HH^2$ in $\partial\HH^3$,
where $\HH^2(=\RR\times \RR_+\subset \CC\times \RR_+=\HH^3)$ 
is the hyperplane of $\HH^3$ invariant by $\Gamma$.
The Kleinian manifold 
$(\HH^3\cup \Omega(\Gamma_0))/\Gamma_0$ is homeomorphic to the product
of $\Sigma$ and the closed interval $[-\infty,\infty]$,
and the convex core is identified with $\Sigma\times 0$.
}
\end{example}

\begin{example}[Quasifuchsian group]
\label{example:limit-set3}
{\rm
The Fuchsian representation $\rho_0:\pi_1(\Sigma) \to \PSL(2,\RR)$ in 
the previous example
admits a nontrivial deformation into a faithful discrete type-preserving
$\PSL(2,\CC)$-representation $\rho$,
such that $(\HH^3\cup\Omega(\Gamma))/\Gamma\cong \Sigma\times [-\infty, \infty]$
where $\Gamma=\rho(\pi_1(\Sigma))$.
This condition is equivalent to the condition that 
the limit set $\Lambda(\Gamma)$ is a topological circle.
A Kleinian group isomorphic to $\pi_1(\Sigma)$
satisfying this condition is called a 
{\it quasifuchsian group}\label{quasifuchsian group}\label{Kleinian group!quasifuchsian} 
and the holonomy representation is called a
{\it quasifuchsian representation}\label{quasifuchsian representation}\label{representation!quasifuchsian}.
Generically, a quasifuchsian group is not conjugate to a Fuchsian group
in $\PSL(2,\CC)$,
and in this case, the circle $\Lambda(\Gamma)$ in $\partial \HH^3$
is very complicated;
in particular its Housdorff dimension is strictly bigger than $1$.
The convex core $C_M$ of the hyperbolic manifold $M=\HH^3/\Gamma$
is identified with $\Sigma\times [-1,1]$ in $\Sigma\times (-\infty,\infty)\cong M$.
Each boundary component $\Sigma\times \{\pm 1\}$ of the convex core
has the structure of \lq\lq hyperbolic surface bent along a geodesic lamination''.
(see \cite[Section 8.5]{Thurston0}, \cite{CEG}).
The domain of discontinuity $\Omega(\Gamma)$ consists of two components
$\Omega_+(\Gamma)$ and $\Omega_-(\Gamma)$,
and the Riemann surfaces $S_{\pm}=\Omega_{\pm}(\Gamma)$
correspond to the boundary components $\Sigma\times\{\pm\}$
of $\Sigma\times[\infty,+\infty]$.
}
\end{example}

\begin{example}[Fiber group]
\label{example:limit-set4}
{\rm
 Let $\hat M=\HH^3/\hat \Gamma$ be a complete hyperbolic manifold of finite volume,
and assume that $\hat M$ has the structure of a $\Sigma$-bundle over $S^1$.
Then the {\it fiber group}\index{fiber group}, $\Gamma$, the subgroup of $\hat \Gamma$
obtained as the image of the fundamental group of a fiber surface $\Sigma$, 
is an infinite normal subgroup.
This implies that $\Lambda(\Gamma)=\Lambda(\hat \Gamma)=\partial \HH^3$
(see \cite[Corollary 8.1.3]{Thurston0}).
Thus the inverse image of a fiber $\Sigma$ in the universal cover $\HH^3$ of $\hat M$
is a topological plane whose closure contains the whole 
ideal boundary $\partial\HH^3$.
It is very difficult to imagine such a plane, and
in this sense,
the fiber group $\Gamma$ is quite different from a quasifucshian group,
though they are all isomorphic to $\pi_1(\Sigma)$.
}
\end{example}

\subsection{Deformation space} 
\label{subsec:deformatin-space}
We continue to denote by $\Sigma$
a closed orientable surface of genus $g$ with $b$ punctures,
which has a negative Euler characteristic. 
By a {\it marked hyperbolic structure}\index{marked hyperbolic structure}\index{hyperbolic structure!marked} on $\Sigma$,
we mean a pair $(S, f)$ of a finite area complete hyperbolic surface $S=\HH^2/\Gamma$
and an orientation-preserving homeomorphism $f:\Sigma \to S$.
Note that the composition of $f_*:\pi_1(\Sigma)\to \pi_1(S)$
and the holonomy representation $\pi_1(S)\to\Gamma<\Isom^+\HH^2$ 
determine a type-preserving discrete faithful representation
$\rho:\pi_1(\Sigma)\to \Isom^+\HH^2$.
Two marked hyperbolic structure $(S_1, f_1)$ and  
$(S_2, f_2)$ on $\Sigma$ are {\it equivalent}
if there is an orientation-preserving isometry 
$h:S_1\to S_2$ such that $h\circ f_1$ is homotopic to $f_2$.
This is equivalent to the condition that the corresponding
representations $\rho_1$ and $\rho_2$ are equal up to conjugation by an element 
of $\Isom^+\HH^2$.
Let $H(\Sigma)$ be the set of all marked hyperbolic structure on $\Sigma$
up to equivalence.

In order to introduce a natural topology on $H(\Sigma)$,
consider the spaces
\begin{align*}
\Hom_{\mathrm{tp}}(\pi_1(\Sigma),\Isom^+\HH^2)
& :=
\{\rho:\pi_1(\Sigma)\to \Isom^+\HH^2 \ | \ \mbox{$\rho$ is type-preserving}\};
\\
\Rep_{\mathrm{tp}}(\Sigma)
& :=
\Hom_{\mathrm{tp}}(\pi_1(\Sigma),\Isom^+\HH^2)/\Isom^+\HH^2.
\end{align*}
By choosing a finite generating set of $\pi_1(\Sigma)$ of cardinality $k$,
$\Hom_{\mathrm{tp}}(\pi_1(\Sigma),\Isom^+\HH^2)$ 
is identified with a subset
of the product (topological) space $(\Isom^+\HH^2)^k$,
and the subspace topology it inherits is independent of the choice of 
a finite generating set.
The group $\Isom^+\HH^2$ acts by conjugation on the space 
$\Hom_{\mathrm{tp}}(\pi_1(\Sigma),\Isom^+\HH^2)$,
and $\Rep_{\mathrm{tp}}(\Sigma)$ is defined to be the quotient space.
The set $H(\Sigma)$ is identified with a subset of $\Rep_{\mathrm{tp}}(\Sigma)$,
and we denote the resulting topological space by $AH(\Sigma)$.

The space $AH(\Sigma)$ is nothing other than 
the {\it Teichm\"uller space}\index{Teichm\"uller space} $\Teich(\Sigma)$ of $\Sigma$.
The {\it Fenchel-Nielsen coordinate}\index{Fenchel-Nielsen coordinate}
gives a homeomorphism from $AH(\Sigma)=\Teich(\Sigma)$ onto 
the Euclidean space $\RR^{6g-6+3b}$
(see \cite{Imayoshi-Taniguchi}, \cite[Theroem 5.3,5]{Thurston0}).
It should be noted that $\Teich(\Sigma)$ can be also identified with
the space of {\it marked Riemann surface structures on $\Sigma$}.

Now we consider hyperbolic structures
on the oriented $3$-manifold $\Sigma\times \RR$.
By a {\it marked hyperbolic structure}\index{marked hyperbolic structure}\index{hyperbolic structure!marked} on $\Sigma\times \RR$,
we mean a pair $(N,f)$
where $N=\HH^3/\Gamma$ is an oriented complete hyperbolic $3$-manifold
and $f:\Sigma\times \RR \to N$ an orientation-preserving homeomorphism
which satisfies the following conditions.
\begin{itemize}
\item[$\circ$]
Let $\rho:\pi_1(\Sigma)\to \Gamma <\Isom^+\HH^3$ be the homomorphism
obtained as the composition of the homomorphism 
$(f\circ j)_*:\pi_1(\Sigma)\to \pi_1(N)$,
where $j:\Sigma\to \Sigma\times 0 \to \Sigma\times \RR$ is the inclusion map,
and the holonomy representation $\pi_1(N)\to \Isom^+(\HH^3)$ 
of the hyperbolic manifold $N$.
Then we require that $\rho$ is type-preserving.
(In other words, we require that the homeomorphism $f$ maps $(\mbox{ends of $\Sigma$})\times \RR$
into the {\it main} cusp of $N$ carrying the parabolic elements
$\rho(\mbox{peripheral elements})$.)
\end{itemize}
Thus we restrict our attention to the hyperbolic structures 
on the {\it pared manifold} 
$(\Sigma_{g,b}\times I, \partial \Sigma_{g,b}\times I)$
with $I=[-\infty,\infty]$ (see \cite[Section 7]{Thurston3}) for the teminology).

Two marked hyperbolic structures 
$(N_1,f_1)$ and $(N_2, f_2)$ on $\Sigma\times \RR$ are 
regarded as {\it equivalent}
if there is an orientation-preserving isometry $h:N_1\to N_2$
such that $h\circ f_1$ is homotopic to $f_2$.
This condition is equivalent to the condition
that the corresponding representations $\rho_1$ and $\rho_2$ are equal
up to conjugation by
an element of $\Isom^+\HH^3$.
Thus the set $H(\Sigma\times \RR)$
of all marked hyperbolic structures on $\Sigma\times \RR$ up to equivalence
is identified with the subset of the space
\[
\Rep_{\mathrm{tp}}(\Sigma\times \RR):=\{\rho:\pi_1(\Sigma)\to \Isom^+\HH^3 \ | \ \mbox{$\rho$ is type-preserving}\}/\Isom^+\HH^3
\]
consisting of (the images of) discrete faithful representations.
The set $H(\Sigma\times \RR)$ with the subspace topology is 
denoted by $AH(\Sigma\times \RR)$.
This topology is called the {\it algebraic topology}\label{algebraic topology}
 of $H(\Sigma\times \RR)$.
It is well-known that $\Rep_{\mathrm{tp}}(\Sigma\times \RR)$ is Hausdorff,
and $AH(\Sigma\times \RR)$ is a closed subset of $\Rep_{\mathrm{tp}}(\Sigma\times \RR)$ 
(cf. \cite[Section 4]{Marden1}).

Let $\QF(\Sigma\times \RR)$ be the subspace of $AH(\Sigma\times \RR)$
consisting of the quasifuchsian representations.
For each quasifuchsian representation $\rho:\pi_1(\Sigma)\to \PSL(2,\CC)$ with 
$\Gamma=\rho(\pi_1(\Sigma))$,
the Kleinian manifold $(\HH^3\cup \Omega(\Gamma))/\Gamma
\cong \Sigma\times [-\infty,\infty]$
is bounded by two marked Riemann surfaces $S_{\pm}=\Omega_{\pm}(\Gamma)/\Gamma$,
where $S_{\pm}$ correspond to 
$\Sigma\times\{\pm\infty\}\subset \Sigma\times [-\infty,\infty]$.
The pair $(S_-,S_+)$
is regarded as a point in the product $\Teich(\bar\Sigma)\times \Teich(\Sigma)$,
where $\bar\Sigma$ is the surface $\Sigma$ with the reverse orientation.
This determines a map
\[
\nu:\QF(\Sigma\times \RR) \to \Teich(\bar\Sigma)\times\Teich(\Sigma).
\]
Bers' simultaneous uniformization theorem\label{Bers' simultaneous uniformization theorem}
says that $\nu$ is a homeomorphism (see \cite{Imayoshi-Taniguchi}).

The positive solution to Thurston's Density Conjecture 
by Brock, Canary and Minsky \cite{BCM},
obtained as a consequence of deep results 
by a number of researchers 
in the deformation theory of Kleinian groups, 
says that $AH(\Sigma\times \RR)$ is equal to the closure of 
its open subset $\QF(\Sigma\times \RR)$:
\[
AH(\Sigma\times \RR)=\overline{\QF(\Sigma\times \RR)}
\] 
Thus any discrete faithful type-preserving $\PSL(2,\CC)$-representation of $\pi_1(\Sigma)$
is a limit of quasifuchsian representations.
In particular, a fiber Kleinian group of a hyperbolic surface bundle over $S^1$
is obtained as the limit of quasifuchsian groups.
Historically, the existence of the fiber Kleinian group 
(and so the existence of a complete hyperbolic structure on surface bundles)
was first proved in the case where $\Sigma$ is a once-punctured torus by J\o rgensen \cite{Jorgensen}:
the simplest case of the figure-eight knot complement
was also proved by Riley \cite{Riley1975}.
Thurston was impressed by these works.
He proved the hyperbolization theorem for surface bundles in \cite{Thurston5} 
(cf. Otal \cite{Otal2})
via his double limit theorem \cite[Theorem 4.1]{Thurston5}.

Cannon and Thurston \cite{Cannon-Thurston} found the following surprising fact.
Let $\rho_0:\pi_1(\Sigma)\to \PSL(2,\CC)$ be a Fuchsian representation, and
let $\rho:\pi_1(\Sigma)\to \PSL(2,\CC)$ be the type-preserving
discrete faithful representation whose image $\Gamma$ 
gives the fiber group of a hyperbolic $\Sigma$-bundle over $S^1$.
Recall that $\Lambda(\Gamma_0)=\partial \HH^2$ and 
$\Lambda(\Gamma)=\partial \HH^3$
(see Example \ref{example:limit-set4}), 
and $\pi_1(\Sigma)$
acts on these sets via $\rho_0$ and $\rho$, respectively.

\begin{theorem}[Cannon--Thurston map]
There is a $(\rho_0,\rho)$-equivariant surjective continuous map
\[
\CT:\partial\HH^2=\Lambda(\Gamma_0) \to \Lambda(\Gamma)=\partial \HH^3.
\]
\end{theorem}
The map $\CT$ is called the {\it Cannon--Thurston map}\index{Cannon--Thurston map}.
This theorem was first proved by Cannon and Thurston \cite{Cannon-Thurston} 
for the closed surface case,
and then proved by Bowditch \cite{Bowditch} for the general case.
Work of many authors has extended the results in various ways
(see the review \cite{Mj}).
For the simplest case where $\Sigma$ is the once-punctured torus,
the computer program {\it OPTi} developed by Wada \cite{Wada_OPTi}
visualizes deformations of the limit sets
of quasifuchsian punctured torus groups
(see \cite{ASWY} for background).
We can also see
a lot of breathtaking pictures related to the Cannon--Thurston maps
(mainly for the once-punctured torus)
in the book Indra's Pearls \cite{MSW}.

\subsection{Nielsen-Thurston classification of surface homeomorphisms and 
geometrization of surface bundles}
We quickly recall the Nielsen-Thurston classification of surface homeomorphisms
(see \cite{Thurston4, Fathi-Laudenbach-Poenaru, Farb-Margalit}).
Let $\MCG(\Sigma)$ be the 
{\it mapping class group}\index{mapping class group} of $\Sigma$
(the closed orientable surface of genus $g$ with $b$ punctures
such that $\chi(\Sigma)<0$),
the group of the orientation-preserving homeomorphisms of $\Sigma$
modulo isotopy.
We do not distinguish between a homeomorphism of $\Sigma$
and the element (mapping class) of $\MCG(\Sigma)$ represented by it,
as long as there is no fear of confusion.
Then Nielsen-Thurston theory says that for any 
$\varphi\in \MCG(\Sigma)$, one of the following holds.
\begin{enumerate}
\item
$\varphi$ is {\it periodic}\index{periodic}, namely
$\varphi$ has finite order in $\MCG(\Sigma)$.
In this case, $\varphi$ is represented by a (periodic) isometry 
with respect to some finite-volume complete hyperbolic structure on $\Sigma$.
\item
$\varphi$ is {\it reducible}\index{reducible},
i.e., there is a nonempty family of mutually disjoint essential simple loops 
whose union is preserved by (a representative of) $\varphi$. 
\item
$\varphi$ is {\it pseudo-Anosov}\index{pseudo-Anosov}.
This means that $\Sigma$ has a \lq \lq half-translation structure''
such that the homeomorphism $\varphi$ is \lq\lq realized by'' a diagonal matrix
$\begin{pmatrix}
\alpha & 0\\
0 & 1/\alpha
\end{pmatrix}$
with $\alpha>1$.
\end{enumerate}
The precise meaning of the last condition is as follows.
A {\it half-translation structure}\index{half-translation structure} 
on $\Sigma$ is
a singular Euclidean metric on $\Sigma$,
with a finite number of conical singularities of cone angle $k\pi$
($k\ge 3$),
and total cone angle $k'\pi$ ($k'\ge 1$) around each puncture.
The surface $\Sigma$ with cone points removed admits an isometric atlas over $\RR^2$
whose transition maps are of the form $(x,y)\mapsto \pm (x,y)+(a,b)$
for some $(a,b)\in\RR^2$.
Then $\varphi$ is pseudo-Anosov
if there is a half-translation structure on $\Sigma$,
such that the homeomorphism $\varphi$ has a local expression
$(x,y)\mapsto (\alpha x, \alpha^{-1}y)$
with respect to isometric atlas of the half-translation structure. 
The constant $\alpha$ is called the 
{\it expansion factor}\index{expansion factor}
 of the map $\varphi$.

This condition is described as follows in Thurston's original paper
\cite[Theorem 4]{Thurston4}:
there is a real number $\alpha>1$ and a pair of 
transverse measured foliations $\foliation^s$ and $\foliation^u$
such that $\varphi(\foliation^s)=\alpha^{-1}\foliation^s$
and $\varphi(\foliation^u)=\alpha\foliation^u$.
Here a {\it measured foliation}\index{measured foliation} 
on $\Sigma$ is a singular foliation
endowed with a measure in the transverse direction,
where only finitely many singularities of \lq\lq $k$-pronged saddle'' 
($k=1$ or $k\ge 3$) are allowed.
The notation $\foliation_1=\alpha\foliation_2$ means that 
$\foliation_1$ and $\foliation_2$ agree as foliations,
but the transverse measure of $\foliation_1$ is $\alpha$ times
that of $\foliation_2$.
With respect to the half-translation structure of $\Sigma$ discussed in the above,
the measured foliations $\foliation^s$ and $\foliation^u$
are the vertical and holizontal foliations, $\lambda^+$ and $\lambda^-$,
equipped with the transverse measures $|dx|$ and $|dy|$ respectively. 
(Note that every straight line segment in $\Sigma$ belongs to 
a unique (singular) foliation by parallel straight lines,
and so the vertical and horizontal foliations make sense.) 
Since $\varphi$ is locally expressed by 
$(x,y)\mapsto (\alpha x, \alpha^{-1}y)$,
it preserves the vertical and horizontal measured foliations
up to the factors $\alpha^{-1}$ and $\alpha$, respectively.

By considering the \lq\lq projective classes'' of measured foliations,
Thurston constructed the 
{\it projective measured foliation space}\index{projective measured foliation space}
$\PMF(\Sigma)$ and proved that it forms the boundary
of a natural compactification of the Teichm\"uller space $\Teich(\Sigma)$.
\[
\overline{\Teich}(\Sigma) =\Teich(\Sigma)\sqcup \PMF(\Sigma)
\cong \interior B^{6g-6+2b}\sqcup \partial B^{6g-6+2b}
\cong B^{6g-6+2b}
\]
The compactification is natural in the following sense. 
The action of 
$\MCG(\Sigma)$ on $\Teich(\Sigma)$
defined by the rule
\[
\varphi (S,f) :=(S, f\circ \varphi^{-1})
\quad \mbox{for $(S,f)\in \Teich(\Sigma)$}
\] 
extends to the action on the compactification,
so that its restriction to the boundary $\PML(\Sigma)$
is the natural action given by
\[
\int_{\gamma}\varphi_*(\foliation)=\int_{\varphi^{-1}(\gamma)}\foliation.
\]
Here $\gamma$ is an arc transverse to the foliation $\varphi(\foliation)$,
and $\int_{\gamma}\varphi_*(\foliation)$ is the measure of $\gamma$
with respect to the measured foliation $\varphi_*(\foliation)$.
It should be noted that 
the set of all essential simple loops in $\Sigma$ up to isotopy
is identified with a dense subset of $\PMF(\Sigma)$
and that the above action is an extension of the natural action of 
$\MCG(\Sigma)$ on $\mathcal{S}$.

By using this natural compactification of Teichm\"uller space,
Thurston established the classification of surface homeomorphisms,
as follows.
For a given $\varphi\in\MCG(\Sigma)$, its action on 
$\overline{\Teich}(\Sigma) \cong B^{6g-6+2b}$ 
has a fixed point, by Brower's fixed point theorem.
If there is a fixed point in $\Teich(\Sigma)$,
then $\varphi$ is periodic.
Suppose there is no fixed points in $\Teich(\Sigma)$
and so all fixed points lie in $\PMF(\Sigma)$.
If the underlying foliation of some fixed point 
contains a closed leaf,
then $\varphi$ is reducible.
Thurston managed to prove that $\varphi$ is pseudo-Anosov in the remaining case.

Now, let 
$M_{\varphi}:=\Sigma \times \RR/(x,t)\sim (\varphi(x), t+1)$
be the $\Sigma$-bundle over $S^1$ with monodromy $\varphi$.
Then it is easy to observe that
if $\varphi$ is periodic then $M_{\varphi}$ is a Seifert fibered space,
and that 
if $\varphi$ is reducible then $M_{\varphi}$ admits a nontrivial torus decomposition.
For the remaining case when $\varphi$ is pseudo-Anosov,
the following theorem was proved by Thurston,
as a special case of the geometrization Theorem \ref{thm:geometric-decomposition-mfd}.

\begin{theorem}
\label{thm:geometrization-bundle}
The surface bundle $M_{\varphi}$ is hyperbolic
if and only if $\varphi$ is pseudo-Anosov.
\end{theorem}

As noted in Subsection \ref{subsec:deformatin-space},
the corresponding fiber group $\rho\in AH(\Sigma\times I)$
is a limit of quasi-fuchsian groups.
Actually, for any $(S_-,S_+)\in \Teich(\bar\Sigma)\times\Teich(\Sigma)$,
$\rho$ is obtained as a limit
of a subsequence of 
the sequence of quasifuchsian groups
$\{\nu^{-1}(\varphi^{-k}(S_-), \varphi^{k}(S_+))\}_{k\ge 0}$
(see McMullen \cite[Theorem 3.8]{McMullen_1996}).

\subsection{Cannon--Thurston maps and veering triangulations}
\label{subsection:CT-Veering}
We now describe the combinatorial structure of the Cannon--Thurston map
associated with the $\Sigma$-bundle $M_{\varphi}$ 
with pseudo-Anosov monodromy $\varphi$.
Let $\rho_0:\pi_1(\Sigma)\to \PSL(2,\CC)$ be a Fuchsian representation
with image $\Gamma_0$, and
let $\rho:\pi_1(\Sigma)\to \PSL(2,\CC)$ be the 
discrete faithful representation whose image $\Gamma$ 
gives the fiber group of $M_{\varphi}$.

Let $j$ be the inclusion map
from $\Sigma=\HH^2/\Gamma_0$
to the infinite cyclic cover $\tilde M_{\varphi}=\HH^3/\Gamma$ of $M_{\varphi}$,
and consider its lift 
$\tilde j:\HH^2\to \HH^3$ to the universal cover.
Then the Cannon Thurston map $\CT:\partial\HH^2\to \partial\HH^3$
is the boundary map of the extension of $\tilde j$
to a map from $\HH^2\cup\partial\HH^2$ to $\HH^3\cup\partial\HH^3$.

In order to describe the combinatorics of the Cannon--Thurston map $\CT$,
let $\tilde \lambda^{\pm}$ be the singular foliations of $\HH^2$
obtained as the lifts of the vertical/horizontal foliations $\lambda^{\pm}$,
invariant by $\varphi$.
Then the endpoints of each leaf of $\tilde \lambda^{\pm}$ 
are mapped by $\CT$ into the same point,
and this turns out to generate the combinatorics of $\CT$.
To be precise, let $\sim^{\pm}$ be the equivalence relation on $\partial \HH^2$
which identifies the endpoints of each leaf of $\tilde \lambda^{\pm}$
by allowing for leaves that pass through singularities.
Let $\sim$ be the equivalence relation on $\partial \HH^2$ 
generated by $\sim^+$ and $\sim^-$.
Here the relations $\sim^+$ and $\sim^-$ 
are \lq\lq almost independent'' in the sense that
if $x\sim y$ then $x\sim^+ y$ or $x\sim^- y$ or else
there is a parabolic fixed point $p$ of $\Gamma_0$
such that either ($x\sim^+ p\sim^- y$) or ($x\sim^- p\sim^+ y$).
Moreover distinct parabolic fixed points of $\Gamma_0$ cannot be equivalent under $\sim$.
It was proved by Bowditch \cite[Section 9]{Bowditch} (cf. \cite[Section 5]{Cannon-Thurston})
that
\[
\mbox{ $\CT(x)=\CT(y)$ if and only if $x\sim y$.}
\]

\begin{figure}[ht]
\begin{center}
 {
\includegraphics[height=10cm]{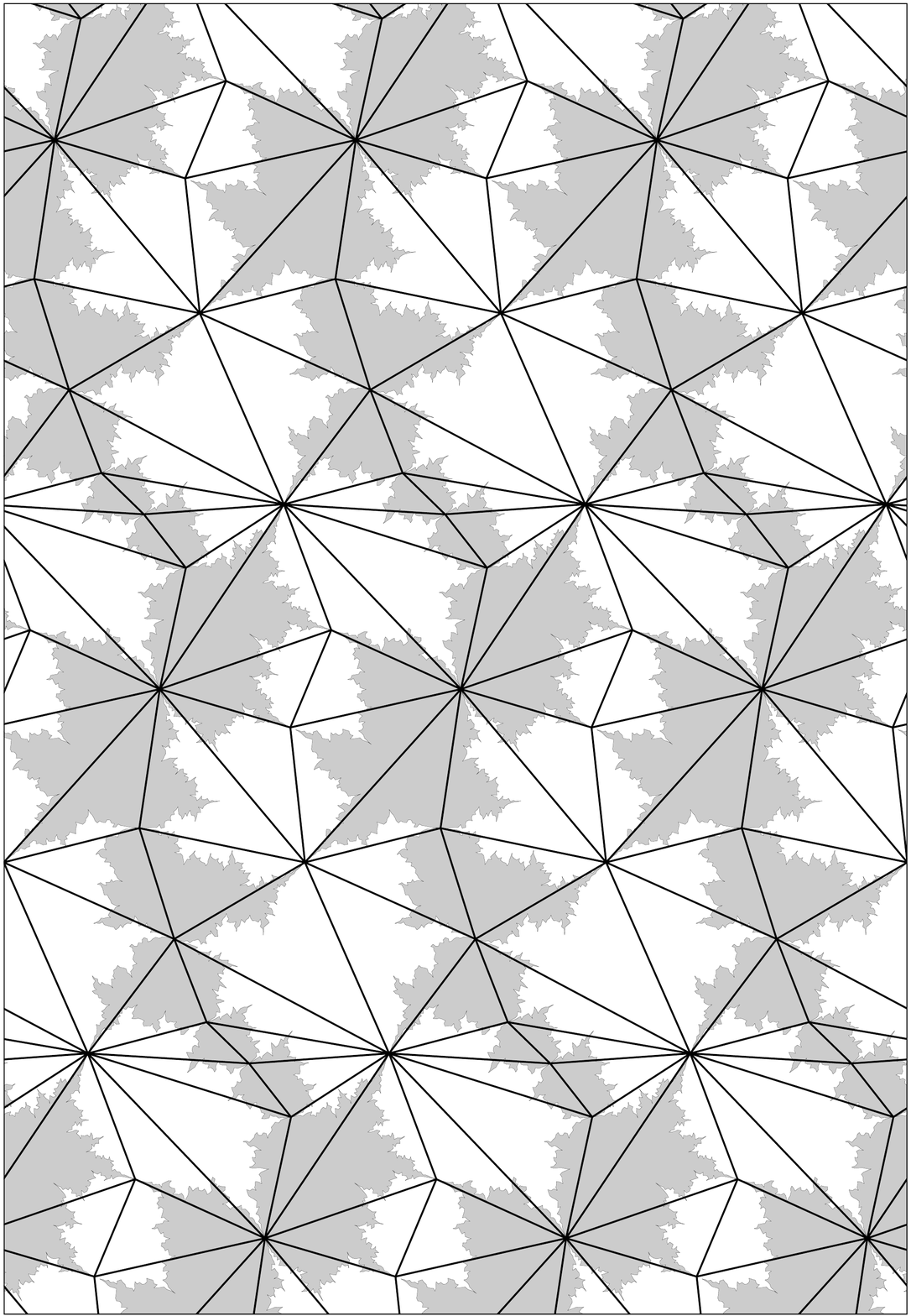}
 }
\end{center}
\caption{
Projected-horosphere
triangulation induced by the canonical decomposition and 
fractal tessellation for a once-punctured torus bundle.
Straight line segments etch the projected-horosphere triangulation while fractal arcs etch the fractal tessellation.
This picture is taken from \cite[Fig. 1]{Dicks-Sakuma}.
}
\label{fig.RLLRRRLLLL}
\end{figure}

In the remainder of this subsection,
we assume that
the singularities of the invariant foliations $\lambda^{\pm}$ occur 
only at punctures of the fiber.
(This condition is satisfied if $\Sigma$ is a once-punctured torus.)
Then it follows that for a point $q\in \partial\HH^3$,
the inverse image $\CT^{-1}(q)$ consists of 1, 2 or countably infinitely many points.
The last case happens if and only if $q$ is a parabolic fixed point of $\Gamma$,
and in this case $\partial\HH^2-\CT^{-1}(q)$ consists of
countably infinitely many open intervals.
Cannon and Dicks \cite{Cannon-Dicks2} studied the way
these intervals are mapped onto the complex plane $\CC\cong \partial\HH^3-\{q\}$,
and constructed a certain fractal tessellation of $\CC$ 
in the case where $\Sigma$ is a once-punctured torus.
Dicks and Sakuma \cite{Dicks-Sakuma} then observed that there is an intimate relation
between the fractal tessellation and the cusp triangulation 
(lifted to the universal cover $\CC$) induced by the canonical triangulation
of the hyperbolic once-punctured torus bundles
(see Figure \ref{fig.RLLRRRLLLL}).

On the other hand, Agol \cite{Agol2011} introduced
{\it veering triangulations}\index{veering triangulation},
which are (topological) ideal triangulations of
cusped hyperbolic $3$-manifolds
with a very special combinatorial structure.
He proved that every hyperbolic surface bundle,
for which the singularities of the invariant foliations $\lambda^{\pm}$ occur 
only at punctures of the fiber,
admits a veering triangulation,
which is canonical in the sense that 
it is uniquely determined by the fiber structure.
(More strongly, it is determined by Thurston's fiber face
to which the fibration belongs \cite{Minsky-Taylor}.)

In the beautiful paper \cite{Gueritaud2016},
Gu\'eritaud revealed an intimate relation
between the veering triangulation and 
the fractal tessellation arising from the Cannon--Thurston map
for every such hyperbolic surface bundle $M_{\varphi}$.
To this end, he gave a natural construction of 
the veering triangulation in terms of 
the invariant foliations.
The construction works in the universal cover $\tilde \Sigma$,
endowed with the half-translation structure
associated with the pseudo-Anosov monodromy.
He considered maximal rectangles in $\tilde \Sigma$ 
whose sides are vertical and horizontal in $\tilde \Sigma$ 
and whose interiors are disjoint from the singularities. 
Such maximal rectangles have one singularities on each side;
connecting these 4 singularities produces the ideal tetrahedra 
of the veering triangulation.
This construction enabled Gu\'eritaud to describe the relation
between the veering triangulation and the fractal tessellation
associated with the Cannon--Thurston map.

Roughly speaking, 
Gu\'eritaud's construction of the veering triangulation
is an analogue of the Delaunay triangulation 
relative to the singular set, 
with respect to the $\ell^{\infty}$-metric arising from 
the half-translation structure.
On the other hand, the canonical decomposition of a cusped hyperbolic manifold
is an analogue of the Delaunay triangulation 
relative to cusps,
with respect to the hyperbolic metric.
For hyperbolic once-punctured torus bundles,
these two decompositions are equal.
However, these two decompositions are quite different in general.
In fact, it was shown by Hodgson, Issa, Ahmad and Segerman \cite{HIS} that 
there exist veering triangulations which are not geometric,
in the sense that they are not isotopic to 
hyperbolic ideal triangulations.
Moreover, it was recently proved by 
Futer, Taylor and Worden \cite{FTW} that
generically veering triangulations are not geometric.
In spite of this defect from the viewpoint of hyperbolic geometry,
nice applications of veering triangulations to the study of curve complexes 
were given by Minsky and Taylor \cite{Minsky-Taylor}.

\section{Representations of $3$-manifold groups}
\label{sec:representation-space}

In Sections \ref{sec:Hyperbolic-Dehn-filling} and \ref{sec:flexibility},
we treated deformations of hyperbolic structures.
In Section \ref{sec:Hyperbolic-Dehn-filling},
we considered complete hyperbolic manifolds of finite volume
and studied deformations into incomplete hyperbolic structures,
whereas in Section \ref{sec:flexibility},
we considered complete hyperbolic manifolds of infinite volume
and studied deformations keeping the completeness.
In both sections, deformations are described in terms of deformations 
of holonomy representations.

One purpose of this section is to present the definition of 
{\it $\SL(2,\CC)$ character varieties},
which forms a common base ground for both treatments in 
Sections \ref{sec:Hyperbolic-Dehn-filling} and \ref{sec:flexibility},
and then to give a description of the hyperbolic Dehn filling theorem
independent of ideal triangulations, following 
Boileau-Heusener-Porti \cite[Appendix B]{Boileau-Porti}.
For another treatment, see 
Hodgson-Kerckhoff \cite[p.49, Remark]{Hodgson-Kerckhoff1998}.

Another purpose of this section is to describe applications
of the character varieties to knot theory and $3$-manifold theory.
We have already observed in Subsection \ref{subsec:universal-method}
that study of representations of knot groups to finite groups
gives us a powerful tool in knot theory.
The character variety, which is essentially
the space of representations of a knot group or a $3$-manifold group into the Lie group 
$\SL(2,\CC)$ up to 
conjugation by elements of $\SL(2,\CC)$,
leads to new versatile tools in knot theory and $3$-manifold theory.
We give a quick review to the Culler-Shalen theory 
\cite{Culler-Shalen, Culler-Shalen2, CGLS} and the $A$-polynomials
due to Cooper, Culler, Gillet, Long and Shalen \cite{CCGLS}.
For further information, 
see the survey Shalen \cite{Shalen}.

\subsection{Character variety}
\label{subsec:character-variety}
Let $M$ be a compact connected manifold, and 
let $R(M)=\Hom(\pi_1(M),\SL(2,\CC))$ be the space of all representations of $\pi_1(M)$ 
into $\SL(2,\CC)$.
This set has the structure of a complex affine algebraic set,
because it is identified with a subspace of
$(\SL(2,\CC))^k\subset \CC^{4k}$,
where $k$ is the cardinality of a generating set of $\pi_1(M)$,
defined by a system of polynomial equations.
For a representation $\rho\in R(M)$,
the function $\chi_{\rho}:\pi_1(M) \to \CC$ defined by 
$\chi_{\rho}(\gamma)=\tr(\rho(\gamma))$
is called the {\em character}\index{character} of $\rho$.
(We don't distinguish between a representation and the
element of $R(M)$.)
The set $X(M)$ of all characters also has the structure of an affine algebraic set,
and it is called the {\it character variety}\index{character variety} of $M$.
This can be seen as follows.
For each $\gamma\in\pi_1(M)$, consider the function
$I_{\gamma}:X(M)\to \CC$,
defined by $I_{\gamma}(\rho)=\chi(\gamma)$.
Then there are finitely many elements $\gamma_1,\dots,\gamma_d$
for which $I_{\gamma_1}\times \dots \times I_{\gamma_d}:X(M)\to \CC^d$
is an embedding,
and its image forms an affine algebraic set
\cite[Corollary 1.4.5]{Culler-Shalen}.

The natural projection from the space $R(M)/\SL(2,\CC)$
of all conjugacy classes of representations onto $X(M)$
fails to be injective
only at the conjugacy classes of reducible representations
(see \cite[Proposition 1.1.1]{Shalen}).
In this sense, $X(M)$ is regarded as the quotient $R(M)//\SL(2,\CC)$
in the category of affine algebraic sets.

If $M$ is a hyperbolic $3$-manifold, 
i.e., if $\interior M$ admits a complete hyperbolic structure,
then the holonomy representation $\rho:\pi_1(M)\to \PSL(2,\CC)$
lifts to an $\SL(2,\CC)$-representation (see \cite{Culler}).
In particular, the space $\Rep_{\mathrm{tp}}(\Sigma\times \RR)$
of conjugacy classes of type-preserving $\PSL(2,\CC)$-representations 
of $\pi_1(\Sigma\times \RR)$
(see Subsection \ref{subsec:deformatin-space})
is covered by a subspace of $X(\Sigma\times\RR)$. 

\subsection{Hyperbolic Dehn filling theorem and character variety}
\label{subsec:DehnFilling-CharacterVariety}

Consider the setting in Subsection \ref{subsec:HDF-theorem},
namely $M$ is a connected compact orientable $3$-manifold 
with $\partial M=\sqcup_{j=1}^m T_j$ a non-empty union of tori,
such that $\interior M$ admits a complete hyperbolic structure.
Let $\{\mu_j,\lambda_j\}$ is a pair of oriented slopes 
in the boundary torus $T_j$,
which forms a generator system of $H_1(T_j;\ZZ)$.
Let $\rho_0$ be a lift of the holonomy representation of
the complete hyperbolic structure of $\interior M$,
and let $\chi_0$ be its character.
Consider the map
$I_{\mu}=(I_{\mu_1},\dots, I_{\mu_m}):X(M)\to \CC^m$. 
Then the following theorem holds
(see \cite[Theorem B.1.2]{Boileau-Porti}).

\begin{theorem}
\label{thm:local-deformation-holonomy}
The map $I_{\mu}:X(M)\to \CC^m$ is locally bianalytic at $\chi_0$.
\end{theorem}

Using this theorem, we can associate each character
in some neighborhood of $\chi_0$ 
with generalized Dehn filling coefficients.
To describe this, recall that
the complex translation length $\LL_A$ of an element $A\in\SL(2,\CC)$
is defined as an element of $\CC/2\pi \sqrt{-1}\ZZ$
up to multiplication by $\pm 1$, by the following formula
(see Subsection\ref{subsec:hyperbolic-space}).
\[
\tr A=2\cosh \frac{\LL_A}{2}
\]
In order to have a well-defined complex translation length 
as an element in $\CC$, 
we consider the $(\ZZ/2\ZZ)^m$-branched covering map 
$\psi:\tilde U\to W$
from a neighborhood $\tilde U\subset \CC^m$ of the origin 
onto a neighborhood $W\subset X(M)$ of $\chi_0$ such that
\[
I_{\mu_j} \psi(\ub)=\epsilon_j \cosh\frac{u_j}{2} \quad 
\mbox{for every $\ub=(u_1,\dots, u_m)\in \tilde U$}
\]
where $\epsilon_j\in \{\pm\}$ is chosen so that 
$I_{\mu_j}(\psi(\mbox{\bf 0}))=\chi_0(\mu_j)=\tr(\rho_0(\mu_j))=\epsilon_j 2$.
Note that Theorem \ref{thm:local-deformation-holonomy}
guarantees the existence of this covering.

One can define a generalized Dehn surgery coefficients map as in 
Subsection \ref{subsec:HDF-theorem}, as
Theorem \ref{thm:hyp-Dehn2} holds in this setting,
where 
a certain open neighborhood, $U\subset \tilde U$, of the origin plays the role of
the open neighborhood $U$ in the theorem (see \cite[Proposition B.1.9]{Boileau-Porti}).
(In fact, by Remark \ref{rem.parameter-space}, 
the space $U$ is bi-holomorphic to
the space $U$ in Theorem \ref{thm:hyp-Dehn2},
when $M$ admits an ideal triangulation.)

In the setting of Subsection \ref{subsec:proof-HDFT},
each parameter $\ub\in U$ corresponds to a parameter $\zb$ 
representing the shapes of ideal tetrahedra,
and so it determines an (incomplete) hyperbolic structure on $\interior M$. 
In the current setting, we appeal to 
the fact that a small deformation of a hyperbolic structure
is parametrized by deformation of the holonomy representation 
(see \cite[Proposition 5.1]{Thurston1} and \cite[Proposition B.1.10]{Boileau-Porti}).
This is an outline of the proof 
the hyperbolic Dehn filling Theorem \ref{thm:hyp-Dehn}
without using an ideal triangulation of $M$,
given by \cite[Appendix B]{Boileau-Porti}.

\medskip

Note that the above proof is not effective in the sense that
it gives no information about the size or shape 
of hyperbolic Dehn surgery space $V\subset (\RR^2\cup\{\infty\})^m$.
In \cite{Hodgson-Kerckhoff2008} (cf. \cite{Hodgson-Kerckhoff1998, Hodgson-Kerckhoff2005}),
Hodgson and Kerckhoff developed a new theory 
of infinitesimal harmonic deformations for compact hyperbolic $3$-manifolds
with \lq\lq tubular neighborhood'',
and gave an effective proof of the hyperbolic Dehn filling theorem;
they proved that all generalized Dehn filling coefficients outside a disc
of \lq\lq uniform'' size yield hyperbolic structures.

\subsection{The Culler-Shalen theory and the cyclic surgery theorem}
\label{subsec:Culler-Shalen-theory}
We give a quick survey of the Culler-Shalen theory\index{Culler-Shalen theory}.
See \cite{Shalen} for a detailed self-contained review.
Let $M$ be a compact, connected, orientable, irreducible $3$-manifold.
Suppose $M$ contains an essential surface $F$.
Then by considering the inverse image $\tilde F$ of $F$
in the universal covering $\tilde M$,
we can construct a tree $T$,
such that the vertices correspond to the components of $\tilde M-\tilde F$
and the edges correspond to the connected components of $\tilde F$,
where the edge corresponding to a component of $\tilde F$
joins the two vertices corresponding to the components of $\tilde M-\tilde F$
abutting the component of $\tilde F$.
The covering transformation group $\pi_1(M)$ acts on the tree $T$ simplicially,
and this action is {\em nontrivial}
(i.e., no vertex is stabilized by the whole group $\pi_1(M)$)
and {\em without inversion}
(i.e., if an element $\gamma\in\pi_1(M)$ leaves an edge invariant,
then $\gamma$ fixes the edge pointwise).
Conversely, it is known that if $\pi_1(M)$ acts simplicially on a tree
nontrivially and without inversion,
then $M$ contains an essential surface.

In \cite{Culler-Shalen}, Culler and Shalen 
established a method for constructing 
such actions of $\pi_1(M)$ on trees,
by using the character variety $X(M)$.
The theory says that if $X(M)$ contains an algebraic curve $\mathcal{C}$,
then each {\em ideal point}\index{ideal point} of the curve gives rise to such an action of $\pi_1(M)$
and hence an essential surface in $M$.
In this theory, Tits-Bass-Serre theory \cite{Serre} on 
the structure of subgroups of $\SL(2,F)$, 
where $F$ is a field with a discrete valuation, plays a key role.
Various applications of this theory are given,
including (a) a simpler proof and generalization of the Smith conjecture \cite{Culler-Shalen}
and (b) a proof of the Neuwirth conjecture which says  
that every nontrivial knot group is a free product of two proper subgroups
amalgamated along a free product \cite{Culler-Shalen2}.

In \cite{CGLS}, Culler, Gordon, Luecke and Shalen
introduced a norm
$||\cdot||:H_1(\partial M;\RR)\to \RR$ 
for a compact orientable hyperbolic $3$-manifold $M$
with a single torus boundary.
A key fact behind this definition is the following:
Let $X_0$ be the irreducible component of the character variety $X(M)$
containing the character $\chi_0$ 
of the (lifted) holonomy representation of the complete hyperbolic structure on 
$\interior M$.
Then 
$X_0$ has complex dimension $1$ (see Theorem \ref{thm:local-deformation-holonomy}).
Let $\hat X_0$ be the projective completion of the affine algebraic curve $X_0$
in which the ideal points are smooth.
Then for each $\gamma\in H_1(\partial M)=\pi_1(\partial M)$,
the restriction to $X_0$ of the function $I_\gamma:X(M)\to \CC$, defined by 
$I_{\gamma}(\chi)=\chi(\gamma)$,
extends to a rational function, $\hat I_{\gamma}: \hat X_0\to \CC\cup\{\infty\}$,
where the ideal points of $\hat X_0$ (i.e., the points in $\hat X_0-X_0$) 
are the poles of this rational function.
The norm $||\cdot||:H_1(\partial M;\RR)\to \RR$ is defined to be the norm
obtained as the continuous extension of the function 
$H_1(\partial M;\ZZ)\to \ZZ$
which associates $\gamma$ with the degree of $\hat I_{\gamma}$.
The norm plays a crucial role in the proof of the cyclic surgery theorem below,
established by Culler, Gordon, Luecke, and Shalen \cite{CGLS}. 
The theorem was proved
by combining (i) arguments using the norm
and (ii) graph-theoretic analysis of 
the intersection of two incompressible, planar surfaces in knot exteriors.

\begin{theorem}[Cyclic surgery theorem]\index{cyclic surgery theorem}
\label{thm:cyclic-surgery}
Let $M$ be a compact, connected, orientable, irreducible $3$-manifold
such that $\partial M$ is a single torus, and
suppose that $M$ is not a Seifert fibered space. 
Let $\alpha$ and $\beta$ be two non-isotopic essential simple loops on $\partial M$,
such that $\pi_1(M(\alpha))$ and $\pi_1(M(\beta))$ are cyclic.
Then the geometric intersection number of $\alpha$ and $\beta$ is equal to $1$.
\end{theorem}

In \cite{Boyer-Zhang1, Boyer-Zhang2},
Boyer and Zhang generalized the above idea and
proved an analogue of the above theorem for finite surgeries.
See \cite{Boyer}, for further information.

The Culler-Shalen theory 
was extended by Morgan and Shalen \cite{Morgan-Shelen1984, Morgan-Shelen1988a, Morgan-Shelen1988b}
to the theory of $\RR$-trees.
Here an $\RR$-tree is a metric space in which any two points are joined 
by a unique topological arc.
The theory plays a key role in Otal's proof \cite{Otal2}
of the double limit theorem.
See the reviews \cite{Bestvina2002, Morgan1992} for further information.

\subsection{$A$-polynomials}
We give a short review of the $A$-polynomial of a knot $K$,
which is introduced by 
Cooper, Culler, Gillet, Long, and Shalen \cite{CCGLS}
by using the character variety $X(M)$ of the knot exterior $M$ of $K$. 
The idea is to consider the restriction map $r:X(M)\to X(\partial M)$
induced by the inclusion of $\pi_1(\partial M)$ into $\pi_1(M)$.
Then even though $X(M)$ is complicated,
its image $r(X(M))$ can be very simple.
Note that $\pi_1(\partial M)$ is the free abelian group
freely generated by 
the longitude $\lambda$ and the meridian $\mu$.
Thus, for any irreducible $1$-dimensional component $\mathcal{C}$ in the image 
$r(X(M))\subset X(\partial M)$,
there is a holomorphic map $f:\mathcal{C}\to \CC\times\CC$
which assigns the pair of the \lq\lq eigen values'' of 
the images of $\lambda$ and $\mu$ by the corresponding representations. 
Then the closure of the image $f(\mathcal{C})$ becomes an algebraic curve in $\CC^2$.
Such a curve is equal to the zero set of a single defining polynomial, $F_{\mathcal{C}}(x,y)$.
Now consider the product 
$\prod_{\mathcal{C}}F_{\mathcal{C}}(x,y)$
of the defining polynomials $F_{\mathcal{C}}(x,y)$
where $\mathcal{C}$ runs over the $1$-dimensional irreducible components
of $r(X(M))$.
Then the {\it $A$-polynomial}\index{$A$-polynomial} of $K$ is defined as
\[
A_K(x,y)=\frac{1}{x-1}\prod_{\mathcal{C}}F_{\mathcal{C}}(x,y)
\]
The reason of dividing out by the factor $x-1$ is that $H_1(M)$ is the free abelian group
generated by $\mu$
and so we always have a component corresponding to abelian representations,
which gives rise to the factor $x-1$.
By normalizing $A_K(x,y)$ so that it is an integral polynomial,
it is defined up to multiplication by $\pm x^ay^b$.

It is obvious that $A_O(x,y)=1$ for the trivial knot $O$,
and it is proved that the converse also holds 
(see Boyer-Zhang \cite{Boyer-Zhang3} and 
Dunfield-Garoufalidis \cite{Dunfield-Garoufalidis}).
The most important properties of the $A$-polynomials 
come from the fact that they encode information about the boundary slopes of the knot,
via the Newton polygon of $A_K(x,y)$.
Recall that a {\it boundary slope}\index{boundary slope} of a knot $K$ is 
a slope (isotopy class of an essential simple loop) in 
the boundary torus of the knot exterior $M$,
such that there is an essential surface in $M$
whose boundary consists of loops representing the slope. 
The Newton polygon of the polynomial $A_K(x,y)$ is the convex hull 
of the finite set:
\[
\{(i,j)\in\ZZ^2 \ | \ \mbox{the coefficient of $x^iy^j$ in $A_K(x,y)$ is non-zero}\}.
\]
The following striking theorem is proved by \cite[Theorem 3.4]{CCGLS}.

\begin{theorem}
Slopes of the edges of the Newton polygon of $A_K(x,y)$
are boundary slopes of the knot $K$.
\end{theorem}

\section{Knot genus and Thurston norm}
\label{sec:Thurston-norm}
By generalizing the genus of a knot,
Thurston \cite{Thurston_norm} defined a (semi-)norm 
on $H^1(M;\RR)\cong H_2(M,\partial M;\RR)$
for a compact orientable $3$-manifold $M$.
It is called the {\it Thurston (semi-)norm} of $M$.
By the work of Gabai \cite{Gabai1983a},
the Thurston norm is identical to the Gromov norm 
on $H_2(M,\partial M;\ZZ)$. 
The Thurston norm can be used to study the set of fiberings of $M$ over the circle,
and the work of Fried and McMullen enabled
a unified treatment of the fiberings of $M$.
After recalling these works,
we explain two {\it Thurstonian connections}
between the topology and geometry of $3$-manifolds, related to Thurston norms.
Namely, we survey (i) the relation of the Thurston norm
with the {\it hyperbolic torsion polynomial} 
due to Dunfield-Friedl-Jackson \cite{Dunfield-Friedl-Jackson} 
and Agol-Dunfield \cite{Agol-Dunfield},
and  (ii) that with the {\it harmonic $L^2$-norm} with respect to the hyperbolic metric
due to Brock-Dunfield \cite{Brock-Dunfield}.

\subsection{Thurston norm}
\label{subsec:Thurston-norm}
Let $M$ be a compact oriented $3$-manifolds
with $\partial M$ a possibly empty union of tori.
For a compact possibly disconnected surface $\Sigma$,
let $\Sigma_0$ be the surface consisting of the components of $\Sigma$  
which are neither homeomorphic to $D^2$ nor $S^2$,
and define its {\it complexity} 
by $\chi_-(\Sigma):=|\chi(\Sigma_0)|$.
For an integral homology class 
$\alpha\in H_2(M,\partial M;\ZZ)$,
define its {\it Thurston norm}\index{Thurston norm} $||\alpha||_{\mathrm{Th}}$ by
\[
||\alpha||_{\mathrm{Th}} = \min\{ \chi_-(\Sigma) \ | \ [\Sigma]=\alpha\}
\]

\begin{theorem}
{\rm (1)}
$||\cdot ||_{\mathrm{Th}}$ extends to a continuous map
$|| \cdot ||_{\mathrm{Th}}:H^1(M;\RR)\cong H_2(M,\partial M;\RR)
\to \RR_{\ge 0}$,
and this gives a semi-norm on $H^1(M;\RR)$.
Moreover, if any compact orientable surface properly embedded in $M$, 
representing a nontrivial homology class, 
has a negative Euler characteristic,
then 
$|| \cdot ||_{\mathrm{Th}}$ is a norm.

{\rm (2)}
Suppose $|| \cdot ||_{\mathrm{Th}}$ is a norm, then the unit ball
\[
B_M=\{\alpha\in H^1(M;\RR)\ | \ ||\alpha ||_{\mathrm{Th}}\le 1\}
\]
is a finite-sided polyhedron
whose vertices are rational points.

{\rm (3)}
Suppose $|| \cdot ||_{\mathrm{Th}}$ is a norm.
Then there are codimension one faces 
$F_1,\dots, F_k$, of $B_M$ satisfying the following conditions.
\begin{enumerate}[(i)]
\item
Any integral cohomology class in the interior of the cone $\RR_+\cdot F_i$
is a fiber class.
\item
Conversely, any fiber class is contained in the interior of 
some cone $\RR_+\cdot F_i$.
\end{enumerate}
\end{theorem}

Here a class $\phi\in H^1(M;\ZZ)$ is called a {\it fibered class}\index{fibered class}
if it is an integral multiple of the cohomology class
represented by a bundle projection $p:M\to S^1$. 
In the above theorem, each $F_i$ is called a {\it fibered face}\index{fibered face}.

The fiber structures contained in the interior of the cone
on a fibered face can be given a unified treatment,
and various interesting results can be obtained.
In particular, building on the results of Fried,
McMullen \cite{McMullen_Teich} proved that
each fibered face $F$ determines a $2$-dimensional \lq\lq lamination'' 
$\mathcal{L}$ of $M$
transverse to every fiber surface $\Sigma$  
with $(\mbox{Poincar\'e dual of $[\Sigma]$})\in \RR_+\cdot F$,
where $\Sigma\cap \mathcal{L}$ is the stable lamination 
of the monodromy of the fibration.
By using this result, he defined 
the {\it Teichm\"uller polynomial} \index{Teichm\"uller polynomial} 
$\theta_F\in \ZZ[H_1(M;\ZZ)/\Tor H_1(M;\ZZ)]$
and proved the following results \cite{McMullen_Teich}.
\begin{itemize}
\item[$\circ$]
The Teichm\"uller polynomial is symmetric,
i.e., if $\theta_F=\sum_g a_g g$ then
$\theta_F=\sum_g a_g g^{-1}$ up to a unit in $\ZZ[H_1(M)/{\Tor}H_1(M)]$.
\item[$\circ$]
For any integral cohomology class $\phi\in \RR_+\cdot F$,
the expansion factor $k(\varphi)$ of the corresponding monodromy $\varphi$ 
is equal to the largest root of the one-variable polynomial 
obtained by evaluating $\theta_F$ by $\phi$.
\item[$\circ$]
The function $\phi\mapsto 1/\log k(\varphi)$ 
extends to a real-analytic function on $\RR_+\cdot F$
which is strictly concave.
\item[$\circ$]
The cone $\RR_+\cdot F$ is dual to a vertex of the Newton polygon 
$\subset H_1(M;\RR)$ of $\theta_F$.
\item[$\circ$]
If the lamination $\mathcal{L}$ is transversely orientable,
then the (multivariable) Alexander polynomial of $M$ divides 
the Teichm\"uller polynomial $\theta_F$
\end{itemize}

\medskip
To end this subsection,
we recall an important result of Gabai \cite{Gabai1983a},
obtained as a corollary of his construction  
of codimension $1$ transversely oriented foliations without Reeb components 
which contain a given Thurston norm minimizing surface as a closed leaf.
To explain this, we consider another (semi-)norm
$||\cdot||_{\mathrm{Th}}^{\mathrm{s}}$ on $H^1(M;\RR)\cong H_2(M,\partial M;\RR)$
for a compact irreducible orientable $3$-manifold $M$,
defined by using immersed surfaces instead of embedded surfaces.
Namely,
for an integral homology class 
$\alpha\in H_2(M,\partial M;\ZZ)$,
define $||\alpha||_{\mathrm{Th}}^{\mathrm{s}}$ to be
the minimum of $\chi_-(\Sigma)$ of a compact oriented surface $\Sigma$ for which
there is a proper immersion $f:(\Sigma,\partial\Sigma)\to (M,\partial M)$
such that $f_*([\Sigma])=\alpha$, namely,
\[
||\alpha||_{\mathrm{Th}}^{\mathrm{s}} = 
\min\{ \chi_-(\Sigma) \ | \
\mbox{$\exists f:(\Sigma,\partial\Sigma)\looparrowright (M,\partial M)$
such that $f_*([\Sigma])=\alpha$}\}.
\]
The new norm $||\cdot||_{\mathrm{Th}}^{\mathrm{s}}$ is defined as a continuous extension
of the above norm on the integral homology.

In addition to this, 
as in Subsection \ref{subsec:Gromov-norm},
the {\it Gromov norm} $||\cdot||_{\mathrm{Gr}}$ is defined by
\[
||\alpha||_{\mathrm{Gr}}:=\inf\{||z|| \ | \ \mbox{$z$ is a singular cycle 
representing the homology class $\alpha$}\},
\]
where, for a (real) singular chain $z=\sum_j a_j\sigma_j$,
its norm $||z||$ is defined as the sum  $\sum_j |a_j|$ of the absolute values of its coefficients.
The following theorem was proved by Gabai \cite{Gabai1983a}.

\begin{theorem}
\label{thm:immersed-genus}
Let $M$ be a connected compact irreducible orientable $3$-manifold
with possibly empty toral boundary.
Then the three norms on $H^1(M;\RR)\cong H_2(M,\partial M;\RR)$
coincide, namely,
\[
||\cdot||_{\mathrm{Th}}=||\cdot||_{\mathrm{Th}}^{\mathrm{s}}=||\cdot||_{\mathrm{Gr}}.
\]
\end{theorem}

In particular, for a knot $K$ in $S^3$, its genus $g(K)$ is equal to
the {\it immersed genus} 
of $K$, 
which is defined as
the minimum of the genus $g(\Sigma)$ 
of a compact connected oriented surface $\Sigma$ such that
there is an immersion $f:\Sigma\to S^3$, 
with $f^{-1}(K) =\partial\Sigma$,
whose singular set is disjoint from $K$.
This is a generalization of Dehn's lemma for higher genus,
and gives a partial affirmative answer to a question raised by
Papakyriakopoulos \cite{Papakyriakopoulos},
who established Dehn's lemma. 

\subsection{Evaluation of Thurston norms in terms of Twisted Alexander polynomials}
\label{subsec:twisted-A-plynomial}
The {\it twisted Alexander polynomials}, 
defined by Lin \cite{Lin} for classical knots
and by Wada \cite{Wada} in the general setting,
give a powerful tool for studying the Thurston norm.
Such a \lq polynomial' $\Delta(M,\phi,\rho)$ depends on 
a class $\phi\in H^1(M;\ZZ)$ and a linear representation
$\rho:\pi_1(M)\to \GL(V)$, where $V$ is a finite dimensional vector space over a field $F$.
Then $\Delta(M,\phi,\rho)$ is  defined as an element of the quotient field
$F(t^{\pm 1})$ of the group ring $F[t^{\pm 1}]$,
and analogies of Theorem \ref{thm:Alexander-polynomial} on the classical Alexander polynomial are obtained by several authors 
(see the surveys \cite{Friedl-Viddusi1, Kitano}).
Friedl and Viddusi \cite{Friedl-Viddusi2, Friedl-Viddusi3} proved the surprising results that the twisted Alexander polynomials can detect fiber classes and the Thurston norms.

When $K$ is a hyperbolic knot in $S^3$,
it is natural to consider the twisted Alexander polynomial
for the representation $\rho:G(K)\to \SL(2,\CC)$ 
which projects to the holonomy representation of the complete hyperbolic structure 
of $S^3-K$.
Though there are precisely two such representations up to conjugacy,
there is unique one for which $\tr\rho(\mu)=+2$,
where $\mu$ a meridian of $K$.
(For the other lift $\rho'$, we have $\tr\rho'(\mu)=-2$.)
Thus we can consider the twisted Alexander ploynomial $\Delta(E(K),\phi,\rho)$,
where $\phi\in H^1(E(K);\ZZ)\cong \ZZ$ is the generator.
The invariant is called the 
{\it hyperbolic torsion polynomial}\index{hyperbolic torsion polynomial} of $K$
and is denoted by $\HypTorP_K(t)$ (see \cite{Dunfield-Friedl-Jackson}).
The artificial choice of the lift $\rho$ is irrelevant,
because if $\rho$ is replaced with the other lift $\rho'$,
then the corresponding polynomial $\HypTorP'_K(t)$ is equal to $\HypTorP_K(-t)$.
As a special case of the general results on the twisted Alexander polynomial,
the following hold for every hyperbolic knot $K$ in $S^3$.
\begin{enumerate}
\item
$4g(K)-2\ge \deg\HypTorP_K(t)$.
\item
If $K$ is fibered, then $\HypTorP_K(t)$ is monic.
\end{enumerate}
These may be regarded as analogies of Theorem \ref{thm:Alexander-polynomial}(2) and (3)
on the classical Alexander polynomial.
Dunfield, Friedl and Jackson \cite{Dunfield-Friedl-Jackson} 
made extensive computer experiments,
and confirmed that for all hyperbolic knots with at most 15 crossings,
the estimate (1) is sharp and that  (2) detect all non-fibered knots.
In particular, 
the hyperbolic torsion polynomial detects that 
the genera of the Kinoshita--Terasaka knot and the Conway knot are $3$ and $5$,
respectively.
(The genera of arborescent links, including these two knots, 
had been determined by
Gabai \cite{Gabai1986} through the topological study of complementary sutured manifolds.)
Thus the hyperbolic torsion polynomials can distinguish knots
which are mutants of each other.
In \cite{Agol-Dunfield}, Agol and Dunfield  
studied the conjecture posed by \cite{Dunfield-Friedl-Jackson}, 
that the estimate (1) is sharp for every hyperbolic knot,
and they verified the conjecture for {\it libroid} hyperbolic knots in $S^3$.
The libroid knots form a broad class of knots,
which is closed under Murasugi sum,
and in particular all arborescent are libroid knots.

\subsection{Harmonic norm and Thurston norm}
\label{sec:harmonic-norm}

Let $M$ be a closed orientable hyperbolic $3$-manifold.
Then in addition to the topologically defined Thurston norm $||\cdot||_{\mathrm{Th}}$,
there is yet another canonically defined geometric norm on $H^1(M;\RR)$.
By the rigidity theorem,
$M$ admits a unique hyperbolic metric,
and by applying Hodge theory to this Riemannian metric,
we can identify $H^1(M;\RR)$ with the space of harmonic $1$-forms.
Thus the {\it harmonic norm}\index{harmonic norm}
$||\cdot||_{L^2}$ determines another norm on 
$H^1(M;\RR)\cong H_2(M;\RR)$.
Here the harmonic norm is the one associated with the usual inner product for $1$-forms:
\[
\langle\alpha,\beta\rangle = \int_M \alpha\wedge *\beta,
\]
where $*$ is the Hodge $*$-operator.
Since it comes from a positive-definite inner product,
the unit ball of $||\cdot||_{L^2}$ is a smooth ellipsoid.
Brock and Dunfield \cite{Brock-Dunfield} proved the following relation
between the topological norm and the geometric norm.

\begin{theorem}
For all closed orientable hyperbolic $3$-manifold $M$ one has
\[
\frac{\pi}{\sqrt{\vol(M)}}||\cdot||_{\mathrm{Th}}
\le
||\cdot||_{L^2}
\le
\frac{10\pi}{\sqrt{\inj(M)}}||\cdot||_{\mathrm{Th}}.
\]
\end{theorem}

In the above theorem, $\inj(M)$ denotes the 
{\it injectivity radius}\index{injectivity radius} of $M$,
i.e., half of the length of the shortest closed geodesic.
Moreover, they also showed that the above estimates are
in some sense sharp, by giving families of examples.

These results were obtained as refinements of
a result of Bergeron, Seng\"un and Venkatesh \cite{Bergeron-Sengun-Venkatesh},
which in turn is preceded by the work by Kronheimer and Mrowka \cite{Kronheimer-Mrowka}
that characterize the Thurston norm as the infimum
(over all possible Rimennian metrics) of certain scaled harmonic metrics.

\section{Finite-index subgroups of knot groups and $3$-manifold groups}
\label{sec:subgroups}

As explained in Subsection \ref{subsec:universal-method},
finite branched/unbranched coverings of knots
are a powerful tool for distinguishing knots.
This fact reflects the richness of finite-index subgroups 
of knot groups. 
In this section, we survey the following topics
which illustrate this richness:
(i) universal groups, which produce all closed orientable $3$-manifolds,
(ii) positive solution of the virtual fibering conjecture,
(iii) Grothendieck rigidity of $3$-manifold groups, and
(iv) mysterious relation 
between the Gromov norm and the homology growth of finite coverings.

\subsection{Universal knots/links and universal groups}

In an unpublished preprint \cite{Thurston1882b}, 
W. Thurston presented a very complicated six component link in $S^3$,
and proved the surprising fact that
every closed orientable $3$-manifold can be expressed as a branched cover 
of the $3$-sphere branched over this link. 
He called links with this property {\it universal links}\index{universal link}.
He asked if a universal knot exists,
and if even the figure-eight knot was universal.
This question was answered affirmatively 
by Hilden, Lozano and Montesinos in
\cite{HLM1985},
where they proved that every hyperbolic $2$-bridge knot and link is universal.

\begin{figure}[ht]
\begin{center}
 {
\includegraphics[height=3.5cm]{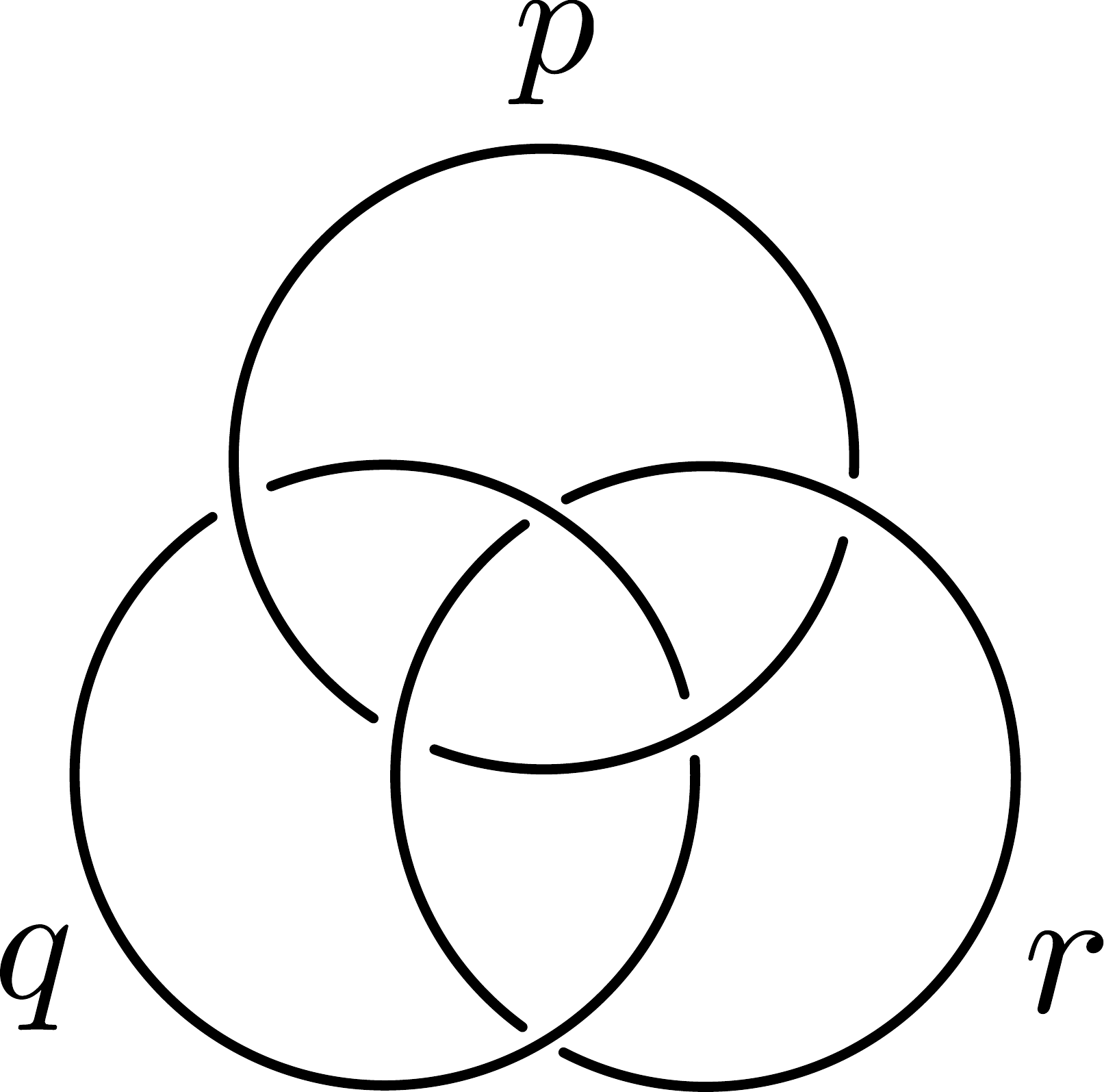}
 }
\end{center}
\caption
{
The {\it Borromean orbifold} $B(p,q,r)$
is a universal orbifold
if $p\ge 3$ and both $q$ and $r$ are even integers $\ge 4$
by \cite{HLM2010}.
}
\label{fig:Borromean}
\end{figure}

Moreover, it was later proved
by Hilden, Lozano, Montesinos and Whitten
in \cite{HLMW1987} that
every closed orientable $3$-manifold $M$ is 
a covering of $S^3$ branced over the Borromean rings
and having branching indices $1$, $2$ and $4$.
This implies that the hyperbolic orbifold 
$\UniversalOrb=\HH^3/\UniversalGp$
with underlying space $S^3$ and with singular set the Borromean ring
where all components have cone angle $\pi/2$,
is a {\it universal orbifold} in the following sense:
for any closed orientable $3$-manifold $M$,
there is a finite orbifold covering 
$\OO\to \UniversalOrb$
with underlying space $|\OO|$ homeomorphic to $M$.
In other words, the orbifold fundamental group
$\UniversalGp=\pi_1^{\mathrm{orb}}(\UniversalOrb)$
is a {\it universal group}\index{universal group}, i.e.,
for any closed orientable $3$-manifold $M$,
there is a finite index subgroup $\Gamma$ of $\UniversalGp$
such that $|\HH^3/\Gamma|\cong M$.
It is surprising that all closed orientable $3$-manifolds are
constructed from a single group $\UniversalGp$
and its finite index subgroups.
Moreover, universal groups seem to be ubiquitous (cf. \cite{HLM2010}).

\subsection{Virtual fibering conjecture}
The positive solution of
Thurston's virtual fibering conjecture
by Agol \cite{Agol-Groves-Manning} 
and the geometric solution of Waldhausen's conjecture for hyperbolic manifolds
due to Kahn-Markovich \cite{Kahn-Markovic},
which play a key role in the proof of the virtual fibering conjecture,
also reflect the richness of subgroups of Kleinian groups.
Please see Bestvina \cite{Bestvina2014}
for a survey of this important topic.

Here, I only recall
Walsh's simple construction \cite{Walsh}
of a nontrivial example of virtual fibering
using knot theory.
Let $K$ be a spherical Montesinos knot/link which is not fibered,
e.g., the $5_2$ knot, the $2$-bridge knot of slope $2/7$.
Then the double branched covering $M_2(K)$ of $S^3$
branched over $K$ is a spherical manifold and so
its universal covering $\tilde M_2(K)$ is the $3$-sphere.
The inverse image, $\tilde K$, of $K$ in the universal cover 
is a great circle link in $S^3$,
because it is the singular set of
the isometric group action of the 
$\pi$-orbifold group of $K$ (cf. Subsection \ref{subsec:Bonahon-Siebenmann-decomposition}).
Pick a component $O$ of $\tilde K$, and observe that
the remaining components form a closed braid around $O$,
because $\tilde K$ consists of great circles.
This shows that the covering $\tilde M_2(K)-\tilde K$ of $S^3-K$
is a punctured disk bundle over the circle,
though $S^3-K$ itself does not admit a fiber structure over the circle.

\subsection{Profinite completions of knot groups and $3$-manifold groups}
As explained in Subsection \ref{subsec:universal-method},
representations of knot groups onto finite groups
serve a powerful tool for distinguishing knots.
Thus it is natural to ask the following question
(cf. \cite{Boileau-Friedl2015}).

\begin{question}
\label{question:finite-quotients}
{\rm 
To what degree does the set of finite quotients of
knot groups distinguish knots?
More generally, what properties of $3$-manifolds are determined by
the set of finite quotients of their fundamental groups?
}
\end{question}

The geometrization theorem and Hempel's argument \cite{Hempel1987} 
show that every $3$-manifold group is 
{\it residually finite}\index{residually finite},
namely, for any nontrivial element $g \in \pi_1(M)$,
where $M$ is a compact connected orientable $3$-manifold,
there is a finite quotient of $\pi_1(M)$ in which $g$ remain nontrivial.
This implies that the above question can be formulated in terms of the
profinite completion of the fundamental group.

Recall that the {\it profinite completion}\index{profinite completion} 
of a group $\Gamma$,
is the inverse limit
of the inverse system of finite quotients of $\Gamma$:
we denote it by $\hat\Gamma$.
(The profinite completion is actually defined to be a topological group endowed with the {\it profinite topology}.
By Nikolov-Segal \cite{Nikolov-Segal},
the topology of any \lq\lq finitely generated profinite group'' 
is determined by the algebraic structure.
So we do not care about the topological structure in this subsection.) 
The natural map $\Gamma\to \hat\Gamma$ is injective
if and only if $\Gamma$ is residually finite.
Let $\Subgp(\Gamma)$ denote the family of finite quotients
of $\Gamma$.
Then the following holds (see \cite[p.88-89]{Ribes-Zaresskii}, 
\cite[Theorem 2.2]{Long-Reid2011}).

\begin{theorem}
For two finitely generated residually finite groups $\Gamma_1$ and $\Gamma_2$,
the equality $\Subgp(\Gamma_1)=\Subgp(\Gamma_2)$ holds
if and only if 
$\hat\Gamma_1\cong \hat\Gamma_2$,
i.e., the profinite completions 
are isomorphic.
\end{theorem}

Thus Question \ref{question:finite-quotients} is 
reformulated by using the profinite completion;
in particular, the following question arises as a special case.

\begin{question}
\label{question:profinite}
{\rm
Let $M_1$ and $M_2$ be connected compact orientable $3$-manifolds,
for which the profinite completions $\widehat{\pi_1(M_1)}$ and
$\widehat{\pi_1(M_2)}$ are isomorphic.
Are $\pi_1(M_1)$ and $\pi_1(M_2)$ isomorphic?
}
\end{question}

The answer to the above question is no.
In fact, Funar \cite{Funar2013} 
and Hempel \cite{Hempel-ax1409} showed 
that the profinite completion of the fundamental group cannot always distinguish 
certain pairs of torus bundles nor certain pairs of Seifert fibered spaces. 
It is still an open question though whether the profinite completion can distinguish any two hyperbolic $3$-manifolds.
Boileau and Friedl \cite{Boileau-Friedl2015} considered a more relaxed
Question \ref{question:finite-quotients}
and obtained various interesting results concerning fiberedness and the Thurston norm,
and have shown that the figure-eight knot and torus knots
are distinguished from other knots by the profinite completions of their knot groups.

On the other hand, the following problem had been posed by 
Grothendieck \cite{Grothendieck}.

\begin{problem}[Grothendieck]
\label{prob:Grothendieck}
{\rm
Let $\varphi:\Gamma_1\to\Gamma_2$ be a homomorphism
of finitely presented residually finite groups
for which the extension $\hat{\varphi}:\hat\Gamma_1\to \hat\Gamma_2$ is 
an isomorphism.
Is $\varphi$ an isomorphism?
}
\end{problem}

If $\hat{\varphi}:\hat\Gamma_1\to \hat\Gamma_2$ is an isomorphism,
then the composition $\Gamma_1\to \hat\Gamma_1\to \hat\Gamma_2$ is an injection
and so $\varphi:\Gamma_1\to\Gamma_2$ must be an injection.
Therefore Grothendieck's problem reduces to the case
where $\Gamma_1$ is a subgroup of $\Gamma_2$.
Long and Reid \cite{Long-Reid2011}
introduced the following terminology.
For a group $G$ and its subgroup $H<G$, the pair
$(G,H)$ is a {\it Grothendieck pair}\index{Grothendieck pair}
if the inclusion $j:H\to G$ provides a negative answer to
Grothendieck's problem.
If for all finitely generated subgroups $H<G$,
$(G,H)$ is never a Grothendieck pair
then $G$ is {\it Grothendieck rigid}\index{Grothendieck rigid}.

The following theorem was proved 
by Long and Reid \cite{Long-Reid2011}
for the case where $M$ is closed
and by Boileau and Friedl \cite{Boileau-Friedl2017}
for general case.

\begin{theorem}
\label{thm:Grothendieck-rigid}
Let $M$ be a connected, orientable, irreducible, compact $3$-manifold.
Then $\pi_1(M)$ is  Grothendieck rigid.
\end{theorem}

In the examples of
Funar \cite{Funar2013} 
and Hempel \cite{Hempel-ax1409},
the isomorphisms between the profinite completions
are not induced by a homomorphism between the $3$-manifold groups.

\subsection{Homology growth}
\label{subsec:homology-growth}
Investigation of
the first homology groups of finite (branched or unbranched) coverings
has a long history (cf. Subsection \ref{subsec:universal-method}).
For the homology of finite abelian coverings of links,
it was proved that  
they are essentially determined by the Alexander invariants of links
(see \cite{Fox1954, Mayberry-Murasugi1982, Sakuma1995}). 
In \cite{Gordon1972},
Gordon studied the asymptotic behavior the homology of 
finite cyclic branched coverings of a knot, and
gave a necessary and sufficient condition for $H_1(M_n(K);\ZZ)$
to be periodic with respect to $n$,
in terms of the Alexander invariants.
This in particular implies that if the Alexander polynomial $\Delta_K(t)$
has a root which is not a primitive root of $1$
then  $H_1(M_n(K);\ZZ)$ cannot be periodic.
In fact, he showed that the order $|H_1(M_n(K);\ZZ)|$ is unbounded
under the same assumption, and then
asked if the order $|H_1(M_n(K);\ZZ)|$ tends to $\infty$.
Riley \cite{Riley1990} and 
Gonz\'alez-Acu\~na and Short \cite{GS1991}, independently,
proved that $|H_1(M_n(K);\ZZ)|$ grows exponentially.
To be precise, the following was proved:
\[
\lim_{n_j\to\infty} \frac{1}{n_j}\log|H_1(M_{n_j}(K);\ZZ)|
=
\log\mathbb{M}(\Delta_K),
\]
where $\{n_j\}$ runs over the natural numbers
such that $|H_1(M_{n_j}(K);\ZZ)|$ is finite,
and $\mathbb{M}(\Delta_K)$ is the {\it Mahler measure}\index{Mahler measure} of 
the Alexander polynomial $\Delta_K(t)$.
The Mahler measure of a polynomial $f(t)$ is defined by
\begin{align*}
\mathbb{M}(f)
&=
\exp
\left(
\int_{S^1}
\log |f(s)|ds
\right)\\
&=
\exp
\left(
\int_0^1 
\log |f(e^{2\pi\sqrt{-1}t})|dt\right)\\
&=
|c|\prod_{f(\omega)=0}\max(|\omega|,1)
\quad \mbox{($c$ is the leading coefficient of $f(t)$).}
\end{align*}
Here, the last equality is a consequence of Jensen's formula 
\cite[p. 205]{Ahlfors1966}. 

This result was extended by Silver and Williams \cite{SlWl2002} to
links in $S^3$,
by using the result of Schmidt \cite{Scm1995}
on the entropy of a certain dynamical system.
Let $L$ be an $m$-component oriented link in $S^3$, 
with the complement $X = S^3-L$. 
For a subgroup $\Lambda\subset H_1(X;\ZZ)\cong \ZZ^m$ of rank $m$,
let $X_{\Lambda}^{\rm{br}}$ be the corresponding branched covering of $Z$.
Set 
\[
\langle \Lambda \rangle=\min\{|x| \ | \ x\in\Lambda-\{0\}\},
\] 
where $|x| = \sqrt{\sum_i |x_i|^2}$ for $x=(x_1,\dots,x_m)\in\ZZ^m$.
Let $\Delta_L\in\ZZ[t_1,\dots,t_m]$ be the ($0$-th) Alexander polynomial of $L$.

\begin{theorem}[Silver-Williams]
\label{thm:homology-abelian-cover}
Under the above setting, suppose that $\Delta_L\ne 0$.
Then the following holds:
\[
\limsup_{\langle \Lambda \rangle\to \infty}
\frac{\ln |\Tor_{\ZZ} H_1(X_{\Lambda}^{\rm{br}};\ZZ)|}{|\ZZ^m/\Lambda|}
=
\log\mathbb{M}(\Delta_L).
\]
\end{theorem}
\noindent
Here $\mathbb{M}(\Delta_L)$ is the Mahler measure of $\Delta_L$,
defined by
\[
\mathbb{M}(\Delta_L)=
\exp
\left(
\int_{T^m}
\log |\Delta_L(s)|ds
\right),
\]
where $T^m:=(S^1)^m\subset \CC^m$ is the multiplicative subgroup in $\CC^m$,
and $ds$ indicates integration with respect to normalized Haar measure on $T^m$.

In \cite{Le2014}, Thang Le solved a conjecture of Schmidt \cite{Scm1995},
and by using the solution,
he extended the above result 
to links in oriented integral homology $3$-spheres,
and to include the case where the $0$-th Alexander polynomial vanishes,
by replacing the $0$-th Alexander polynomial with the first non-vanishing
Alexander polynomial.
He also proved that the same formula holds for unbranched abelian coverings.

In \cite{Le2018}, Le also studied the asymptotic behavior of 
the homology of non-abelian coverings,
by using the result on $L^2$-torsion by L\"uck 
\cite[Theorems 4.3 and 4.9]{Luck2002}.
Let $X$ be an irreducible compact orientable $3$-manifold
with infinite fundamental group
with (possibly empty) toral boundary.
For a subgroup $\Gamma$ of $\pi_1(X)$ of finite index,
let $X_{\Gamma}$ be the corresponding finite covering of $X$.
A sequence $\{\Gamma_k\}$ 
of subgroups of $\pi_1(X)$ of finite index
is said to be {\it nested},
if $\Gamma_{k+1}<\Gamma_k$.
It is said to be {\it exhaustive} if $\cap_k \Gamma_k=\{1\}$.

\begin{theorem}[Le]
\label{thm:Le2}
Under the above setting, the following holds for 
any nested exhaustive sequence $\{\Gamma_k\}$ 
of normal subgroups of $\pi_1(X)$ of finite index:
\[
\limsup_{k\to\infty}\frac{\ln |\Tor H_1(X_{\Gamma_k};\ZZ)|}{[\pi_1(X):\Gamma_k]}
\le \frac{V_{\mathrm{tet}}||X||}{6\pi},
\]
where $||X||$ is the Gromov norm of $X$.
\end{theorem}

For a knot $K$ in $S^3$ with exterior $X$
and a finite index subroup $\Gamma<G(K)$,
let $X_{\Gamma}^{\mathrm{br}}$ be the corresponding branched covering of $S^3$ 
branched over $K$.

\begin{theorem}[Le]
\label{thm:Le3}
Under the above setting, 
the following holds
for any nested exhaustive sequence $\{\Gamma_k\}$ 
of normal subgroups of $G(K)=\pi_1(X)$ of finite index.
\[
\limsup_{k\to\infty}\frac{\ln |\Tor H_1(X_{\Gamma_k}^{\mathrm{br}};\ZZ)|}{[\pi_1(X):\Gamma_k]}
\le \frac{V_{\mathrm{tet}}||X||}{6\pi},
\]
where $||X||$ is the Gromov norm of $X$.
\end{theorem}

For the sake of simplicity,
we stated Le's theorem only for regular coverings.
However, the actual statement of his theorem is much more general
and it does not 
restrict to regular coverings.
For a precise statement, see the original paper \cite{Le2018}.
Moreover, he conjectures that the identity holds in both theorems.

The homology of finite (branched/unbranched)
coverings is a common and well-known invariant in knot theory.
It is impressive that the asymptotic behavior of this familiar invariant
reflects the deep geometric structure of the knot.

\bibstyle{plain}

\end{document}